\documentclass[12pt,leqno,fleqn]{amsart}
\usepackage{amsmath,amstext,amsthm,amssymb,amsxtra}
\usepackage{txfonts} %also pxfonts
\usepackage[colorlinks,citecolor=red,pagebackref,hypertexnames=false]{hyperref}
\usepackage{pgf,tikz}
\usepackage{pdfsync}
\allowdisplaybreaks
\newcommand{\RR}{{\mathbb{R}}}
\newcommand{\NN}{{\mathbb{N}}}
\newcommand{\ZZ}{{\mathbb{Z}}}

\newcommand{\Lmn}{{\mathcal L}^{m,n}}
\newcommand{\Lmno}{{\mathcal L}^{m,n}_o}

\newcommand{\eps}{\varepsilon}
\newcommand{\po}{\partial\Omega}

\newcommand{\bp}{\noindent {\it Proof}.\,\,}
\newcommand{\ep}{\hfill$\Box$ \vskip 0.08in}

\def\ring{\mathaccent"0017 }

\newcommand{\todo}[1]{\vspace{5 mm}\par \noindent
\marginpar{\textsc{}} \framebox{\begin{minipage}[c]{0.95
\textwidth} \tt #1
\end{minipage}}\vspace{5 mm}\par}

\newtheorem{proposition}{Proposition}[section]
\newtheorem{theorem}[proposition]{Theorem}
\newtheorem{lemma}[proposition]{Lemma}
\newtheorem{corollary}[proposition]{Corollary}

\usepackage{color}

\newcommand\comment[1]{{\color{purple}#1}}
\long\def\comment#1{}\long\def\notcomment#1{}

\begin{document}

\author{Svitlana Mayboroda and  Vladimir Maz'ya}

\title[Regularity of solutions to the polyharmonic equation]{Regularity of solutions to the polyharmonic equation in general domains}

%Boundedness and continuity at the boundary of the gradient
%for the solutions of the biharmonic equation

%\date{ }

\begin{abstract} 
%The maximum principle would normally warrant that for any elliptic operator of order $2m$ any solution with nice data has bounded $(m-1)$-st derivatives. It is known though that for higher order equations ($m>2$) in non-smooth domains the maximum principle may fail  in dimensions $n\geq 4$.

The present paper establishes boundedness of $\left[m-\frac n2+\frac 12\right]$ derivatives for the solutions to the polyharmonic equation of order $2m$ in arbitrary bounded open sets of $\RR^n$, $2\leq n\leq 2m+1$, without any restrictions on the geometry of the underlying domain. It is shown that this result is sharp and cannot be improved in general domains. Moreover, it is accompanied by sharp estimates on the polyharmonic Green function. 

\end{abstract}

\maketitle

\tableofcontents

\section{Introduction}\label{s1}
\setcounter{equation}{0}

%Higher order elliptic boundary value problems, despite their  tremendous importance in physics and     engineering \cite{Meleshko}, are still rather poorly understood. Numerical methods geared towards  higher order PDEs  are  demanding and often target only a few specific  situations  \cite{CiarletMethods}, \cite{Dong}. At the same time, the world of higher order operators on general non-smooth domains remains largely a terra incognita. 

Higher order differential operators are important in physics and in engineering \cite{Meleshko} and  have been integrated in many areas of mathematics, including conformal geometry (Paneitz operator, $Q$-curvature \cite{AliceChang}, \cite{AliceChang2}), free boundary problems \cite{AdamsPotAnal}, non-linear elasticity \cite{Skrypnik},  \cite{CiarletElasticity123}, \cite{Antman2005}, to mention just a few. Unfortunately,  in spite of  evident demand, the properties of higher order PDEs on {\it general} domains remained largely beyond  reach. Their investigation brought challenging hypotheses and surprising counterexamples, and few general positive results. For instance, Hadamard's 1908 conjecture regarding positivity of the biharmonic Green function \cite{Hadamard} was actually refuted in 1949 (see \cite{Duffin}, \cite{Garabedian}, \cite{Shapiro}), and later on the weak maximum principle was proved to fail as well, at least in high dimensions \cite{MR}, \cite{PVpoly}. Another curious feature is a paradox of passage to the limit for solutions under approximation of a smooth domain by polygons \cite{Bab}, \cite{MN86}. 

For a long time, almost all results for higher order PDEs pertained to piecewise smooth domains (see, e.g., \cite{KMR0}, \cite{KMR}, \cite{MRbook}). The recent two decades have witnessed a great burst of activity in the study of boundary value problems on Lipschitz domains, which brought  some powerful and delicate estimates for the bi-Laplacian and for the polyharmonic operator due to B.\,Dahlberg, C.\,Kenig, J.\,Pipher, G.\,Verchota, Z.\,Shen, M.\,Mitrea, and others (\cite{DKV}, \cite{Ver90}, \cite{PV91}, \cite{PVLp}, \cite{PVmax}, \cite{KenigBook}, \cite{PVpoly},  \cite{PVdilation}, \cite{Ver96}, \cite{DKPV97}, \cite{AdPi}, \cite{Ver05}, \cite{She06A},  \cite{She06B}, \cite{She06C}, \cite{She07A}, \cite{Mit10}, \cite{MazMS10}, \cite{Ver10}, \cite{KS11A}, \cite{MitMW11}). Unfortunately, none of these results could be extended to domains of more complicated geometry. 

{\it The present paper establishes sharp pointwise estimates on the solutions to the polyharmonic equation and their derivatives in arbitrary bounded open sets, without any restrictions on the geometry of the underlying domain. It is shown that  these estimates are sharp and can not be improved. }

Let us recall the second order case. 
One of the fundamental results of elliptic theory is the maximum principle for harmonic functions. It holds in {\it
arbitrary} domains and guarantees that every solution to the
Dirichlet problem for the Laplace's equation, with bounded data,
is bounded. A similar statement is valid for a broad class of second order
 elliptic operators.

In 1960 an appropriate, weak, version of the maximum principle was obtained for  higher order equations
on {\it smooth} domains (\cite{Agmon}, see also \cite{Mir48}, \cite{Mir58}). Roughly speaking, it
amounts to the estimate
\begin{equation}\label{eq1.1}
\|\nabla^{m-1}u\|_{L^\infty(\overline\Omega)}\leq
C\|\nabla^{m-1}u\|_{L^\infty(\po)},
\end{equation}
\noindent where $u$ is a solution of an elliptic differential
equation of order $2m$ such that $|\nabla^{m-1}u|$ is continuous
up to the boundary, $\nabla^{m-1}u=\{\partial^\alpha
u\}_{|\alpha|=m-1}$ is a vector of all partial derivatives of $u$
of order $m-1$ and we adopt the usual convention that the zeroth order derivative of $u$ is $u$ itself. The formulation in non-smooth domains is somewhat trickier (see, e.g., \cite{Necas}). However, in any setting, the weak maximum principle would always guarantee that the solution with ``nice" data has bounded derivatives of order $m-1$.  In the early 1990s, \eqref{eq1.1} was extended to
three-dimensional domains diffeomorphic to a polyhedron (\cite{KMR}, \cite{MR1}) or having a Lipschitz boundary (\cite{PVmax}, \cite{PVpoly}). In general domains, no direct analog of the maximum principle exists  (see Problem~4.3, p.~275, in
Ne\v cas's book \cite{Necas}).

Moreover, it turns out that   for every elliptic operator
of order {\it greater than two} the maximum principle can be violated, in a four-dimensional cone  (\cite{MR}, see also \cite{MazNP83},  \cite{PVLp}, \cite{PVpoly}). It has been shown, in particular, that in dimensions $n\geq 4$ there are solutions to the polyharmonic equation with unbounded $(m-1)$-st derivatives (cf. \eqref{eq1.1}). This phenomenon raises two fundamental questions: whether the boundedness of the $(m-1)$-st derivatives remains valid in dimensions $n\leq 3$, and whether there are some other, possibly lower-order, estimates that characterize the solutions when $n\geq 4$. The main result of the present paper is as follows. 

\begin{theorem}\label{t1.1}
Let $\Omega$ be a bounded domain in $\RR^n$, $2\leq n \leq 2m+1$,
and
\begin{equation}\label{eq1.2}
(-\Delta)^m u=f \,\,{\mbox{in}}\,\,\Omega, \quad f\in
C_0^{\infty}(\Omega),\quad u\in \ring
W^{m,2}(\Omega).
\end{equation}

\noindent Then  the
solution to the boundary value problem \eqref{eq1.2}
satisfies
\begin{equation}\label{eq1.3}
\nabla^{m-n/2+1/2} u\in L^\infty(\Omega)\,\,\mbox{when $n$ is odd\quad and \quad} \nabla^{m-n/2} u\in L^\infty(\Omega)\,\,\mbox{when $n$ is even}.
\end{equation}

\noindent In particular, 
\begin{equation}\label{eq1.4}
\nabla^{m-1} u\in L^\infty(\Omega)\,\,\mbox{when $n=2,3$}.
\end{equation}

\end{theorem}

\noindent Here the space $\ring W^{m,2}(\Omega)$, is, as usually, a
completion of $C_0^\infty(\Omega)$ in the norm given by $\|u\|_{\ring
W^{m,2}(\Omega)}=\|\nabla^m u\|_{L^2(\Omega)}$. We note that $\ring W^{m,2}(\Omega)$ embeds into $C^k(\Omega)$ only when $k$ is strictly smaller than $m-\frac n2$, $n<2m$. Thus, whether the dimension is even or odd, Theorem~\ref{t1.1} gains one derivative over the outcome of Sobolev embedding.

The results of Theorem~\ref{t1.1} are sharp, in the sense that the solutions  do not exhibit higher smoothness than warranted by \eqref{eq1.3}--\eqref{eq1.4} in general domains. Indeed, assume that $n\in [3,2m+1]\cap\NN$ is odd and let $\Omega\subset{\mathbb{R}}^n$ be the punctured unit ball
$B_1\setminus\{O\}$, where $B_r=\{x\in{\mathbb{R}}^n:\,|x|<r\}$. Consider
a function $\eta\in C_0^\infty(B_{1/2})$ such that $\eta=1$ on
$B_{1/4}$. Then let
\begin{equation}\label{int6}
u(x):=\eta(x)\,\partial_x^{m-\frac n2-\frac 12}(|x|^{2m-n}),\qquad x\in B_1\setminus\{O\},
\end{equation}

\noindent where $\partial_x$ stands for a derivative in the direction of $x_i$ for some $i=1,...,n$.
It is straightforward to check that $u\in \ring W^{m,2}(\Omega)$ and $(-\Delta)^m u\in
C_0^\infty (\Omega)$. While $\nabla^{m-\frac n2+\frac 12} u$ is bounded, the derivatives of the order $m-\frac n2+\frac 32$ are not, and moreover, $\nabla^{m-\frac n2+\frac 12} u$
is not 
continuous at the origin. Therefore, the estimates \eqref{eq1.3}--\eqref{eq1.4} are optimal in general domains. 

As for the case when $n$ is even, the results in \cite[Section~10.4]{KMR} demonstrate  that in an exterior of a ray there is an $m$-harmonic function behaving as $|x|^{m-\frac n2+\frac 12}$. Thus, upon truncation by the aforementioned cut-off $\eta$, one obtains a solution to \eqref{eq1.2} in $B_1\setminus\{x_1=0, ..., x_{n-1}=0, 0\leq x_n<1\}$, whose derivatives of order $m-\frac n2+1$ are not bounded. Therefore, in even dimensions \eqref{eq1.3} is a sharp property as well.

Furthermore, Theorem~\ref{t1.1} has several quantitative manifestations, providing specific estimates on the solutions to \eqref{eq1.2}. Most importantly, {\it we establish sharp pointwise estimates on Green's function of the polyharmonic operator and its derivatives}, once again without any restrictions on the geometry of the domain. The full list of the estimates is quite extended. For the purposes of the introduction, let us just highlight the highest order case. 

As customary, we denote by $G(x,y)$, $x,y\in\Omega$, Green's
function for the polyharmonic equation and by $\Gamma$ its fundamental solution, so that, in particular, $G(x,y)-\Gamma(x-y)$, $x,y\in\Omega$, is the regular part of the Green function.  By definition, for every
fixed $y\in\Omega$ the function $G(\cdot ,y)$ satisfies
\begin{equation}\label{eq8.1}
(-\Delta_x)^m G(x,y)=\delta(x-y), \qquad x\in\Omega,
\end{equation}

\noindent in the space $\ring W^{m,2}(\Omega)$. Here  $\Delta_x$ stands for the Laplacian in $x$ variable, and
similarly we use the notation $\Delta_y$, $\nabla_y$, $\nabla_x$ for
the Laplacian and gradient in $y$, and gradient in $x$,
respectively. By $d(x)$ we denote the distance from $x\in\Omega$ to
$\po$.

\begin{theorem}\label{t1.2} Let $\Omega\subset \RR^n$ be an arbitrary bounded
domain.  If $n\in [3,2m+1] \cap\NN$ is odd then
\begin{equation}\label{eq8.8.1}
\left|\nabla_x^{m-\frac n2+\frac 12}\nabla_y^{m-\frac n2+\frac 12} (G(x,y)-\Gamma(x-y))\right|\leq \frac{C}{\max\{d(x),d(y),|x-y|\}}, \quad \mbox{for
every $x,y\in\Omega$, }
\end{equation}

\noindent and, in particular, 
\begin{equation}\label{eq8.5.1}
\left|\nabla_x^{m-\frac n2+\frac 12}\nabla_y^{m-\frac n2+\frac 12} G(x,y)\right|\leq  \frac{C}{|x-y|}, \quad \mbox{for
every $x,y\in\Omega$. }\end{equation}

If $n\in [2,2m]\cap\NN$ is even, then 
\begin{equation}\label{eq8.10.1}
\left|\nabla_x^{m-\frac n2}\nabla_y^{m-\frac n2} (G(x,y)-\Gamma(x-y))\right| \leq C\, \log \left(1+\frac{{\rm diam}\,\Omega}{\max\{d(x),d(y),|x-y|\}}\right),\,\, \mbox{for
every $x,y\in\Omega$, }\end{equation}

\noindent and
\begin{equation}
\label{eq8.7.1}
\left|\nabla_x^{m-\frac n2}\nabla_y^{m-\frac n2} G(x,y)\right| \leq  C  \log \left(1+\frac{\min\{d(x),d(y)\}}{|x-y|}\right),\quad \mbox{for
every $x,y\in\Omega$. }\end{equation}

The constant $C$ in \eqref{eq8.8.1}--\eqref{eq8.7.1}  depends on $m$ and $n$ only. In particular, it does not depend on the size or the geometry of the domain $\Omega$.

\end{theorem}

We mention that the pointwise bounds  on the absolute value of Green's function itself have been treated previously in dimensions $2m+1$ and $2m+2$ for  $m>2$ and dimensions $5,6,7$ for $m=2$ in \cite[Section~10]{Maz99B} (see also \cite{M1}). Also, as for the solutions, there exist results in smooth domains \cite{DASw}, \cite{Kras}, \cite{Sol1}, \cite{Sol2}, in conical domains \cite{MazPlam}, \cite{KMR}, and in polyhedra  \cite{MR1}. The estimates on the derivatives of Green's function in {\it arbitrary} bounded domains, provided by Theorem~\ref{t1.2},  are new. 

Furthermore, using standard techniques, the Green's function estimates can be employed to establish the bounds on the solution to \eqref{eq1.2} for general classes of data $f$, such as $L^p$ for a certain range of $p$, Lorentz spaces etc. We defer the detailed discussion of those to the body of the paper.

Here we just would like to point out that until recently, the properties stated in Theorems~\ref{t1.1}--\ref{t1.2} seemed accessible exclusively under heavy restrictions on the geometry of the domain (see the references above). The present paper rests on a new method, based on some intricate weighted integral identities. The biggest challenge, and the core idea, is the proper choice of the weight function $w$. It is very subtle and finely tuned to the underlying elliptic operator in such a way that the positivity, or, rather, suitable bounds from below, could be obtained for expressions akin to $\int_\Omega (-\Delta)^m u(x) u(x)\, w(x)\, dx.$ When such integrals have been considered before (in particular, we have established the three-dimensional biharmonic version of the results of this paper earlier in  \cite{MayMaz2}), the difficulties of handling general $m$ and $n$  seemed insurmountable. One of the main technical achievement of this paper is the novel systematic construction of the weight leading to sharp bounds for the solutions. It invokes numerous new aspects: employing induction in eigenvalues of the Laplace-Beltrami operator on the sphere, preservation of some positivity properties under a change of underlying higher order operator, exploiting delicate peculiarities of $(-\Delta)^m$ depending on the parity of $m,n, m-n/2$, and others. 
The construction appears for the first time in this manuscript and is likely to be applicable to general classes of elliptic equations.

\vskip 0.08 in {\it Acknowledgements}. We are greatly indebted to Marcel Filoche for the idea relating certain positivity properties of one-dimensional differential operators to particular configurations of the
roots of associated polynomials. It has been reflected in Section~\ref{s3} and it has ultimately significantly influenced our technique. 

\section{Integral inequalities and global estimate: the case of odd dimension. Part I: power weight}\label{s2}
\setcounter{equation}{0}

Let us start with a list of notation and conventions used throughout the paper.

For any domain $\Omega\subset \RR^n$ a function  $u\in
C_0^\infty(\Omega)$ can be extended by zero to $\RR^n$ and we will
write $u\in C_0^\infty(\RR^n)$ whenever convenient. Similarly, the
functions in $\ring W_2^m(\Omega)$, $m\in\NN$, will be extended by zero and
treated as functions on $\RR^n$ or other open sets containing $\Omega$ without further comments.

The symbols $B_r(Q)$ and $S_r(Q)$ denote,
respectively, the ball and the sphere with radius $r$ centered at
$Q$ and $C_{r,R}(Q)=B_R(Q)\setminus\overline{B_r(Q)}$. When the
center is at the origin, we write $B_r$ in place of $B_r(O)$, and
similarly $S_r:=S_r(O)$ and $C_{r,R}:=C_{r,R}(O)$.

Let $(r,\omega)$ be spherical coordinates in $\RR^n$, $n\geq 2$, i.e. $r=|x|\in
(0,\infty)$ and $\omega=x/|x|$ is a point of the unit sphere $S^{n-1}$.
In fact, it will be more convenient to use $e^{-t}$, $t\in\RR$, in place of $r$, so that $t=\log r^{-1}=\log |x|^{-1}$. Then by  
$\varkappa$ we denote the mapping
\begin{equation}\label{eq2.2}
\RR^n \ni
x\,\stackrel{\varkappa}{\longrightarrow} \,(t,\omega)\in\RR\times
S^{n-1},\quad n\geq 2.
\end{equation}

\noindent The symbols $\delta_\omega$ and $\nabla_\omega$ refer,
respectively, to the Laplace-Beltrami operator and the gradient on
$S^{n-1}$.

Finally, by $C$, $c$, $C_i$ and $c_i$, $i\in\NN$, we generally denote some
constants, possibly depending on the order of operator $m$ and the dimension $n$ but not on any other variables and not on the domain, unless explicitly stated otherwise. Their exact values are of no importance and can change from line to line. Also, we write
$A\approx B$, if $C^{-1}\,A\leq B\leq C\,A$ for some $C>0$.

\begin{theorem}\label{t2.1} Assume that $m\in\NN$ and $n\in [3,2m+1]\cap \NN$
is odd. Let $\Omega$ be a bounded domain in $\RR^n$, $O\in
\RR^n\setminus\Omega$, $u\in C_0^\infty(\Omega)$ and
$v=e^{\left(m-\frac n2+\frac 12\right)t}(u\circ \varkappa^{-1})$.
Then
\begin{eqnarray}\label{eq2.3}
\hskip -1cm &&\int_{\RR^n}(-\Delta)^m u(x)\,u(x)|x|^{-1}\,dx\nonumber\\[4pt]
\hskip -1cm &&\quad \geq C\sum_{k=1}^m \int_{\RR}\int_{S^{n-1}}
\left(\partial_t^k v\right)^2\,d\omega dt+C
\int_{\RR}\int_{S^{n-1}} v\,\prod_{p=-\frac n2+\frac 32}^{m-\frac
n2+\frac 12}\left(-\delta_\omega-p\,(p+n-2)\right) v\,d\omega dt,
\end{eqnarray}

\noindent where $C>0$ is some constant depending on $m$ and $n$
only.
\end{theorem}

\vskip 0.08in

\bp {\bf Step I.} In the system of coordinates $(t,\omega)$
the polyharmonic operator can be written as
\begin{equation}\label{eq2.4}
(-\Delta)^m=(-1)^me^{2mt}\prod_{j=0}^{m-1}\Bigl((-\partial_t-2j)(-\partial_t-2j+n-2)+\delta_\omega\Bigr).
\end{equation}

\noindent Then
\begin{equation}\label{eq2.5}
\int_{\RR^n}(-\Delta)^m
u(x)\,u(x)|x|^{-1}\,dx=\int_{\RR}\int_{S^{n-1}} {\mathcal
L}^{m,n}(\partial_t,\delta_\omega)v(t,\omega)\,v(t,\omega)\,d\omega dt,
\end{equation}

\noindent with
\begin{equation}\label{eq2.6}
{\mathcal L}^{m,n}(\partial_t,\delta_\omega)=(-1)^m\prod_{j=0}^{m-1}
\Bigg(\Bigl(-\partial_t+m-\frac n2+\frac
12-2j\Bigr)\Bigl(-\partial_t+m+\frac n2-\frac
32-2j\Bigr)+\delta_\omega\Bigg).
\end{equation}

\noindent Denote by $v_{pl}$ the coefficients of the expansion of
$v$ into spherical harmonics:
\begin{equation}\label{eq2.7}
v(t,\omega)=\sum_{p=0}^\infty \sum_{l=-p}^{p} v_{pl}(t)
Y_l^p(\omega), \qquad t\in\RR, \,\,\omega\in S^{n-1}.
\end{equation}

\noindent Then we can write the expression on the right-hand side
of (\ref{eq2.5}) as
\begin{equation}\label{eq2.8}
\sum_{p=0}^\infty \sum_{l=-p}^{p} \int_{\RR} {\mathcal
L}^{m,n}(\partial_t,-p\,(p+n-2)) v_{pl}(t) \,v_{pl}(t)\,dt.
\end{equation}

\noindent As usually, denote by $\widehat v$ the Fourier transform
of $v$, i.e.
\begin{equation}\label{eq2.9}
\widehat v(\gamma)=\frac{1}{\sqrt{2\pi}} \int_\RR e^{-i\gamma\,
t}v(t)\,dt,\qquad\gamma\in\RR.
\end{equation}

\noindent By the Plancherel's identity (\ref{eq2.8}) is equal to
\begin{eqnarray}\label{eq2.10}\nonumber
&&\sum_{p=0}^\infty \sum_{l=-p}^{p} \int_{\RR} {\mathcal
L}^{m,n}(i\gamma,-p\,(p+n-2))
\left|\widehat{v_{pl}}(\gamma)\right|^2\,d\gamma\\[4pt]
&&\qquad\qquad = \sum_{p=0}^\infty \sum_{l=-p}^{p} \int_{\RR} \Re
e\, {\mathcal L}^{m,n}(i\gamma,-p\,(p+n-2))
\left|\widehat{v_{pl}}(\gamma)\right|^2\,d\gamma.
\end{eqnarray}

\noindent Indeed, the imaginary part of the polynomial ${\mathcal
L}^{m,n}(i\gamma,-p\,(p+n-2))$ is odd, while $\widehat{v_{pl}}$ is
even for any real-valued function $v_{pl}$, hence the imaginary
part of the first integral above is equal to zero.

\vskip 0.08in

From now we shall carry out all calculations for $n=3$ (Steps
II-VII) and then show that the general case can be reduced to the
three-dimensional one (Step VIII). For brevity, ${\mathcal L}^{m,3}$
will be denoted ${\mathcal L}^{m}$.

\vskip 0.08in \noindent  {\bf Step II.} We claim that ${\mathcal
L}^m$, $m\in\NN$, satisfy the following relations:
\begin{equation}\label{eq2.11}
{\mathcal L}^m(i\gamma,-p\,(p+1))=a_p+i\,b_p,
\end{equation}

\noindent where
\begin{eqnarray}\label{eq2.12}
a_0 & = & \prod_{k=0}^{m-1} \left(\gamma^2+k^2\right),\\[4pt]
a_1 & = & (\gamma^2+m^2+1)\prod_{k=0}^{m-2}
 \left(\gamma^2+k^2\right),\label{eq2.13}\\[4pt]
a_p &=&
\frac{(2p-1)a_{p-1}+\left(\gamma^2+(p-1+m)^2\right)a_{p-2}}{\gamma^2+(m-p)^2},\quad
p\geq 2,\label{eq2.14}
\end{eqnarray}

\noindent and
\begin{eqnarray}\label{eq2.15}
b_0 & = & \gamma
 m\prod_{k=1}^{m-1} \left(\gamma^2+k^2\right),\\[4pt]
 b_1 & = & \gamma
 m(\gamma^2+m^2-1)\prod_{k=1}^{m-2}
 \left(\gamma^2+k^2\right),\label{eq2.16}\\[4pt]
b_p &=&
\frac{-(2p-1)b_{p-1}+\left(\gamma^2+(p-1+m)^2\right)b_{p-2}}{\gamma^2+(m-p)^2},\quad
p\geq 2.\label{eq2.17}
\end{eqnarray}

\noindent When $p=m$ and $\gamma=0$, the expressions in
(\ref{eq2.14}) and (\ref{eq2.17}) should be understood in the sense
of the corresponding limits. Formulas (\ref{eq2.12})--(\ref{eq2.17})
show that each of the polynomials $a_{m-1}$, $a_m$, $b_{m-1}$,
$b_m$ contains the factor $\gamma^2$, hence the aforementioned
limits are finite. Below we shall not consider separately the case
$p=m$, $\gamma=0$, as all arguments can be justified for this case
by a simple limiting procedure.

Let us denote ${\mathcal L}^m(i\gamma,-p\,(p+1))$ by ${\mathcal L}_p$
throughout this argument. According to (\ref{eq2.6})
\begin{eqnarray}\label{eq2.18}\nonumber
{\mathcal L}_p&=&(-1)^m\prod_{j=0}^{m-1}
\Bigl((-i\gamma+m-1-2j)(-i\gamma+m-2j)-p\,(p+1)\Bigr)\\[4pt]
&=& (-1)^m\prod_{j=0}^{m-1} (-i\gamma+m-1-2j-p)(-i\gamma+m-2j+p).
\end{eqnarray}

\noindent It remains to show that the polynomials above satisfy
(\ref{eq2.11})--(\ref{eq2.17}). First,
\begin{eqnarray}\label{eq2.19}\nonumber
{\mathcal L}_0&=&(-1)^m\prod_{j=0}^{m-1}
(-i\gamma+m-1-2j)(-i\gamma+m-2j)=(-1)^m\prod_{k=-m+1}^{m}
(-i\gamma+k)\\[4pt]&=&(-1)^m\,(-i\gamma)(m-i\gamma)\prod_{k=1}^{m-1}
(-i\gamma+k)(-i\gamma-k)=(\gamma^2+i\gamma
m)\prod_{k=1}^{m-1}(\gamma^2+k^2)
\end{eqnarray}

\noindent and
\begin{eqnarray}\label{eq2.20}\nonumber
{\mathcal L}_1&=&(-1)^m\prod_{j=0}^{m-1}
(-i\gamma+m-2-2j)(-i\gamma+m-2j+1)\\[4pt] \nonumber &=&(-1)^m\,(m+1-i\gamma)(m-1-i\gamma)(-i\gamma)(-m-i\gamma)\prod_{k=1}^{m-2}
(-i\gamma+k)(-i\gamma-k) \\[4pt]\nonumber &=&(m^2-2mi\gamma-\gamma^2-1)(mi\gamma-\gamma^2)\prod_{k=1}^{m-2}
(\gamma^2+k^2)\\[4pt] &=&\Bigl(\gamma^2(\gamma^2+m^2+1)+mi\gamma(\gamma^2+m^2-1)\Bigr)\prod_{k=1}^{m-2}
(\gamma^2+k^2),
\end{eqnarray}

\noindent as desired. As for $p\geq 2$, we have to show that
\begin{eqnarray}\label{eq2.21}\nonumber
{\mathcal
L}_p&=&\frac{1}{\gamma^2+(m-p)^2}\Bigl((2p-1)(a_{p-1}-ib_{p-1})
+\left(\gamma^2+(p-1+m)^2\right)(a_{p-2}+ib_{p-2}) \Bigr)\\[4pt]
&=&\frac{1}{\gamma^2+(m-p)^2}\Bigl((2p-1)\overline{{\mathcal L}_{p-1}}
+\left(\gamma^2+(p-1+m)^2\right){\mathcal L}_{p-2} \Bigr).
\end{eqnarray}

\noindent It follows from (\ref{eq2.18}) that
\begin{eqnarray}\label{eq2.22}\nonumber
{\overline{{\mathcal L}_{p-1}}}&=&(-1)^m\prod_{j=0}^{m-1}
(i\gamma+m-2j-p)(i\gamma+m-2j+p-1)\\[4pt]\nonumber
&=&(-1)^m\prod_{j=0}^{m-1}
(i\gamma+m-2(m-1-j)-p)(i\gamma+m-2(m-1-j)+p-1)\\[4pt]
&=&(-1)^m\prod_{j=0}^{m-1} (m-2j-i\gamma-p-1)(m-2j-i\gamma+p-2),
\end{eqnarray}

\noindent and
\begin{eqnarray}\label{eq2.23}\nonumber
{\mathcal L}_{p-2}&=&(-1)^m\prod_{j=0}^{m-1}
(-i\gamma+m-2j-p+1)(-i\gamma+m-2j+p-2)\\[4pt]
&=&(-1)^m \prod_{j=0}^{m-1}(-i\gamma+m-2j+p-2)\prod_{k=-1}^{m-2}
(-i\gamma+m-2k-p-1).
\end{eqnarray}

\noindent Hence,
\begin{eqnarray}\label{eq2.24}\nonumber
&&\frac{1}{\gamma^2+(m-p)^2}\Bigl((2p-1)\overline{{\mathcal L}_{p-1}}
+\left(\gamma^2+(p-1+m)^2\right){\mathcal L}_{p-2} \Bigr)\\[4pt]\nonumber
&&\qquad =\frac{1}{\gamma^2+(m-p)^2}\,(-1)^m
\prod_{j=0}^{m-1}(-i\gamma+m-2j+p-2) \prod_{k=0}^{m-2}
(-i\gamma+m-2k-p-1)\times\\[4pt]\nonumber
&&\qquad \,\,\times \Bigl((2p-1)(-m+1-i\gamma-p)
+\left(\gamma^2+(p-1+m)^2\right)(m+1-i\gamma-p) \Bigr)\\[4pt]\nonumber
&&\qquad =\frac{1}{\gamma^2+(m-p)^2}\,(-1)^m
\prod_{k=1}^{m}(-i\gamma+m-2k+p) \prod_{k=0}^{m-2}
(-i\gamma+m-2k-p-1)\times\\[4pt]
&&\qquad \,\,\times (-m+1-i\gamma-p) \Bigl((2p-1)
-(m-1-i\gamma+p)(m+1-i\gamma-p) \Bigr).
\end{eqnarray}

\noindent The latter expression is equal to
\begin{eqnarray}\label{eq2.25}
\nonumber &&\qquad {\mathcal L}_p\,\,\frac{1}{\gamma^2+(m-p)^2}\,\,
\frac{(-i\gamma-m+p)}{(-i\gamma+m+p)}\,\,\frac{1}{(-m+1-i\gamma-p)}
\times\\[4pt]
&&\qquad\quad \,\,\times (-m+1-i\gamma-p)
(-1)(-i\gamma+m+p)(-i\gamma+m-p)\,=\,{\mathcal L}_p.
\end{eqnarray}

\noindent This finishes the proof of (\ref{eq2.11})--(\ref{eq2.17}),
in particular, we have
\begin{equation}\label{eq2.26}
\Re e\, {\mathcal L}^m(i\gamma,-p\,(p+1))=a_p \mbox{ for every
}p\in\NN\cup\{0\},
\end{equation}

\noindent and hence
\begin{equation}\label{eq2.27}
\int_{\RR^3}(-\Delta)^m u(x)\,u(x)|x|^{-1}\,dx= \sum_{p=0}^\infty
\sum_{l=-p}^{p}\int_{\RR} a_p(\gamma)
\left|\widehat{v_{pl}}(\gamma)\right|^2\,d\gamma.
\end{equation}

\vskip 0.08in \noindent  {\bf Step III.} In order to prove
(\ref{eq2.3}) for $n=3$ one has to show that
\begin{equation}\label{eq2.28}
a_q(\gamma) \geq C\sum_{k=1}^m \gamma^{2k} +C
\prod_{p=0}^{m-1}\Bigl(q(q+1)-p\,(p+1)\Bigr) \mbox{ for every
$q\in\NN\cup\{0\}$, $\gamma\in\RR$},
\end{equation}

\noindent with some constant $C>0$ independent of $q$ and
$\gamma$. Indeed, if (\ref{eq2.28}) is satisfied, then
\begin{eqnarray}\label{eq2.29}
\nonumber &&\sum_{q=0}^\infty \sum_{l=-q}^{q}\int_{\RR}
a_q(\gamma)
\left|\widehat{v_{ql}}(\gamma)\right|^2\,d\gamma\\[4pt]
&& \qquad \geq C \sum_{q=0}^\infty \sum_{l=-q}^{q}\int_{\RR}
\left(\sum_{k=1}^m\left|(i\gamma)^k\,
\widehat{v_{ql}}(\gamma)\right|^2+
{\overline{\widehat{v_{ql}}(\gamma)}}
\,\prod_{p=0}^{m-1}\Bigl(q(q+1)-p\,(p+1)\Bigr)\widehat{v_{ql}}(\gamma)\right)\,dt
\nonumber \\[4pt]
&& \qquad = C \sum_{q=0}^\infty \sum_{l=-q}^{q}\int_{\RR}
\left(\sum_{k=1}^m\left(\partial_t^k v_{ql}(t)\right)^2+ v_{ql}(t)
\prod_{p=0}^{m-1}\Bigl(q(q+1)-p\,(p+1)\Bigr)v_{ql}(t)\right)\,dt
\nonumber \\[4pt]
&& \qquad = C \int_{\RR}
\int_{S^2}\left(\sum_{k=1}^m\left(\partial_t^k
v(t,\omega)\right)^2+ v(t,\omega)
\prod_{p=0}^{m-1}\Bigl(-\delta_\omega-p\,(p+1)\Bigr)v(t,\omega)\right)\,dt
\end{eqnarray}

\noindent which combined with (\ref{eq2.27}) gives (\ref{eq2.3})
 for $n=3$. Let us now concentrate on (\ref{eq2.28}).

\vskip 0.08in \noindent  {\bf Step IV.} The formulas
(\ref{eq2.12})--(\ref{eq2.14}) show that all $a_q$ are nonnegative.
Moreover, since
\begin{equation}\label{eq2.30}
(q-1+m)^2>(m-q)^2 \mbox{ for every $q\in\NN$ and $m\in\NN$},
\end{equation}

\noindent we have $a_q\geq a_{q-2}$ for all $q\geq 2$, and
therefore,
\begin{equation}\label{eq2.31}
a_q(\gamma)\geq \min\{a_0(\gamma),a_1(\gamma)\} \geq
\prod_{k=0}^{m-1}(\gamma^2+k^2) \mbox{ for every
$q\in\NN\cup\{0\}$, $\gamma\in\RR$}.
\end{equation}

\noindent Hence,
\begin{eqnarray}\label{eq2.32}
a_q(\gamma) \geq \sum_{k=1}^m \gamma^{2k} \mbox{ for every
$q\in\NN\cup\{0\}$, $\gamma\in\RR$}.
\end{eqnarray}

\vskip 0.08in \noindent  {\bf Step V.} As for the second term on
the right hand side of (\ref{eq2.28}), it is clear that
\begin{equation}\label{eq2.33}
a_q(\gamma) \geq \prod_{p=0}^{m-1}\Bigl(q(q+1)-p\,(p+1)\Bigr)
\mbox{ for $q=0,1,..., m-1$, $\gamma\in\RR$},
\end{equation}

\noindent since for such $q$ the product in (\ref{eq2.33}) is equal
to 0. Now assume that $q\geq m$. If $|\gamma|\geq q+1$, then by
(\ref{eq2.31})
\begin{equation}\label{eq2.34}
a_q(\gamma)\geq  \prod_{k=0}^{m-1}\Bigl((q+1)^2+k^2\Bigr) \geq
\prod_{p=0}^{m-1}\Bigl(q(q+1)-p\,(p+1)\Bigr),
\end{equation}

\noindent and it remains to consider $|\gamma|\leq q+1$, $q\geq
m$.

\vskip 0.08in \noindent  {\bf Step VI.} By definition
\begin{eqnarray}\label{eq2.35}
a_q(\gamma) &=& (-1)^m \,\Re e \,\prod_{j=0}^{m-1}
(-i\gamma+m-1-2j-q)(-i\gamma+m-2j+q)\nonumber\\[4pt]
&=& \Re e \,\prod_{j=0}^{m-1} \Bigl(\gamma^2
+(-m+1+2j+q)(m-2j+q)+i\gamma(2m-4j-1)\Bigr).
\end{eqnarray}

\noindent Since $q\geq m$, we have $(-m+1+2j+q)>0$ and
$(m-2j+q)>0$. Therefore,
\begin{eqnarray}\label{eq2.36}
a_q(\gamma) &=& \prod_{j=0}^{m-1} \Bigl(\gamma^2
+(-m+1+2j+q)(m-2j+q)\Bigr)\nonumber\\[4pt]&& \times \,\Re e \,\prod_{j=0}^{m-1}
\Bigl(1+i\gamma\frac{(2m-4j-1)}{\gamma^2
+(-m+1+2j+q)(m-2j+q)}\Bigr).
\end{eqnarray}

\noindent Let us show that there exists $C(m)\geq m$ such that for
every $q\geq C(m)$ and $|\gamma|\leq q+1$

\begin{equation}\label{eq2.37}
\Re e \,\prod_{j=0}^{m-1} \Bigl(1+i\gamma\frac{(2m-4j-1)}{\gamma^2
+(-m+1+2j+q)(m-2j+q)}\Bigr)\geq \frac 12.
\end{equation}

\noindent A general idea is that the expression above can be
written as $1+\Re e \,R$, where $|R|$ is sufficiently small when $q$
is large. Indeed, $R$ is a sum of at most $2^m-1$ terms, each of
those being a product of (at least one, at most $m$) elements

\begin{equation}\label{eq2.38}
c_j:=i\gamma\frac{(2m-4j-1)}{\gamma^2 +(-m+1+2j+q)(m-2j+q)}, \quad
j=0,...,m-1.
\end{equation}

\noindent If $q\geq 2m$, then $|2j-m|\leq q/2$ and

\begin{equation}\label{eq2.39}
\left| c_j\right|\leq (q+1)\,\frac{2m}{(q/2+1)(q/2)}\leq
\frac{8m}{q}.
\end{equation}

\noindent Moreover, if $q> 8m$, the expression above is less than
1. Hence,

\begin{equation}\label{eq2.40}
|\Re e \,R|\leq |R|\leq (2^m-1)\,\frac{8m}{q}.
\end{equation}

\noindent Set now $C(m):= 2^{m+4} m$. Then $|\Re e R|\leq 1/2$ for
every $q\geq C(m)$, and hence (\ref{eq2.37}) holds for such $q$.
Combining this with (\ref{eq2.36}), we deduce that

\begin{eqnarray}\label{eq2.41}
a_q(\gamma) &\geq& \frac 12 \prod_{j=0}^{m-1} \Bigl(\gamma^2
+(-m+1+2j+q)(m-2j+q)\Bigr)\nonumber\\[4pt]
&\geq& \frac 12 \prod_{j=0}^{m-1} (-m+1+2j+q)(m-2j+q)\geq
\frac 12\, \Bigg(\frac{(q+2)q}{4}\Bigg)^m \nonumber\\[4pt]
&\geq&\frac 1{2^{2m+1}} \prod_{p=0}^{m-1}\Bigl(q(q+1)-p\,(p+1)\Bigr)
\end{eqnarray}

\noindent whenever $|\gamma|\leq q +1$, $q\geq  2^{m+4} m$.

\vskip 0.08in \noindent  {\bf Step VII.} Finally, we claim that
there exists $D(m)>0$ such that
\begin{equation}\label{eq2.42}
a_q(\gamma)\geq D(m)
\prod_{p=0}^{m-1}\Bigl(q(q+1)-p\,(p+1)\Bigr)\mbox{ if $m\leq q
\leq 2^{m+4} m$ and $|\gamma|\leq q+1$.}
\end{equation}

\noindent Since $a_q\geq a_{q-2}$ for $q\geq 2$, it is enough to
show that for some $D>0$

\begin{equation}\label{eq2.43}
\min\{a_m(\gamma),a_{m+1}(\gamma)\}\geq D, \qquad |\gamma|\leq
q+1,
\end{equation}

\noindent as (\ref{eq2.43}) implies (\ref{eq2.42}) with

\begin{equation}\label{eq2.44}
D(m)=D\left( \prod_{p=0}^{m-1}\Bigl(2^{m+4} m(2^{m+4}
m+1)-p\,(p+1)\Bigr) \right)^{-1}.
\end{equation}

However, the formula (\ref{eq2.18}) implies

\begin{equation}\label{eq2.45}
a_m(0)=\prod_{j=0}^{m-1} (1+2j)(2m-2j)\geq 2 \quad\mbox{and}\quad
a_{m+1}(0)=\prod_{j=0}^{m-1} (2+2j)(2m+1-2j)\geq 2,
\end{equation}

\noindent so that for some $\eps_0=\eps_0(m)>0$ we have
$\min\{a_m(\gamma),a_{m+1}(\gamma)\}\geq 1$ when
$|\gamma|<\eps_0$. On the other hand, by (\ref{eq2.31})

\begin{equation}\label{eq2.46}
a_q(\gamma)\geq \prod_{k=0}^{m-1}(\eps_0^2+k^2) \mbox{ for every
$q\in\NN\cup\{0\}$, $|\gamma|>\eps_0$}.
\end{equation}

\noindent Hence, (\ref{eq2.43}) holds with $D$ given by the minimum
of 1 and the right hand side of (\ref{eq2.46}). This concludes the
proof of (\ref{eq2.42}) and the proof of the Theorem for the case
$n=3$.

\vskip 0.08in \noindent  {\bf Step VIII.} Let us now consider any
odd $n\geq 3$. Similarly to (\ref{eq2.28})--(\ref{eq2.29}) in
order to prove (\ref{eq2.3}) one has to show that
\begin{equation}\label{eq2.47}
\Re e \,{\mathcal L}^{m,n}(i\gamma, -q(q+n-2)) \geq C\sum_{k=1}^m
\gamma^{2k} +C \prod_{p=-\frac n2+\frac 32}^{m-\frac n2+\frac
12}\Bigl(q(q+n-2)-p\,(p+n-2)\Bigr)
\end{equation}

\noindent for every $q\in\NN\cup\{0\}$, $\gamma\in\RR$. However,

\begin{eqnarray}\label{eq2.48}\nonumber
&&{\mathcal L}^{m,n}(i\gamma, -q\,(q+n-2))\\[4pt]\nonumber
&&\quad =(-1)^m\prod_{j=0}^{m-1} \Bigl(-i\gamma+m-\frac n2+\frac
12-2j-q\Bigr)\Bigl(-i\gamma+m+\frac n2-\frac 32-2j+q\Bigr)\\[4pt]\nonumber
&&\quad =(-1)^m\prod_{j=0}^{m-1}
(-i\gamma+m-1-2j-p)(-i\gamma+m-2j+p)\\[4pt]
&&\quad= {\mathcal L}^{m,3}(i\gamma, -p\,(p+1)),
\end{eqnarray}

\noindent for $p=q +\frac n2-\frac 32$. Hence,

\begin{eqnarray}\label{eq2.49}\nonumber
&&\Re e\,{\mathcal L}^{m,n}(i\gamma, -q\,(q+n-2))=a_{q +\frac n2-\frac 32}(\gamma)\\[4pt]\nonumber
&&\quad \geq C\sum_{k=1}^m \gamma^{2k} +C
\prod_{p=0}^{m-1}\Bigg(\Bigl(q +\frac n2-\frac 32\Bigr)\Bigl(q +\frac n2-\frac 12\Bigr)-p\,(p+1)\Bigg)\\[4pt]\nonumber
&&\quad = C\sum_{k=1}^m \gamma^{2k} +C \prod_{s=-\frac n2+\frac
32}^{m-\frac n2+\frac 12}\Bigg(\Bigl(q +\frac n2-\frac
32\Bigr)\Bigl(q +\frac n2-\frac 12\Bigr)-\Bigl(s +\frac n2-\frac
32\Bigr)\Bigl(s+\frac n2-\frac 12\Bigr)\Bigg)\\[4pt]
&&\quad = C\sum_{k=1}^m \gamma^{2k} +C \prod_{s=-\frac n2+\frac
32}^{m-\frac n2+\frac 12}\Bigl(q(q+n-2)-s(s+n-2)\Bigr)
\end{eqnarray}

\noindent where we used (\ref{eq2.28}) and substitution
$s:=p-\frac n2+\frac 32$.  As in \eqref{eq2.28} the constant $C>0$ above is independent of $q$ and $\gamma$. Now (\ref{eq2.49}) leads to
(\ref{eq2.47}) and finishes the argument.
 \ep
 
 \comment{\todo{Note to self: In even dimensions it also can be expected that $m-n/2+1/2$ derivatives are bounded, at least the counterexamples don't say otherwise: when we are outside a ray $m-n/2+1/2$ is the first eigenvalue (see \cite{KMR}) and outside a point we can only prove that $m-n/2+1$ derivatives are not continuous, because it has to be an integer number (of derivatives of fundamental solution) strictly bigger than $m-n/2$ -- checked! Anyways, it is still possible that the sharp {\it non-integer} result is not $m-n/2$ in even dimensions.
 
BUT the approach above can not work to prove that for even dimensions we can bound $m-n/2+1/2$ derivatives too. Using formula (\ref{eq2.6}) we can write explicitly ${\mathcal L}^{m,n}(0,-q(q+n-2))$ and to see that it alternates sign through $q\in\NN\cup\{0\}$ when $n$ is even. Just take $q=...m-n/2-1,m-n/2,m-n/2+1...$ Inserting $\ln (L+t)$ will not change the situation either.}
}
\section{Preservation of positivity for solutions of differential operators} \label{s3}
\setcounter{equation}{0}

%This section is aimed at proving the estimate \eqref{eq2.11.6-n7.1} in the assumptions of Lemma~\ref{l2.3}. In other words, we establish the following result. 

\begin{proposition}\footnote{The idea of the proof has been suggested by Marcel Filoche. }\label{p3.1} Assume that $m\in\NN$ and $n\in [3,2m+1]\cap \NN$
is odd.
Let us define $h$ as the unique solution of
\begin{align}\label{eq3.1}
\mathcal L^{m,n}\left(-\partial_t,0\right)\,h = \delta,
\end{align}

\noindent which is bounded and vanishes at $+\infty$.
Then, 
\begin{equation}\label{eq3.2}
 {\mathcal L}^{m,n}\left(-\partial_t,-p(p+n-2)\right)\,h\geq 0,
 \end{equation}
for all $p \in [0,m-n/2+1/2]$.
\end{proposition}

Let us start with the operators of the second order.

Let $a$, $b$, $c$, and $d$ four real numbers such that $b<0\leq a$ and  $d<0\leq c$ and let $f$ be a (tempered) distribution on $\RR$. We define  $h$ as a solution of
\begin{equation}\label{eq3.3}
\left(\partial_t - a\right)\left(\partial_t - b \right) h = f,
\end{equation}
and let the distribution $g$ be given by
\begin{equation}\label{eq3.4}
g := \left(\partial_t - c\right)\left(\partial_t -d\right) h
\end{equation}

\noindent We want to investigate under what condition on $a$, $b$, $c$, and $d$ the positivity of $f$ entails the positivity of $g$.

%When both $a$ and $b$ are not equal to zero, the assumption that $h$ is bounded defines it uniquely (since the solution to the corresponding homogeneous problem has exponential growth at either $+\infty$ or $-\infty$). When $a=0$, even a bounded solution is defined modulo a constant. 

First, one has to pick a suitable solution (clearly, by subtracting some solution of a homogeneous equation one can always destroy positivity).
Also, $f$ has to be reasonable in order for a solution to exist. 

To keep the discussion at full generality and, at the same time, ensure existence and uniqueness, in the preliminary results of this section we simply explicitly assume that the solution is obtained by convolution with the bounded fundamental solution vanishing at $+\infty$ (and, in particular, that $f$ is such that the convolution integral converges), and later we will check that this hypothesis fits the argument proving Proposition~\ref{p3.1}. 

Also, we record a few conventions. As usually, a distribution $f$ has the singular support at ${0}$ if $f\in C^\infty(\RR\setminus \{0\})$. Furthermore, we say that such an $f$ has an exponential decay at $+\infty$ if $|f(t)|\leq C e^{-\alpha t}, \,\, t>0,$ for some $\alpha, C>0$, and similarly, $f$ has an exponential decay at $-\infty$ if $|f(t)|\leq C e^{\beta t}, \,\, t<0$ for some $\beta, C>0$. Within this section all distributions of interest are bounded in a punctured neighborhood of $0$ and thus, there is no loss of generality in considering all $t>0$ (or, respectively, $t<0$) in the definitions above.

%\todo{Note to self: a distribution decaying at $+\infty$ can be defined, for instance, using shifts of Schwartz functions: for every fixed Schwartz function for every $\eps$ there is a shift right such that when our distribution acts on the shift, the result is less that $\eps$. I think it is well-posed - to be checked. But it is perhaps classics anyways.}

\begin{proposition} \label{p3.2}
Fix  $a, b,c,d\in\RR$ such that $b< 0\leq a$ and $d<0\leq c$. Let $f$ and $h$ be two distributions on $\RR$ such that $h$ is a convolution of $f$ with the bounded fundamental solution of the operator $\left(\partial_t - a\right)\left(\partial_t - b\right)$, which  vanishes at $+\infty$. In particular, 
\begin{align}\label{eq3.5}
\left(\partial_t - a\right)\left(\partial_t - b\right) h = f.
\end{align}

Assume that $f$ is positive. If $[b,a] \subset [d,c]$, then 
\begin{align}\label{eq3.6}
g := \left(\partial_t - c\right)\left(\partial_t - d\right) h
\end{align}
is also a positive distribution.

Furthermore, if $f$ has the singular support at ${0}$ and exponential decay at $+\infty$ and $-\infty$, then $g$ also has the singular support at ${0}$ and exponential decay at $+\infty$ and $-\infty$  provided that either $a,c \neq 0$ or $a=c=0$.
\end{proposition}

\bp Let us  start with the case $f=\delta$. Problem~\eqref{eq3.5} therefore writes
\begin{align}\label{eq3.7}
\left(\partial_t - a\right)\left(\partial_t - b \right) h = \delta.
\end{align}
The solution of \eqref{eq3.7} which is bounded and vanishes at $+\infty$ is given by
\begin{align}\label{eq3.8}
h_0(t) = \kappa
\begin{cases}
e^{at} &\qquad {\rm if}\, t \leq 0\\
e^{bt} &\qquad {\rm if}\, t > 0
\end{cases}
\qquad {\rm with}\quad \kappa = \frac{1}{b-a}.
\end{align}
This solution is unique because any other bounded solution vanishing at $+\infty$, $h_1$, would be such that $h_1 - h_0$ solves the homogeneous problem corresponding to \eqref{eq3.5}. Thus, $h_1 - h_0$ would have an exponential growth at $-\infty$ and/or violate the assumption of vanishing at $+\infty$.

Applying the operator in \eqref{eq3.6} to this solution $h_0$ yields 
\begin{align}\label{eq3.9}
g_0(t) = \delta + 
\begin{cases}
\kappa~(a-c)(a-d)\,e^{at} &\qquad {\rm if}\, t \leq 0,\\
\kappa~(b-c)(b-d)\,e^{bt} &\qquad {\rm if}\, t > 0,
\end{cases}
\end{align}
i.e., $g_0$ is the sum of the Dirac delta function, coming from the second derivative appearing in \eqref{eq3.6}, and a piecewise continuous function. The constant $\kappa$ being always negative, it follows immediately  that $g_0$ will be positive if and only if $(a-c)(a-d)$ and $(b-c)(b-d)$ are both non-positive, i.e., if both $a$ and $b$ belong to the interval $[d,c]$. On the contrary, due to the Dirac distribution present in $g_0$, if either $a$ or $b$ is outside the interval $[d,c]$ then $g_0$ has no specific sign.

We now turn back to the initial problem \eqref{eq3.5}. As per our assumptions, we consider the solution of this problem obtained from a convolution product of $f$ with the fundamental solution $h_0$ defined in \eqref{eq3.8}, that is,
\begin{align}\label{eq3.10}
h(t) = \int_\RR f(u)\,h_0(t-u)\,du.
\end{align}

\noindent  Applying the operator in  \eqref{eq3.6} to this $h$ yields 
\begin{align}\label{eq3.11}
g(t) = \int_\RR f(u)\,g_0(t-u)\,du
\end{align}
We see here that the positivity of $g$ for any positive $f$ follows from  the positivity of $g_0$, hence $a$ and $b$ belonging to the interval $[d,c]$. Given that $b < 0 \leq a$, this amounts to $[b,a] \subset [d,c]$. 

Furthermore, if $f$ has an exponential decay at both $+\infty$ and $-\infty$ and either $a,c\neq 0$  or $a=c=0$ then, due to formula \eqref{eq3.9}, $g_0$ also has an exponential decay at both $+\infty$ and $-\infty$, and hence, so does their  convolution in \eqref{eq3.11}, as desired. 
\ep 
%To summarize, if $h$ is a reasonable solution of \eqref{eq1.0} for a positive distribution $f$, then one obtains another positive distribution $g$ when replacing in this equation $a$ and $b$ by $c$ and $d$ such that $[d,c] \supset [b,a]$.

\begin{corollary}\label{c3.3}
Let $P$ and $Q$ be two polynomials of order $2m$, with all their roots real and their higher order coefficient equal to 1. Therefore $P$ and $Q$ write
\begin{align}\label{eq3.12}
P(x) = \prod_{i=0}^{m-1} (x-a_i)(x-b_i) \qquad {\rm and} \qquad Q(x) = \prod_{i=0}^{m-1} (x-c_i)(x-d_i).
\end{align}
Let us assume further that both $P$ and $Q$ have $m$ non-negative roots, noted respectively $(a_i)$ and $(c_i)$ ($0 \leq i \leq m-1$), and $m$ negative roots, noted respectively $(b_j)$ and $(d_j)$ ($1 \leq j \leq m$), such that $\forall i \in [0,m-1], [b_i,a_i] \subset [d_i,c_i]$. Moreover, either none of the roots of both polynomials is zero, or both $P$ and $Q$ have exactly one zero root $a_0=c_0=0$.

Let $f$ be a positive distribution with the singular support at $0$   and exponential decay at $+\infty$ and $-\infty$, and  $h$ be the unique bounded solution of
\begin{align}\label{eq3.13}
P(\partial_t)\,h = f,
\end{align}
which vanishes at $+\infty$.
Then 
\begin{align}\label{eq3.14}
g := Q(\partial_t)\,h,
\end{align}
is also a positive distribution with the singular support at $0$  and exponential decay at $+\infty$ and $-\infty$.
\end{corollary}

\noindent {\it Remark.} The gist of this Proposition, as well as Proposition~\ref{p3.2}, is preservation of the positivity under appropriate change of a differential operator, or more precisely, under the change that widens the distances between the roots. The exponential decay is not a necessary condition, but rather a convenient restriction that ensures the choice of a suitable solution at each step.

\comment{
\todo{Note to self: the special condition about zero roots is not necessary either, it just makes an exposition cleaner}}

\vskip 0.08in

\bp First of all, let us mention that there exists a solution $h$ of \eqref{eq3.13} which is bounded and vanishes at $+\infty$. It is the solution obtained by the convolution of $f$ with the fundamental solution of the operator $P(\partial_t)$ which is bounded and vanishes at $+\infty$.

The Corollary will be obtained by successively applying Proposition~\ref{p3.2} for all pairs $[b_i,a_i]$. The function $f_0$ is defined as $f_0 = f$, and $f_{i+1}$ is deduced from $f_i$ by
\begin{align}\label{eq3.15}
f_{i+1} & := \left(\partial_t - c_i\right) \left(\partial_t - d_i\right) h_i,
\end{align}
where $h_i$ is the  solution of
\begin{align}\label{eq3.16}
\left(\partial_t - a_i\right) \left(\partial_t - b_i\right) h_i = f_i,
\end{align}

\noindent obtained by convolution with the bounded fundamental solution vanishing at $+\infty$, following Proposition~\ref{p3.2}. Then it is immediate to see that, by Proposition~\ref{p3.2}, all $f_i$ are positive distributions  with the singular support at $0$  and exponential decay at $+\infty$ and $-\infty$.

We claim that 
\begin{equation}\label{eq3.17}
h_i=\prod_{j=0}^{i-1}(\partial_t-c_j)(\partial_t-d_j) \prod_{j=i+1}^{m-1}(\partial_t-a_j)(\partial_t-b_j) h, \quad i=1,...,m-1,
\end{equation}

\noindent and
\begin{equation}\label{eq3.18}
f_{i+1}=\prod_{j=0}^{i}(\partial_t-c_j)(\partial_t-d_j) \prod_{j=i+1}^{m-1}(\partial_t-a_j)(\partial_t-b_j) h, \quad i=0,...,m-1,
\end{equation}

\noindent where in the definitions of $h_0$, $h_{m-1}$, $f_m$ one of the products is void. 

First of all, letting $f_0=f$, we see that 
\begin{equation}\label{eq3.19}
h_0':= \prod_{j=1}^{m-1}(\partial_t-a_j)(\partial_t-b_j) h, \quad i=1,...,m-2,
\end{equation}

\noindent satisfies equation \eqref{eq3.16}. The only question is whether it is the same solution as $h_0$. However, much as for the operators of the second order in Proposition~\ref{p3.2}, $h_0'$ 
is bounded and vanishes at $+\infty$ for a particular $h$ that we are working with (see the beginning of the proof). Clearly, the same is true for $h_0$ obtained from \eqref{eq3.16} by convolution of $f$ with the bounded fundamental solution of $\left(\partial_t - a_0\right) \left(\partial_t - b_0\right)$ which vanishes at $+\infty$. However, then the difference $h_0-h_0'$ must solve the homogeneous equation for the operator $\left(\partial_t - a_0\right) \left(\partial_t - b_0\right)$ and at the same time, $h_0-h_0'$ must be bounded and vanish at $+\infty$, hence, $h_0-h_0'=0$. Note, further, that $f_1$ can now be defined equivalently by \eqref{eq3.18} and \eqref{eq3.15} and, moreover, $f_1$ has an exponential decay at both $+\infty$ and $-\infty$. The latter follows from \eqref{eq3.15}--\eqref{eq3.16} and the last statement of Proposition~\ref{p3.2}. 

Following the same lines, we verify \eqref{eq3.17}--\eqref{eq3.18} identifying them with functions in \eqref{eq3.15}--\eqref{eq3.16} as the unique solutions of the corresponding differential equations satisfying aforementioned decay properties. 

As we already pointed out, by virtue of Proposition~\ref{p3.2}, all $f_i$ are positive distributions  with the singular support at $0$  and exponential decay at $+\infty$ and $-\infty$. It is now clear from \eqref{eq3.18} that $f_{m} = g$. This finishes the proof of the Corollary.
\ep

\noindent {\it Proof of Proposition~\ref{p3.1}}.
 We now focus on the following differential operator
\begin{align}
& \mathcal L^{m,n}(-\partial_t,-p(p+n-2)) := \nonumber \\ & (-1)^m \prod_{j=0}^{m-1} \left(\partial_t + m -\frac{n}{2} + \frac{1}{2} -2j - p\right)\left(\partial_t + m +\frac{n}{2} - \frac{3}{2} -2j + p\right).\label{eq3.20}
\end{align}
Let us introduce the numbers 
\begin{align}\label{eq3.21}
c_j = 2j - \left(m - \frac{n-1}{2}\right), \quad 0\leq j \leq m-1.
\end{align}

\noindent Note that 
$$m+\frac n2-\frac 32 -2j+p=-m+\frac n2+\frac 12+2(m-1-j)+p.$$
The operator $\mathcal L^{m,n}$ thus writes
\begin{align}\label{eq3.22}
\mathcal L^{m,n} \left(-\partial_t,-p(p+n-2)\right) = (-1)^m \prod_{j=0}^{m-1} \left(\partial_t - c_j -p \right) \left(\partial_t + c_j + 1 + p\right).
\end{align}

\noindent The roots of the corresponding polynomial are $c_j + p$ and $-c_j-1-p$, $j=0,...,m-1$.  

 Consider an integer $p$ such that $0 \leq p \leq m - \displaystyle \frac{n-1}{2} = :k$ and let us examine the following set of roots:
\begin{align}\label{eq3.23}
S_0 = \{c_j\} \cup \{-c_j-1\} = & \{-k,-k-2,...,2(m-1)-k\} 
\\ & \cup \{-2(m-1)+k-1,-2(m-1)+k+1,...,k-1\}.\label{eq3.24}
\end{align}
This set contains $2m$ elements. The largest of the $\{c_j\}$ is larger than the largest of the $\{-c_j-1\}$:
\begin{align}\label{eq3.25}
2(m-1)-k > k-1 \qquad {\rm because} \qquad 2(m-1) \geq 2k > 2k-1.
\end{align}
Similarly, the smallest of the $\{-c_j-1\}$ is smaller than the smallest of the $\{c_j\}$:
\begin{align}\label{eq3.26}
-2(m-1)+k-1 < -k \qquad {\rm because} \qquad 2(m-1) \geq 2k > 2k-1.
\end{align}
Moreover, the $\{c_j\}$ and the $\{-c_j-1\}$ have different parities. So when these two subsets of roots overlap, the overlap covers all integers, odd and even. Hence, one can rewrite the set $S$ as the union of 3 sets, the  center interval and two ``lateral" sets, these two being respectively composed of only negative and positive numbers:
\begin{align}\label{eq3.27}
S_0 = &\{-2(m-1) + k - 1,...,-k-5,-k-3\} \,\cup \\ & \Bigl([-k-1,k]\cap\ZZ\Bigr) \,\cup \label{eq3.28}\\& \{k+2,k+4,...,2(m-1)-k\}.\label{eq3.29}
\end{align}
The first and the third set might possibly be empty. 

Next, the set of roots $S_p = \{c_j+p\} \cup \{-c_j-p-1\}$ can be described in similar terms replacing $k$ by $k-p$:
\begin{align}\label{eq3.30}
S_p = &\{-2(m-1) + k-p - 1,...,-k+p-5,-k+p-3\} \cup \\ \label{eq3.31}& \Bigl([-k+p-1,k-p]\cap\ZZ\Bigr) \cup \\\label{eq3.32} & \{k-p+2,k-p+4,...,2(m-1)-k+p\}.
\end{align}
To compare $S_0$ and $S_p$, one rewrites $S_0$ in the following way
\begin{align}\label{eq3.33}
S_0 = &\{-2(m-1) + k - 1,...,-k-5,-k-3\} \,\cup \,\Bigl([-k-1,-k+p-2] \cap\ZZ\Bigr)\\ \label{eq3.34}& \Bigl([-k+p-1,k-p]\cap\ZZ\Bigr) \,\cup \\ \label{eq3.35} & \Bigl([k-p+1,k]\cap\ZZ\Bigr) \,\cup\,\{k+2,k+4,...,2(m-1)-k\}.
\end{align}
Now the center interval is identical in both sets and composed of consecutive integers, $k-p+1$ negative numbers and $k-p+1$ non negative numbers. In particular, this interval always contains $\{-1,0\}$. Also, compared one-by-one, the roots in the first subset of $S_p$ \eqref{eq3.30},  are smaller than the roots in the first subset of $S_0$ \eqref{eq3.33}, and similarly, the roots in the last subset of $S_p$ \eqref{eq3.32},  are larger than the roots in the last subset of $S_0$ \eqref{eq3.35}.
 
We now can name the sequence of roots of $S_0$ (in the order they appear in \eqref{eq3.33}--\eqref{eq3.35}) 
\begin{equation}\label{eq3.36}
b_{m-1}<b_{m-2}<...<b_0=-1<a_0=0<a_1<...<a_m,
\end{equation}

\noindent and, respectively, the roots of $S_p$ (in the order of appearance in \eqref{eq3.30}--\eqref{eq3.32}) 
\begin{equation}\label{eq3.37}
d_{m-1}<d_{m-2}<...<d_0=-1<c_0=0<c_1<...<c_m.
\end{equation}

Then $[b_i,a_i]\subset [d_i,c_i]$ for all $i=0,...,m-1$, and moreover, $a_0=c_0=0$ are the only zero roots. Thus, we are in the setting of Corollary~\ref{c3.3}, which finishes the argument. (Note that when $m$ is odd, we apply Corollary~\ref{c3.3} with $-h$ in place of $h$).
\ep

\section{Integral inequalities and global estimate: the case of odd dimension. Part II: weight $g$}\label{s4}
\setcounter{equation}{0}

Incorporating the results of Section~\ref{s3}, in this section we introduce a new weight function $g$. Together with the inequalities in Section~\ref{s2}, this  yields further improvement of the key integral estimates and ultimately,  improved pointwise estimates on the solution of the polyharmonic equation. 

We start with the following auxiliary result which provides an explicit formula for the solution of \eqref{eq3.1}.

\begin{lemma}\label{l4.1} Assume that $m\in\NN$ and $n\in [3,2m+1]\cap \NN$
is odd. Consider the equation
\begin{equation}
\Lmn(-\partial_t,0)\,h
 =\delta,  \label{eq4.1}
\end{equation}

\noindent where $\delta$ stands for the Dirac delta function. A
unique solution to {\rm (\ref{eq4.1})} which is bounded and
vanishes at $+\infty$ has a form
\begin{equation}\label{eq4.2}
h(t)=\left\{\begin{array}{l}
\sum_{j=1}^m \nu_j \,e^{-\alpha_j t},\qquad \qquad\qquad t>0,\\[4pt]
\sum_{j=1}^m \mu_j \,e^{\beta_j t},\qquad \qquad\qquad\,\,t<0.\\[4pt]
\end{array}
\right.
\end{equation}

\noindent Here $\alpha_j>0$, $j=1,2,...,m$, $\beta_j>0$ for
$j=2,...,m$ and $\beta_1=0$ are such that

\begin{equation}\label{eq4.3}
\{-\alpha_j\}_{j=1}^{m}\bigcup
\{\beta_j\}_{j=1}^{m}=\Bigl\{-m+\frac n2-\frac
12+2j\Bigr\}_{j=0}^{m-1}\cup \Bigl\{-\frac n2 -\frac
12+m-2j\Bigr\}_{j=0}^{m-1},
\end{equation}

\noindent and with the notation
\begin{equation}\label{eq4.4}
\vec\gamma=(-\alpha_1,...,-\alpha_m, \beta_1,...,\beta_m),\qquad
\vec\kappa=(\nu_1,...,\nu_m, -\mu_1,...,-\mu_m)
\end{equation}
\noindent the coefficients $\nu_j,\mu_j\in\RR$ satisfy
\begin{equation}\label{eq4.5}
\kappa_i=(-1)^{m+1}\Bigl(\prod_{j\neq i}(\gamma_j-\gamma_i)\Bigr)^{-1}.
\end{equation}
\end{lemma}

\vskip 0.08in

\bp Recall the representation formula \eqref{eq3.22} and let $p=0$. Similarly to the analysis in the proof of Proposition~\ref{p3.1}, one can see that the characteristic polynomial of the differential equation
(\ref{eq4.1}) has $2m$ distinct roots given by
\begin{equation} c_j:=-m+\frac n2-\frac 12+2j, \quad
j=0,1,...,m-1, \label{eq4.6}
\end{equation}

\noindent and

\begin{equation}\label{eq4.7}
-\frac n2 -\frac 12+m-2j=-1-\Bigl(-m+\frac n2-\frac
12+2j\Bigr)=-1-c_j,\quad j=0,1,...,m-1.
\end{equation}

\noindent Following \eqref{eq3.27}--\eqref{eq3.29}, \eqref{eq3.33}--\eqref{eq3.35} and \eqref{eq3.36}, we conclude that  there are $m$ distinct negative roots, 1
root equal to zero and $m-1$ distinct positive roots for the characteristic
polynomial. Hence, a bounded solution vanishing at infinity must
have the form (\ref{eq4.2}). In addition, $h$ should be such that
that $\partial_t^kg$ are continuous for $k=0,1,...,2m-2$ and
$\lim_{t\to 0^+}\partial_t^{2m-1} h(t)-\lim_{t\to
0^-}\partial_t^{2m-1} h(t)=(-1)^m$. This gives rise to the system of
equations
\begin{eqnarray}\label{eq4.8}
&& \sum_{j=1}^m \nu_j(-\alpha_j)^k-\sum_{j=1}^m \mu_j(\beta_j)^k=0,\quad k=0,1,...,2m-2,\\[4pt]\label{eq4.9}
&& \sum_{j=1}^m \nu_j(-\alpha_j)^{2m-1}-\sum_{j=1}^m
\mu_j(\beta_j)^{2m-1}=(-1)^m.
\end{eqnarray}

\noindent With the notation (\ref{eq4.4}) and the convention
$0^0=1$, (\ref{eq4.8})--(\ref{eq4.9}) can then be written as
$A\,\vec\kappa^\bot=(-1)^m(0,...,0,1)^\bot$, where $A$ is a matrix of
elements $A_{ji}=\gamma_i^{j-1}$, $i,j=1,...,2m$. Now using the
formula for the Vandermonde determinant, we find that the solution
is $\vec \kappa=(\kappa_j)_{j=1}^{2m}$,
\begin{equation}\label{eq4.10}
\kappa_j=(-1)^{j+m}\prod_{i<k,\, i\neq j,\, k\neq
j}(\gamma_j-\gamma_i)\,\,
\Bigl(\prod_{i<k}(\gamma_j-\gamma_i)\Bigr)^{-1}=(-1)^{m+1}\Bigl(\prod_{j\neq
i}(\gamma_j-\gamma_i)\Bigr)^{-1}.
\end{equation}
 \ep

\begin{theorem}\label{t4.2} Assume that $m\in\NN$ and  $n\in [3,2m+1]\cap \NN$
is odd. Let $\Omega$ be a bounded domain in $\RR^n$, $O\in
\RR^n\setminus\Omega$, $u\in C_0^\infty(\Omega)$ and
$v=e^{\left(m-\frac n2+\frac 12\right)t}(u\circ \varkappa^{-1})$. 
Then for
every $\xi\in\Omega$ and $\tau=\log |\xi|^{-1}$ 
\begin{equation}\label{eq4.11}
\int_{S^{n-1}}v^2(\tau,\omega)\,d\omega \leq  C
\int_{\RR^n}(-\Delta)^m u(x)\,u(x)g(\log |x|^{-1},\log |\xi|^{-1})\,dx,
\end{equation}

\noindent with
\begin{equation}\label{eq4.12}
g(t,\tau)=e^{t}\left(C_1 h(t-\tau)+C_2\right),\qquad t,\tau\in\RR,
\end{equation}

\noindent and $h$ given by Lemma~\ref{l4.1}. Here  $C,C_1,C_2$ are some constants depending on $m$ and $n$ only.
\end{theorem}

\bp {\bf Step I.}\,First of all, let us pass to $(t,\omega)$ coordinates. We aim to show that for some positive constants $C$ and $C'$
\begin{eqnarray}\label{eq4.13}\nonumber
\int_{S^{n-1}}v^2(\tau,\omega)\,d\omega &\leq & C
\int_{\RR}\int_{S^{n-1}} {\mathcal
L}^{m,n}(\partial_t,\delta)v(t,\omega)\,v(t,\omega)h(t-\tau)\,d\omega
dt\\[4pt]
&&\qquad +C' \int_{\RR}\int_{S^{n-1}} {\mathcal
L}^{m,n}(\partial_t,\delta)v(t,\omega)\,v(t,\omega)\,d\omega dt.
\end{eqnarray}

\noindent Clearly, \eqref{eq4.13} implies 
(\ref{eq4.11}). 

For future reference, observe that due to (\ref{eq4.1})
\begin{equation}\label{eq4.14}
{\mathcal L}^{m,n}(-\partial_t,0)h(t-\tau)=\delta(t-\tau),\qquad
t,\tau\in\RR.
\end{equation}

\vskip 0.08in
\noindent   {\bf Step II.}\, Going further, let us concentrate on the first integral on the right-hand side of \eqref{eq4.13}:
\begin{eqnarray}\label{eq4.15}
&& \int_{\RR}\int_{S^{n-1}} {\mathcal
L}^{m,n}(\partial_t,\delta)v(t,\omega)\,v(t,\omega)h(t-\tau)\,d\omega
dt\nonumber\\[4pt]
&&\qquad = \int_{\RR}\int_{S^{n-1}} {\mathcal
L}^{m,n}(0,\delta)v(t,\omega)\,v(t,\omega)h(t-\tau)\,d\omega
dt\nonumber\\[4pt]
&&\qquad\quad + \int_{\RR}\int_{S^{n-1}} \left({\mathcal
L}^{m,n}(\partial_t,\delta)-{\mathcal
L}^{m,n}(0,\delta)\right)v(t,\omega)\,v(t,\omega)h(t-\tau)\,d\omega
dt.
\end{eqnarray}

All terms of the operator $\left({\mathcal
L}^{m,n}(\partial_t,\delta)-{\mathcal
L}^{m,n}(0,\delta)\right)$  contain $\partial_t^k$ for some $k\geq 1$, and therefore, we can write

\begin{equation}\label{eq4.16}
{\mathcal L}^{m,n}(\partial_t,\delta)-{\mathcal
L}^{m,n}(0,\delta)=\sum_{\stackrel{k\geq 1,\,i\geq 0}{2i+k\leq
2m}}d_{ik}(-\delta)^i\partial_t^k,\qquad \mbox{for some}\quad d_{ik}\in\RR.
\end{equation}

\noindent  Hence,

\begin{eqnarray}\label{eq4.17}\nonumber
&& \int_{\RR}\int_{S^{n-1}} ({\mathcal
L}^{m,n}(\partial_t,\delta)-{\mathcal
L}^{m,n}(0,\delta))v(t,\omega)\,v(t,\omega)h(t-\tau)\,d\omega
dt\\[4pt]
&&\qquad = \sum_{\stackrel{k\geq 1,\,i\geq 0}{2i+k\leq
2m}}d_{ik}\int_{\RR}\int_{S^{n-1}}
\partial_t^k\nabla_\omega^i v(t,\omega)\,\nabla_\omega^i
v(t,\omega)h(t-\tau)\,d\omega dt.
\end{eqnarray}

We claim that

\begin{equation}\label{eq4.18}
\int_{\RR}\partial_t^kv\,vh\,dt=\frac 12
\int_{\RR}v^2\,(-\partial_t)^kh\,dt+\sum_{\stackrel{i\geq
1,\,j\geq 0}{2i+j\leq
k}}b_{ij}^k\int_{\RR}(\partial_t^iv)^2\,\partial_t^jh\,dt,\quad
b_{ij}^k\in\RR,
\end{equation}

\noindent for any $k\geq 1$. This can be proved by induction. For
$k=1$ and $k=2$ we have, respectively,

\begin{eqnarray}\nonumber
&& \int_{\RR}\partial_tv\,vh\,dt=\frac 12
\int_{\RR}v^2\,(-\partial_t)h\,dt,\nonumber\\[4pt]
&& \int_{\RR}\partial_t^2v\,vh\,dt=\frac 12
\int_{\RR}v^2\,(-\partial_t)^2h\,dt-\int_{\RR}(\partial_tv)^2\,h\,dt.\nonumber
\end{eqnarray}

\noindent Let us now assume that (\ref{eq4.18}) holds for
$k=1,2,...,l-1$, and prove it for $k=l$. First,
\begin{equation}\label{eq4.19}
\int_{\RR}\partial_t^lv\,vh\,dt=-
\int_{\RR}\partial_t^{l-2}(\partial_t v)\,(\partial_t v)h\,dt -
\int_{\RR}\partial_t^{l-1}v\,v\,\partial_th\,dt.
\end{equation}

\noindent Next, using (\ref{eq4.18}) for $l-1$ and $l-2$, we
deduce that the expression above is equal to

\begin{eqnarray}\label{eq4.20}\nonumber
&&-\frac 12
\int_{\RR}(\partial_tv)^2\,(-\partial_t)^{l-2}h\,dt-\sum_{\stackrel{i\geq
1,\,j\geq 0}{2i+j\leq
l-2}}b_{ij}^{l-2}\int_{\RR}(\partial_t^{i+1}v)^2\,\partial_t^jh\,dt\\[4pt]
&&\qquad +\frac 12
\int_{\RR}v^2\,(-\partial_t)^{l}h\,dt+\sum_{\stackrel{i\geq
1,\,j\geq 0}{2i+j\leq
l-1}}b_{ij}^{l-1}\int_{\RR}(\partial_t^iv)^2\,\partial_t^{j+1}h\,dt,
\end{eqnarray}

\noindent and (\ref{eq4.20}) can be written in the form
(\ref{eq4.18}) for $k=l$.

Then, using \eqref{eq4.18} for $v$ and $\nabla_\omega^i v$, \eqref{eq4.17} leads to the representation
\begin{eqnarray}\label{eq4.21}\nonumber
&& \int_{\RR}\int_{S^{n-1}} ({\mathcal
L}^{m,n}(\partial_t,\delta)-{\mathcal
L}^{m,n}(0,\delta))v(t,\omega)\,v(t,\omega)h(t-\tau)\,d\omega
dt\\[4pt]
&&\qquad = \frac 12 \sum_{\stackrel{k\geq 1,\,i\geq 0}{2i+k\leq
2m}}d_{ik}\int_{\RR}\int_{S^{n-1}}
(\nabla_\omega^i
v(t,\omega))^2(-\partial_t)^k h(t-\tau)\,d\omega dt\nonumber \\[4pt]
&&\qquad + \sum_{\stackrel{k\geq 1,\,i\geq 0}{2i+k\leq
2m}}\sum_{\stackrel{l\geq
1,\,j\geq 0}{2l+j\leq
k}}b_{lj}^{ki}\int_{\RR}\int_{S^{n-1}}
(\partial_t^l\nabla_\omega^iv)^2\,\partial_t^jh(t-\tau)\,d\omega dt=:I_1+I_2.
\end{eqnarray}

\vskip 0.08in
\noindent   {\bf Step II.}\,Consider $I_1$ first. Let us decompose $v$ into spherical harmonics. Then 
\begin{eqnarray}\label{eq4.22}\nonumber
&& I_1=\frac 12 \sum_{\stackrel{k\geq 1,\,i\geq 0}{2i+k\leq
2m}}d_{ik}\int_{\RR}\int_{S^{n-1}}
(\nabla_\omega^i
v(t,\omega))^2(-\partial_t)^k h(t-\tau)\,d\omega dt\nonumber \\[4pt]
&&\qquad =\frac 12 \sum_{\stackrel{k\geq 1,\,i\geq 0}{2i+k\leq
2m}}d_{ik}\sum_{p=0}^\infty \sum_{l=-p}^{p} \int_{\RR}  p^i(p+n-2)^i
v_{pl}^2(t)(-\partial_t)^k h(t-\tau)\,  dt
\nonumber \\[4pt]
&&\qquad =\frac 12\sum_{p=0}^\infty \sum_{l=-p}^{p} \int_{\RR} 
v_{pl}^2(t) \sum_{\stackrel{k\geq 1,\,i\geq 0}{2i+k\leq
2m}}d_{ik} p^i(p+n-2)^i (-\partial_t)^k h(t-\tau)\,  dt
\nonumber \\[4pt]
&&\qquad =\frac 12\sum_{p=0}^\infty \sum_{l=-p}^{p} \int_{\RR} 
v_{pl}^2(t) \left({\mathcal L}^{m,n}(-\partial_t,-p(p+n-2))-{\mathcal
L}^{m,n}(0,-p(p+n-2))\right) h(t-\tau)\,  dt, \end{eqnarray}

\noindent where we employed \eqref{eq4.16}. 

Now, 
\begin{eqnarray}\label{eq4.23} 
&&\frac 12\sum_{p=0}^\infty \sum_{l=-p}^{p} \int_{\RR} 
v_{pl}^2(t) {\mathcal
L}^{m,n}(0,-p(p+n-2))h(t-\tau)\,  dt\nonumber \\[4pt]
&&\qquad = \frac 12\int_{\RR}\int_{S^{n-1}}  {\mathcal
L}^{m,n}(0,\delta)
v(t,\omega)v(t,\omega)h(t-\tau)\,d\omega dt.
\end{eqnarray}

%Let us now separate the case $p=0$ from the rest of the sum above. Note that ${\mathcal L}^{m,n}(0,0)=0$. 

%\todo{Note to self: this is easy to see from \eqref{eq2.11.3.0}}

\noindent Also, taking into account \eqref{eq4.14}, we see that 

\begin{eqnarray}\label{eq4.24}\frac 12\sum_{p=0}^\infty \sum_{l=-p}^{p} \int_{\RR} 
v_{pl}^2(t) {\mathcal L}^{m,n}(-\partial_t,0) h(t-\tau)\,d\omega dt= \frac 12 \int_{S^{n-1}} 
v^2(\tau, \omega) \,d\omega.
 \end{eqnarray}

\noindent Hence, 
\begin{eqnarray}\label{eq4.25}\nonumber
&& I_1 =\frac 12\sum_{p=0}^\infty \sum_{l=-p}^{p} \int_{\RR} 
v_{pl}^2(t) \left({\mathcal L}^{m,n}(-\partial_t,-p(p+n-2))-{\mathcal L}^{m,n}(-\partial_t,0)\right) h(t-\tau)\,  dt \nonumber \\[4pt]
&&\qquad -\frac 12\int_{\RR}\int_{S^{n-1}}  {\mathcal
L}^{m,n}(0,\delta)
v(t,\omega)v(t,\omega)h(t-\tau)\,d\omega dt
\nonumber \\[4pt]
&&\qquad +\frac 12 \int_{S^{n-1}} 
v^2(\tau, \omega) \,d\omega. 
\end{eqnarray}

Let us discuss first in rough terms the bounds on \eqref{eq4.25} and then present the details.
The third term above isolates $v^2$ on the sphere, as desired (cf. \eqref{eq4.13}), and the second term will be combined with the first integral on the right-hand side of \eqref{eq4.15}, and ultimately estimated by \eqref{eq2.3}. 
As for the first term, we split it into the sum over $p\leq m-n/2+1/2$ and the complementary one. The sum over $p\geq m-n/2+3/2$ will be estimated by \eqref{eq2.3}. The sum over $p\leq m-n/2+1/2$ must be positive by itself. This amounts to showing that 
\begin{equation}\label{eq4.26}
\left({\mathcal L}^{m,n}(-\partial_t,-p(p+n-2))-{\mathcal L}^{m,n}(-\partial_t,0)\right) h(t)\geq 0, \quad \mbox{ for all } \quad t\neq 0,\quad 0\leq p\leq  m-\frac n2+\frac 12,
\end{equation}

\noindent or, equivalently (see \eqref{eq4.14}),
\begin{equation}\label{eq4.27}
{\mathcal L}^{m,n}(-\partial_t,-p(p+n-2)) h(t)\geq 0, \quad \mbox{ for all } \quad t\neq 0 \quad \mbox{and}\quad 0\leq p\leq  m-\frac n2+\frac 12.
\end{equation}

This, however, is exactly the result of Proposition~\ref{p3.1}.

\vskip 0.08in
\noindent   {\bf Step III.}\,Combining all of the above, we have 
\begin{eqnarray}\label{eq4.28}
&& \int_{\RR}\int_{S^{n-1}} {\mathcal
L}^{m,n}(\partial_t,\delta)v(t,\omega)\,v(t,\omega)h(t-\tau)\,d\omega
dt\nonumber\\[4pt]
&&\qquad = \frac 12 \int_{S^{n-1}} 
v^2(\tau, \omega) \,d\omega+\frac 12 \int_{\RR}\int_{S^{n-1}} {\mathcal
L}^{m,n}(0,\delta)v(t,\omega)\,v(t,\omega)h(t-\tau)\,d\omega
dt\nonumber\\[4pt]
&&\qquad\quad +\,\frac 12\sum_{p=0}^\infty \sum_{l=-p}^{p} \int_{\RR}
v_{pl}^2(t) \left({\mathcal L}^{m,n}(-\partial_t,-p(p+n-2))-{\mathcal L}^{m,n}(-\partial_t,0)\right) h(t-\tau)\, dt
\nonumber\\[4pt]
&&\qquad\quad + \sum_{\stackrel{k\geq 1,\,i\geq 0}{2i+k\leq
2m}}\sum_{\stackrel{l\geq
1,\,j\geq 0}{2l+j\leq
k}}b_{lj}^{ki}\int_{\RR}\int_{S^{n-1}}
(\partial_t^l\nabla_\omega^iv)^2\,\partial_t^jh(t-\tau)\,d\omega dt 
\nonumber\\[4pt]
&&\qquad \geq  \frac 12 \int_{S^{n-1}} 
v^2(\tau, \omega) \,d\omega+\frac 12 \int_{\RR}\int_{S^{n-1}} {\mathcal
L}^{m,n}(0,\delta)v(t,\omega)\,v(t,\omega)h(t-\tau)\,d\omega
dt
\nonumber\\[4pt]
&&\qquad\quad + \sum_{\stackrel{k\geq 1,\,i\geq 0}{2i+k\leq
2m}}\sum_{\stackrel{l\geq
1,\,j\geq 0}{2l+j\leq
k}}b_{lj}^{ki}\int_{\RR}\int_{S^{n-1}}
(\partial_t^l\nabla_\omega^iv)^2\,\partial_t^jh(t-\tau)\,d\omega dt\nonumber\\[4pt]
&&\qquad\quad +\,\frac 12\sum_{p=m-\frac n2+\frac 32}^\infty \sum_{l=-p}^{p} \int_{\RR}
v_{pl}^2(t) \left({\mathcal L}^{m,n}(-\partial_t,-p(p+n-2))-{\mathcal L}^{m,n}(-\partial_t,0)\right) h(t-\tau)\, dt.
\end{eqnarray}

\noindent Note that ${\mathcal
L}^{m,n}(0,\delta)$ is exactly the operator showing up in the last expression in \eqref{eq2.3}. Indeed,  
\begin{equation}\label{eq4.29}{\mathcal
L}^{m,n}(0,\delta)=(-1)^m\prod_{j=0}^{m-1}
\Bigg(\Bigl(m-\frac n2+\frac
12-2j\Bigr)\Bigl(m+\frac n2-\frac
32-2j\Bigr)+\delta_\omega\Bigg)=\prod_{p=-\frac n2+\frac 32}^{m-\frac
n2+\frac 12}\left(-\delta_\omega-p\,(p+n-2)\right).
\end{equation}

\noindent Then the absolute value of the second integral on the right-hand side of \eqref{eq4.28} can be majorized by the right-hand side of \eqref{eq2.3}. This is easy to see, for example, after the decomposition into spherical harmonics:
\begin{eqnarray}\label{eq4.30}\nonumber
&&\left|\sum_{q=0}^\infty \sum_{l=-q}^{q} \int_{\RR} {\mathcal
L}^{m,n}(0,-q\,(q+n-2)) v_{ql}(t) \,v_{ql}(t)h(t-\tau)\,dt\right|\\[4pt]\nonumber 
&&\qquad \leq  \sum_{q=0}^\infty \sum_{l=-q}^{q} \int_{\RR}\prod_{p=-\frac n2+\frac 32}^{m-\frac
n2+\frac 12}\left(q(q+n-2)-p\,(p+n-2)\right) (v_{ql}(t))^2\,|h(t-\tau)|\,dt\\[4pt]
&&\qquad \leq C \sum_{q=0}^\infty \sum_{l=-q}^{q} \int_{\RR}\prod_{p=-\frac n2+\frac 32}^{m-\frac
n2+\frac 12}\left(q(q+n-2)-p\,(p+n-2)\right) (v_{ql}(t))^2\,dt\nonumber \\[4pt]
&&\qquad =C \sum_{q=0}^\infty \sum_{l=-q}^{q} \int_{\RR} {\mathcal
L}^{m,n}(0,-q\,(q+n-2)) v_{ql}(t) \,v_{ql}(t)\,dt,
\end{eqnarray}

\noindent since $h$ is bounded by a constant and the polynomial 
\begin{equation}\label{eq4.31}\prod_{p=-\frac n2+\frac 32}^{m-\frac
n2+\frac 12}\left(q(q+n-2)-p\,(p+n-2)\right)\geq 0\quad\mbox{for all}\quad q\in\NN\cup\{0\}.
\end{equation}

Hence,
(\ref{eq4.11}) follows from (\ref{eq4.28}) once we prove that

\begin{eqnarray}\label{eq4.32}\nonumber
&&  \sum_{\stackrel{k\geq 1,\,i\geq 0}{2i+k\leq
2m}}\sum_{\stackrel{l\geq
1,\,j\geq 0}{2l+j\leq
k}}|b_{lj}^{ki}|\int_{\RR}\int_{S^{n-1}}
(\partial_t^l\nabla_\omega^iv)^2\,|\partial_t^jh(t-\tau)|\,d\omega dt\\[4pt]
&&\qquad\qquad\qquad\qquad\qquad\qquad \leq C
\int_{\RR^n}(-\Delta)^m u(x)\,u(x)|x|^{-1}\,dx
\end{eqnarray}

\noindent and 
\begin{eqnarray}\label{eq4.33}\nonumber
&&  \sum_{p=m-\frac n2+\frac 32}^\infty \sum_{l=-p}^{p} \int_{\RR} 
v_{pl}^2(t) \left|\left({\mathcal L}^{m,n}(-\partial_t,-p(p+n-2))-{\mathcal L}^{m,n}(-\partial_t,0)\right) h(t-\tau)\right|\,  dt\\[4pt]
&&\qquad\qquad\qquad\qquad\qquad\qquad \leq C
\int_{\RR^n}(-\Delta)^m u(x)\,u(x)|x|^{-1}\,dx
\end{eqnarray}

\noindent for some $C>0$. 

\vskip 0.08in
\noindent   {\bf Step IV.}\,
Since all derivatives of $h$ are bounded
(see (\ref{eq4.2})),

\begin{eqnarray}\label{eq4.34}\nonumber
&&  \sum_{\stackrel{k\geq 1,\,i\geq 0}{2i+k\leq
2m}}\sum_{\stackrel{l\geq
1,\,j\geq 0}{2l+j\leq
k}}|b_{lj}^{ki}|\int_{\RR}\int_{S^{n-1}}
(\partial_t^l\nabla_\omega^iv)^2\,|\partial_t^jh(t-\tau)|\,dt\\[4pt]
&&\qquad\leq C \sum_{\stackrel{k\geq 1,\,i\geq 0}{ i+k\leq
m}}\int_{\RR}\int_{S^{n-1}} (\partial_t^k\nabla_\omega^i
v(t,\omega))^2 \,d\omega dt.
\end{eqnarray}

\noindent Furthermore, for all
$k,i$ as above

\begin{eqnarray}\label{eq4.35}
&&\int_{\RR}\int_{S^{n-1}} (\partial_t^k\nabla_\omega^i
v(t,\omega))^2 \,d\omega dt =\sum_{p=0}^\infty \sum_{l=-p}^{p}
\Bigl(p\,(p+n-2)\Bigr)^i \int_{\RR}(\partial_t^k v_{pl}(t))^2 \,dt.
\end{eqnarray}

\noindent We break down the sum above into two parts,
corresponding to the cases $p\leq m-\frac n2+\frac 12$ and $p\geq m-\frac n2+\frac 32$, respectively. In the first case,
\begin{eqnarray}\label{eq4.36}
&&\sum_{p=0}^{m-\frac n2+\frac 12} \sum_{l=-p}^{p}
\Bigl(p\,(p+n-2)\Bigr)^i \int_{\RR}(\partial_t^k v_{pl}(t))^2
\,dt\leq C_{m,n} \sum_{p=0}^{m-\frac n2+\frac 12} \sum_{l=-p}^{p}
\int_{\RR}(\partial_t^k v_{pl}(t))^2 \,dt\nonumber\\[4pt]
&&\quad\leq C_{m,n}
\int_{\RR}\int_{S^{n-1}} (\partial_t^k v(t,\omega))^2 \,d\omega dt
 \leq C_{m,n} \int_{\RR^n}(-\Delta)^m u(x)\,u(x)|x|^{-1}\,dx,
\end{eqnarray}

\noindent where we used (\ref{eq2.3}) for the last inequality (note that $k\geq 1$ since we are in the range of \eqref{eq4.35}).

As for $p\geq m-\frac n2+\frac 32$, we claim that
\begin{equation}\label{eq4.37}
\Bigl(p\,(p+n-2)\Bigr)^m\leq C \prod_{s=-\frac n2+\frac
32}^{m-\frac n2+\frac
12}\Bigl(p\,(p+n-2)-s\,(s+n-2)\Bigr),\quad\mbox{ for every $p\geq
m-\frac n2+\frac 32$},
\end{equation}

\noindent where $C>0$ depends on $m$ and $n$ only. Indeed,  one can choose $C$ such that 
\begin{equation}\label{eq4.38}
\frac 1C\leq \left(1-\frac{\left(m-\frac n2+\frac 12\right)\left(m+\frac n2-\frac 32\right)}{\left(m-\frac n2+\frac 32\right)\left(m+\frac n2-\frac 12\right)}\right)^m\leq \prod_{s=-\frac n2+\frac 32}^{m-\frac n2+\frac
12}\left(1-\frac{s\,(s+n-2)}{p\,(p+n-2)}\right),
\end{equation}

\noindent for every $p\geq
m-\frac n2+\frac 32$.

However, since $i+k\leq m$ then by Young's inequality
\begin{eqnarray}\label{eq4.39}
&& \Bigl(p\,(p+n-2)\Bigr)^i \int_{\RR}(\partial_t^k v_{pl}(t))^2
\,dt= \int_{\RR}\Bigl(p\,(p+n-2)\Bigr)^i \gamma^{2k}
|\widehat{v_{pl}}(\gamma)|^2 \,d\gamma\nonumber\\[4pt] &&\qquad \leq
\int_{\RR}\Bigl(\Bigl(p\,(p+n-2)\Bigr)^m +\gamma^{2m}\Bigr)
|\widehat{v_{pl}}(\gamma)|^2 \,d\gamma\nonumber\\[4pt] &&\qquad \leq\Bigl(p\,(p+n-2)\Bigr)^m
\int_{\RR}(v_{pl}(t))^2 \,dt +\int_{\RR}(\partial_t^m v_{pl}(t))^2
\,dt.
\end{eqnarray}

\noindent Therefore,
\begin{eqnarray}\label{eq4.40}
&&\sum_{p=m-\frac n2+\frac 32}^{\infty} \sum_{l=-p}^{p} \Bigl(p\,(p+n-2)\Bigr)^i
\int_{\RR}(\partial_t^k v_{pl}(t))^2 \,dt \nonumber\\[4pt] &&\qquad \leq\sum_{p=m-\frac n2+\frac 32}^{\infty} \sum_{l=-p}^{p}
\,\,\prod_{s=-\frac n2+\frac 32}^{m-\frac n2+\frac
12}\Bigl(p\,(p+n-2)-s\,(s+n-2)\Bigr)
\int_{\RR}(v_{pl}(t))^2 \,dt \nonumber\\[4pt]&&\qquad\qquad +\sum_{p=m-\frac n2+\frac 32}^{\infty}
\sum_{l=-p}^{p}\int_{\RR}(\partial_t^m v_{pl}(t))^2 \,dt \leq C
\int_{\RR^n}(-\Delta)^m u(x)\,u(x)|x|^{-1}\,dx,
\end{eqnarray}

\noindent by Theorem~\ref{t2.1}.

Analogously, the left-hand side of \eqref{eq4.33} will give rise to terms 
\begin{eqnarray}\label{eq4.41}
&&\sum_{p=m-\frac n2+\frac 32}^{\infty} \sum_{l=-p}^{p} \Bigl(p\,(p+n-2)\Bigr)^i
\int_{\RR}v_{pl}^2(t) \,dt 
\end{eqnarray}

\noindent for $0\leq i\leq m$, which are bounded by

\begin{eqnarray}\label{eq4.42}
&&C\sum_{p=m-\frac n2+\frac 32}^{\infty} \sum_{l=-p}^{p}
\,\,\prod_{s=-\frac n2+\frac 32}^{m-\frac n2+\frac
12}\Bigl(p\,(p+n-2)-s\,(s+n-2)\Bigr)
\int_{\RR}(v_{pl}(t))^2 \,dt \nonumber\\[4pt]&&\qquad\qquad  \leq C
\int_{\RR^n}(-\Delta)^m u(x)\,u(x)|x|^{-1}\,dx,
\end{eqnarray}

\noindent as above.  \ep

\section{Integral identity and global estimate: the case of even dimension. Part I: power-logarithmic weight}\label{s5}
\setcounter{equation}{0}

\begin{theorem}\label{t5.1} Assume that $m\in\NN$ and $n\in [2,2m]\cap \NN$
is even. Let $\Omega$ be a bounded domain in $\RR^n$, $O\in
\RR^n\setminus\Omega$, $u\in C_0^\infty(\Omega)$ and
$v=e^{\left(m-\frac n2\right)t}(u\circ \varkappa^{-1})$. Furthermore, let $R$ be a positive constant such that the support of $u$ is contained in $B_{2R}$, $C_R:=\log (4R)$, and let $\psi$ be a weight function such that either $\psi(t)=C_R+t$ for all $t\in\RR$ or $\psi(t)=1$ for all $t\in\RR$.
Then whenever  $m$ is even,
\begin{eqnarray}\label{eq5.1}
\hskip -1cm &&\int_{\RR^n}(-\Delta)^m u(x)\,u(x)\,\psi(\log |x|^{-1})\,dx\geq C\sum_{k=1}^m \sum _{i=0}^{m-k}\int_{\RR}\int_{S^{n-1}}
\left(\partial_t^k \nabla_\omega^i v\right)^2\,\psi(t)\,d\omega dt\nonumber\\[4pt]
\hskip -1cm &&\quad +\,C
\int_{\RR}\int_{S^{n-1}} v\,\prod_{p}\left(-\delta_\omega-p\,(p+n-2)\right)^2 v\,\psi(t)\,d\omega dt,
\end{eqnarray}
\noindent where the product  is over $p=-n/2+2, -n/2+4,..., m-n/2-2, m-n/2$, that is, $p=-n/2+2j$ with $j=1,2,...,m/2$. If $m$ is odd,
\begin{eqnarray}\label{eq5.2}
\hskip -1cm &&\int_{\RR^n}(-\Delta)^m u(x)\,u(x)\,\psi(\log |x|^{-1})\,dx\geq C\sum_{k=1}^m \sum _{i=0}^{m-k}\int_{\RR}\int_{S^{n-1}}
\left(\partial_t^k \nabla_\omega^i v\right)^2\,\psi(t)\,d\omega dt\nonumber\\[4pt]
\hskip -1cm &&\quad +\,C
\int_{\RR}\int_{S^{n-1}} v\,\prod_{p}\left(-\delta_\omega-p\,(p+n-2)\right)^2 \left(-\delta_{\omega}+(n/2-1)^2\right) v\,\psi(t)\,d\omega dt,
\end{eqnarray}

\noindent where the product is over $p=-n/2+3, -n/2+5,..., m-n/2-2, m-n/2$, that is, $p=-n/2+1+2j$ with $j=1,2,...,(m-1)/2$. In both cases $C>0$ is some constant depending on $m$ and $n$
only.
\end{theorem}

\bp Passing to the coordinates $(t,\omega)$ and decomposing $v$ into spherical harmonics,
we see that
\begin{eqnarray}\label{eq5.3}
&&\int_{\RR^n}(-\Delta)^m
u(x)\,u(x)\,\psi(\log |x|^{-1})\,dx=\int_{\RR}\int_{S^{n-1}} {\mathcal
L}^{m,n}_o(\partial_t,\delta_\omega)v(t,\omega)\,v(t,\omega)\,\psi(t)\,d\omega dt\nonumber\\[4pt]
&&\qquad =\sum_{q=0}^{\infty}\sum_{l=-q}^q\int_{\RR}{\mathcal
L}^{m,n}_o(\partial_t,-q(q+n-2))v_{ql}(t)\,v_{ql}(t)\,\psi(t)\, dt,
\end{eqnarray}

\noindent where
\begin{eqnarray}\label{eq5.4}
{\mathcal L}^{m,n}_o(\partial_t,\delta_\omega)&=&(-1)^m\prod_{j=0}^{m-1}
\Bigg(\Bigl(-\partial_t+m-n/2-2j\Bigr)\Bigl(-\partial_t+m+n/2-2-2j\Bigr)+\delta_\omega\Bigg)\nonumber\\[4pt]
&=&\prod_{j=0}^{m-1}
\Bigg(-\Bigl(\partial_t-m+2j+1\Bigr)^2+\Bigl(n/2-1\Bigr)^2-\delta_\omega\Bigg).
\end{eqnarray}

\noindent Let us denote $A:=\Bigl(n/2-1\Bigr)^2-\delta_\omega$. Then the expression above is equal to 
\begin{eqnarray}\label{eq5.5}
&&\prod_{j=0}^{m-1}
\Bigl(\sqrt A -\partial_t+m-2j-1\Bigr)\Bigl(\sqrt A +\partial_t-m+2j+1\Bigr)\nonumber\\[4pt]
&&\qquad = \prod_{j=0}^{m-1}
\Bigl(\sqrt A -\partial_t+m-2j-1\Bigr)\Bigl(\sqrt A +\partial_t+m-2j-1\Bigr)
 = \prod_{j=0}^{m-1}
\Bigl(-\partial_t^2+(\sqrt A+m-2j-1)^2\Bigr).
\end{eqnarray}

\noindent Hence, 
\begin{equation}\label{eq5.6}
{\mathcal L}^{m,n}_o(\partial_t,-q(q+n-2))=\prod_{j=0}^{m-1}
\left(-\partial_t^2+{\mathcal B}_j(q)^2\right),
\end{equation}

\noindent with
\begin{equation}\label{eq5.7}
{\mathcal B}_j(q)^2=\Bigl(\sqrt{(n/2-1)^2+q(q+n-2)}+m-2j-1\Bigr)^2=\Bigl(q+n/2+m-2j-2\Bigr)^2,\quad q\in\NN\cup\{0\}.
\end{equation}

We claim that 
\begin{equation}\label{eq5.8}
{\mathcal B}_j(q)^2\geq C \,q(q+n-2), \mbox{when $q\in\NN\cup\{0\}$ is such that}\,\, q\neq 2j-m-n/2+2,
\end{equation}

\noindent for some $C>0$ depending on $m$ and $n$  only, and 
\begin{equation}\label{eq5.9}
{\mathcal B}_j(q)=0\quad \mbox{if}\quad  q=2j-m-n/2+2.
\end{equation}

\noindent Indeed, if $q>m-n/2$ then \eqref{eq5.8} holds with any $C\in (0,1)$ satisfying
\begin{equation}\label{eq5.10}
(1-C)\,q(q+n-2)>(m-n/2)(m+n/2-2) \quad\forall\,q>m-n/2, \,\,q\in\NN.
\end{equation}

\noindent Furthermore, one can directly check that ${\mathcal B}_j(q)=0$ for the values of $q$ described in \eqref{eq5.9} and otherwise ${\mathcal B}_j(q)^2>0$. Hence, if we choose $C$  to be the minimum of ${\mathcal B}_j(q)^2/\left(q(q+n-2)\right)$ over all $j=0,...,m-1,$ $0< q \leq m-n/2,$ such that $q\neq 2j-m-n/2+2$, we obtain \eqref{eq5.8} for the values of $q\leq m-n/2$.
This finishes the proof of (\ref{eq5.8}).

Now one can show that 
\begin{eqnarray}\label{eq5.11}
&&\int_{\RR}{\mathcal L}^{m,n}_o(\partial_t,-q(q+n-2))v_{ql}(t)\,v_{ql}(t)\,\psi(t)\, dt\nonumber\\[4pt]
&&\qquad =\sum_{j=1}^m\sum_{k_1<...<k_j}\int_{\RR}(-\partial_t^2)^{m-j}{\mathcal B}_{k_1}(q)^2...{\mathcal B}_{k_j}(q)^2v_{ql}(t)\,v_{ql}(t)\,\psi(t)\, dt+\int_{\RR}(-\partial_t^2)^{m}v_{ql}(t)\,v_{ql}(t)\,\psi(t)\, dt\nonumber\\[4pt]
&&\qquad =\sum_{j=1}^m\sum_{k_1<...<k_j}\int_{\RR}{\mathcal B}_{k_1}(q)^2...{\mathcal B}_{k_j}(q)^2 (\partial_t^{m-j}v_{ql}(t))^2\,\psi(t)\, dt+\int_{\RR}(-\partial_t^2)^{m}v_{ql}(t)\,v_{ql}(t)\,\psi(t)\, dt\nonumber\\[4pt]
&&\qquad\geq C\sum_{\stackrel{i\geq 0,\,\,k\geq 1}{i+k=m}}\int_{\RR}\left(q(q+n-2)\right)^i (\partial_t^{k}v_{ql}(t))^2\,\psi(t)\, dt+ \int_{\RR}\prod_{j=0}^{m-1}{\mathcal B}_{j}(q)^2v_{ql}^2(t)\,\psi(t)\, dt\nonumber\\[4pt]
&&\qquad\geq C\sum_{\stackrel{i\geq 0,\,\,k\geq 1}{i+k\leq m}}\int_{\RR}\left(q(q+n-2)\right)^i (\partial_t^{k}v_{ql}(t))^2\,\psi(t)\, dt+ \int_{\RR}{\mathcal L}^{m,n}_o(0,-q(q+n-2))v_{ql}^2(t)\,\psi(t)\, dt.
\end{eqnarray}

\noindent Note that the second equality uses the particular form of the weight function. Indeed, when $\psi\equiv 1$, this is just integration by parts, and when $\psi(t)=C_R+t$, $t\in\RR$, we prove by induction in $k\in\NN\cup\{0\}$ that 
\begin{eqnarray}\label{eq5.12}
\int_{\RR}(-\partial_t^2)^{k}v^2(t)\,(C_R+t)\, dt =\int_{\RR}(\partial_t^{k}v)^2(t)\,(C_R+t)\, dt,\qquad k\in\NN.
\end{eqnarray}

Finally, a straightforward calculation shows that whenever $m$ is even,
\begin{equation}\label{eq5.13}
{\mathcal L}^{m,n}_o(0,-q(q+n-2))=\prod_{p}\left(q(q+n-2)-p\,(p+n-2)\right)^2,\quad q\in\NN\cup\{0\},
\end{equation}

\noindent where the product above is over $p=-n/2+2, -n/2+4,..., m-n/2-2, m-n/2$, and if $m$ is odd,
\begin{equation}\label{eq5.14}
{\mathcal L}^{m,n}_o(0,-q(q+n-2))=\prod_{p}\left(q(q+n-2)-p\,(p+n-2)\right)^2 \left(q(q+n-2)+(n/2-1)^2\right),\end{equation}

\noindent  for every $q\in\NN\cup\{0\}$, with the product above over $p=-n/2+3, -n/2+5,..., m-n/2-2, m-n/2$. This finishes the argument. \ep

\section{Integral identity and global estimate: the case of even dimension. Part II: weight $g$}\label{s6}
\setcounter{equation}{0}

\begin{lemma}\label{l6.1} Assume that $m\in\NN$, $n\in [2,2m]\cap \NN$
is even and $m-n/2$ is even. Recall that in this case
\begin{equation}
\Lmno(-\partial_t,0)=\prod_{j=0}^{m-1} \left(-\partial_t^2+\Bigl(m-\frac n2-2j\Bigr)^2\right)
 \label{eq6.1}
\end{equation}

\noindent and consider the equation
\begin{equation}
\prod_{j=0}^{m-1} \left(-\partial_t^2+\Bigl(m-\frac n2-2j\Bigr)^2\right)h
 =\delta,  \label{eq6.2}
\end{equation}

\noindent where $\delta$ stands for the Dirac delta function. A
unique solution to {\rm (\ref{eq6.2})} which 
vanishes at $+\infty$ and has at most linear growth or decay at $-\infty$ has a form
\begin{equation}\label{eq6.3}
h(t)=\left\{\begin{array}{l}
\sum_{i=1}^{(m-n/2)/2} \nu_i^{(1)} \,e^{-2i t}+\sum_{i=1}^{(m-n/2)/2} \nu_i^{(2)} \,te^{-2i t}  + \sum_{i=1}^{n/2-1} \nu_i^{(3)} \,e^{-(m-n/2+2i) t},\qquad \qquad\,\, t>0,\\[6pt]
\sum_{i=1}^{(m-n/2)/2} \mu_i^{(1)} \,e^{2i t}+\sum_{i=1}^{(m-n/2)/2} \mu_i^{(2)} \,te^{2i t}  + \sum_{i=1}^{n/2-1} \mu_i^{(3)} \,e^{(m-n/2+2i) t}+\mu^{(4)}t+\mu^{(5)},\,\, t<0.\\[6pt]
\end{array}
\right.
\end{equation}

\noindent Here $\nu_i^{(1)}, \nu_i^{(2)}, \mu_i^{(1)} , \mu_i^{(2)}$, $i=1,...,(m-n/2)/2$, $\nu_i^{(3)}, \mu_i^{(3)}$, $i=1,...,n/2-1$, and $\mu^{(4)}, \mu^{(5)}$ are some real numbers depending on $m$ and $n$ only.

\end{lemma}

\vskip 0.08in

\bp The characteristic polynomial of the differential equation
(\ref{eq6.2}) has $n/2-1$ single roots given by 
\begin{equation} \mbox{even numbers\quad  from }\quad -m-n/2+2\quad \mbox{ to }\quad -m+n/2-2, \label{eq6.4}
\end{equation}

\noindent then $m-n/2+1$ double roots given by 
\begin{equation}\label{eq6.5}
\mbox{even numbers\quad  from }\quad -m+n/2\quad \mbox{ to }\quad m-n/2, 
\end{equation}

\noindent (including 0), and then another $n/2-1$ single roots given by 
\begin{equation} \mbox{even numbers\quad  from }\quad m-n/2+2\quad \mbox{ to }\quad m+n/2-2. \label{eq6.6}
\end{equation}

\noindent This determines the structure of the solution (\ref{eq6.3}).
Furthermore, the solution to (\ref{eq6.2})  $h$ should be such that
that $\partial_t^kh$ are continuous for $k=0,1,...,2m-2$, i.e., $\lim_{t\to 0^+}\partial_t^{k} h(t)-\lim_{t\to
0^-}\partial_t^{k} h(t)=0$ for $k=0,1,...,2m-2$, and
$\lim_{t\to 0^+}\partial_t^{2m-1} g(t)-\lim_{t\to
0^-}\partial_t^{2m-1} h(t)=(-1)^m$. Hence, the coefficients in (\ref{eq6.3}) satisfy the linear system 
\begin{equation}\label{eq6.7}
A\,\vec\kappa^\bot=(-1)^m(0,...,0,1)^\bot,
\end{equation}
\noindent where the vector on the right-hand side is of the length $2m$ (it has first $2m-1$ entries equal to 0 and the last entry equal to 1), the vector $\vec\kappa$ is also of the length $2m$ with the entries (in order) $\nu_i^{(1)}, \nu_i^{(2)}$,  $i=1,...,(m-n/2)/2$,  $\nu_i^{(3)}$, $i=1,...,n/2-1$, and $-\mu_i^{(1)} , -\mu_i^{(2)}$, $i=1,...,(m-n/2)/2$, $-\mu_i^{(3)}$, $i=1,...,n/2-1$, $-\mu^{(4)}, -\mu^{(5)}$. Finally, $A$ is a $2m\times 2m$ matrix built as follows: 
\begin{eqnarray}\nonumber
&& A_{jk}=(-2i)^{j-1},\,\,\mbox{for}\,\,i=k,\,\,j=1,...,2m,\,\,k= 1,...,(m-n/2)/2, \\[4pt]\nonumber
&& A_{0k}=0,\,\,\mbox{for}\,\,k= (m-n/2)/2+1,...,m-n/2, \\[4pt]\nonumber
&& A_{jk}=(j-1)(-2i)^{j-2},\,\,\mbox{for}\,\, i=k-(m-n/2)/2,\,\,j=2,...,2m,\,\,k= (m-n/2)/2+1,...,m-n/2, \\[4pt]\nonumber
&& A_{jk}=(-m+n/2-2i)^{j-1},\,\,\mbox{for}\,\,i=k-m+n/2,\,\,j=1,...,2m,\,\,k= m-n/2+1,...,m-1, \\[4pt]\nonumber
&& A_{jk}=(2i)^{j-1},\,\,\mbox{for}\,\,i=k-m+1,\,\,j=1,...,2m,\,\,k= m,...,(3m-n/2)/2-1, \\[4pt]\nonumber
&& A_{0k}=0,\,\,\mbox{for}\,\,k= (3m-n/2)/2,...,2m-n/2-1, \\[4pt]\nonumber
&& A_{jk}=(j-1)(2i)^{j-2},\,\,\mbox{for}\,\, i=k-(3m-n/2)/2+1,\,\,j=2,...,2m,\,\,k= (3m-n/2)/2,...,2m-n/2-1, \\[4pt]\nonumber
&& A_{jk}=(m-n/2+2i)^{j-1},\,\,\mbox{for}\,\,i=k-2m+n/2+1,\,\,j=1,...,2m,\,\,k= 2m-n/2,...,2m-2, \\[4pt]
\nonumber
&& A_{jk}=0,\,\,\mbox{for}\,\,j=1\,\,\mbox{and}\,\,j=3,...,2m,\,\,\mbox{and}\,\,A_{2k}=1,\,\,\mbox{when}\,\,k= 2m-1,\\[4pt] 
\nonumber
&& A_{jk}=0,\,\,\mbox{for}\,\,j=1,...,2m,\,\,\mbox{and}\,\,A_{1k}=1,\,\,\mbox{when}\,\,k= 2m.
\nonumber
\end{eqnarray}

It remains to show that the system (\ref{eq6.7}) has a unique solution, i.e., that the determinant of $A$ is not 0. This, however, is a standard argument in the theory of ordinary differential equations since determinant of $A$ is Wronskian at $0$ for a complete system of solutions for an ODE of order $2m$. We omit the details. \ep

\comment{
\todo{Note to self: I found a complete and well-written proof of the fact that the determinant of the matrix $A$ is not 0  in Bruce Shapiro's lecture notes on Differential equations on his website: see

http://bruce-shapiro.com/math351/351notes/9.HOConstantCoefficients.pdf

or 

http://bruce-shapiro.com/math351/351notes/

It should be in some other ODE books too, but I do not have another clean reference.}
}

\begin{lemma}\label{l6.2} Assume that $m\in\NN$, $n\in [2,2m]\cap \NN$
is even and $m-n/2$ is odd. Recall that in this case
\begin{equation}
\Lmno(-\partial_t,1-n)=\Lmno(-\partial_t, -1(1+n-2))=\prod_{j=0}^{m-1} \left(-\partial_t^2+\Bigl(m+\frac n2-2j-1\Bigr)^2\right)
 \label{eq6.8}
\end{equation}

\noindent and consider the equation
\begin{equation}
\prod_{j=0}^{m-1} \left(-\partial_t^2+\Bigl(m+\frac n2-2j-1\Bigr)^2\right)h
 =\delta,  \label{eq6.9}
\end{equation}

\noindent where $\delta$ stands for the Dirac delta function. A
unique solution to {\rm (\ref{eq6.9})} which 
vanishes at $+\infty$ and has at most linear growth or decay at $-\infty$ has a form
\begin{equation}\label{eq6.10}
h(t)=\left\{\begin{array}{l}
\sum_{i=1}^{(m-n/2-1)/2} \nu_i^{(1)} \,e^{-2i t}+\sum_{i=1}^{(m-n/2-1)/2} \nu_i^{(2)} \,te^{-2i t}  + \sum_{i=1}^{n/2} \nu_i^{(3)} \,e^{-(m-n/2-1+2i) t},\qquad \qquad\,\, t>0,\\[6pt]
\sum_{i=1}^{(m-n/2-1)/2} \mu_i^{(1)} \,e^{2i t}+\sum_{i=1}^{(m-n/2-1)/2} \mu_i^{(2)} \,te^{2i t}  + \sum_{i=1}^{n/2} \mu_i^{(3)} \,e^{(m-n/2-1+2i) t}+\mu^{(4)}t+\mu^{(5)},\,\, t<0.\\[6pt]
\end{array}
\right.
\end{equation}

\noindent Here $\nu_i^{(1)}, \nu_i^{(2)}, \mu_i^{(1)} , \mu_i^{(2)}$, $i=1,...,(m-n/2-1)/2$, $\nu_i^{(3)}, \mu_i^{(3)}$, $i=1,...,n/2$, and $\mu^{(4)}, \mu^{(5)}$ are some real numbers depending on $m$ and $n$ only.

\end{lemma}

We would like to remark a curious difference in $h$ for three considered situations: $n$ odd; $n$ even and $m-n/2$ even; $n$ even and $m-n/2$ odd.
While formula \eqref{eq6.3} is rather different from \eqref{eq4.2}, the underlying differential equation in \eqref{eq6.1}--\eqref{eq6.2} is a quite natural analogue of \eqref{eq4.1}. Lemma~\ref{l6.2}, however, deals with a different equation \eqref{eq6.8}. Just as in Section~\ref{s4}, we will show below how one can build a suitable weight from the corresponding function $h$.

\vskip 0.08in

\bp The characteristic polynomial of the differential equation
(\ref{eq6.9}) has $n/2$ single roots given by 
\begin{equation} \mbox{even numbers\quad  from }\quad -m-n/2+1\quad \mbox{ to }\quad -m+n/2-1, \label{eq6.11}
\end{equation}

\noindent then $m-n/2$ double roots given by 
\begin{equation}\label{eq6.12}
\mbox{even numbers\quad  from }\quad -m+n/2+1\quad \mbox{ to }\quad m-n/2-1, 
\end{equation}

\noindent (including 0), and then another $n/2$ single roots given by 
\begin{equation} \mbox{even numbers\quad  from }\quad m-n/2+1\quad \mbox{ to }\quad m+n/2-1. \label{eq6.13}
\end{equation}

In fact, these are the same collections of numbers as those appearing in \eqref{eq6.4}--\eqref{eq6.6}, in the following sense. Let us denote the sequence of numbers described in \eqref{eq6.4}--\eqref{eq6.6} by $S_{m,n}^e$ and the one described by \eqref{eq6.11}--\eqref{eq6.13} by $S_{m,n}^o$. Then 
$$ \left\{S_{m,n}^e\right\}_{m\in\NN,\,n\in [2,2m]\cap \NN,\,m-n/2\,\textrm{even}}\supset \left\{S_{m,n}^o\right\}_{m\in\NN,\,n\in [2,2m]\cap \NN,\,m-n/2\,\textrm{odd}},$$
\noindent via the identification $S_{m,n}^o=S_{m,n+2}^e$ for $m, n$ as above. This identification immediately leads to the observation that the solution to \eqref{eq6.9} must be given by the formula \eqref{eq6.3} with $n$ replaced by $n+2$. \ep

\begin{theorem}\label{t6.3} Assume that $m\in\NN$ and  $n\in [2,2m]\cap \NN$
is even. Let $\Omega$ be a bounded domain in $\RR^n$, $O\in
\RR^n\setminus\Omega$, $u\in C_0^\infty(\Omega)$ and
$v=e^{\left(m-\frac n2\right)t}(u\circ \varkappa^{-1})$. Let $R$ be a positive constant such that the support of $u$ is contained in $B_{2R}$.
Then there exist positive constants $C,$ $C'$, $C''$, depending on $m$ and $n$ only, such that for
every $\xi\in B_{2R}$ and $\tau=\log |\xi|^{-1}$ we have
\begin{eqnarray}\label{eq6.14}
\int_{S^{n-1}}v^2(\tau,\omega)\,d\omega \leq  C
\int_{\RR^n}(-\Delta)^m u(x)\,u(x)g(\log
|x|^{-1}, \log |\xi|^{-1})\,dx
\end{eqnarray}

\noindent  where $C_R=\log (4R)$ and $g$ is defined by 
\begin{equation}\label{eq6.15}
g(t,\tau) = h(t-\tau)+\mu^{(4)}(C_R+\tau)+C'+C''(C_R+t),
\end{equation}
where $h$ and $\mu^{(4)}$ are given by \eqref{eq6.3} and \eqref{eq6.10}, depending on the parity of $m-\frac n2.$

\end{theorem} 

\bp {\bf Step I, the preliminary choice of weight}. Recall the function $h$ defined in \eqref{eq6.3} and \eqref{eq6.10} for the cases when $m-n/2$ is even and odd, respectively. Let 
\begin{equation}\label{eq6.16}
\widetilde g (t,\tau):=h(t-\tau)+\mu^{(4)}(C_R+\tau).\qquad t,\tau\in\RR.
\end{equation}

\noindent In particular, when $m-n/2$ is even, we have 
\begin{equation}\label{eq6.17}
\widetilde g (t,\tau)=\left\{\begin{array}{l}
\sum_{i=1}^{(m-n/2)/2} \nu_i^{(1)} \,e^{-2i (t-\tau)}+\sum_{i=1}^{(m-n/2)/2} \nu_i^{(2)} \,(t-\tau)e^{-2i (t-\tau)}  \\[4pt]\qquad\qquad\qquad + \sum_{i=1}^{n/2-1} \nu_i^{(3)} \,e^{-(m-n/2+2i) (t-\tau)}+\mu^{(4)}(C_R+\tau),\qquad\,\, t>\tau,\\[8pt]
\sum_{i=1}^{(m-n/2)/2} \mu_i^{(1)} \,e^{2i (t-\tau)}+\sum_{i=1}^{(m-n/2)/2} \mu_i^{(2)} \,(t-\tau)e^{2i (t-\tau)}  \\[4pt]\qquad\qquad\qquad  + \sum_{i=1}^{n/2-1} \mu_i^{(3)} \,e^{(m-n/2+2i) (t-\tau)}+\mu^{(4)}(C_R+t)+\mu^{(5)},\,\,\qquad  t<\tau,\\[6pt]
\end{array}
\right.
\end{equation}

\noindent and in the case when $m-n/2$ is odd, we have
\begin{equation}\label{eq6.18}
\widetilde g (t,\tau)=\left\{\begin{array}{l}
\sum_{i=1}^{(m-n/2-1)/2} \nu_i^{(1)} \,e^{-2i (t-\tau)}+\sum_{i=1}^{(m-n/2-1)/2} \nu_i^{(2)} \,(t-\tau)e^{-2i (t-\tau)}  \\[4pt]\qquad\qquad\qquad + \sum_{i=1}^{n/2} \nu_i^{(3)} \,e^{-(m-n/2-1+2i) (t-\tau)}+\mu^{(4)}(C_R+\tau),\qquad\,\, t>\tau,\\[8pt]
\sum_{i=1}^{(m-n/2-1)/2} \mu_i^{(1)} \,e^{2i (t-\tau)}+\sum_{i=1}^{(m-n/2-1)/2} \mu_i^{(2)} \,(t-\tau)e^{2i (t-\tau)}  \\[4pt]\qquad\qquad\qquad  + \sum_{i=1}^{n/2} \mu_i^{(3)} \,e^{(m-n/2-1+2i) (t-\tau)}+\mu^{(4)}(C_R+t)+\mu^{(5)},\,\,\qquad  t<\tau.\\[6pt]
\end{array}
\right.
\end{equation}

For future reference, we record a few estimates. First, by our assumptions on $\xi$ and the support of $u$, the discussion will be naturally restricted to the case $t,\tau \geq \log (2R)^{-1}$. Hence, both $C_R+t$ and $C_R+\tau$ are positive. We remark that we do not claim positivity of $h$ or that of $\widetilde g$. However, 
\begin{equation}\label{eq6.19} 
|\widetilde g(t,\tau)|\leq C_0(m,n)+|\mu^{(4)}|\,(C_R+t), \quad t,\tau\geq \log (2R)^{-1}, \end{equation}

\noindent and 
\begin{equation}\label{eq6.20} 
|\partial_t^l\,\widetilde g(t,\tau)|\leq C_1(m,n),\quad t,\tau\geq \log (2R)^{-1}, \quad 1\leq l\leq 2m,
\end{equation}

\noindent for some constants $C_0(m,n), C_1(m,n)>0$ depending on $m,n$ only. We note that $\partial_t^l\,\widetilde g$ can be defined at $t=\tau$ for all $l<2m$ by continuity, and for $l=2m$ one assumes $t\neq \tau$ in \eqref{eq6.20}.

\vskip 0.08 in \noindent {\bf Step II, the set-up}. We commence similarly to \eqref{eq4.15}:
\begin{eqnarray}\label{eq6.21}
&&
\int_{\RR^n}(-\Delta)^m u(x)\,u(x)\widetilde g(\log
|x|^{-1},\log |\xi|^{-1})\,dx
= \int_{\RR}\int_{S^{n-1}} {\mathcal
L}^{m,n}_o(\partial_t,\delta)v(t,\omega)\,v(t,\omega)\widetilde g(t,\tau)\,d\omega
dt\nonumber\\[4pt]
&&\qquad = \int_{\RR}\int_{S^{n-1}} {\mathcal
L}^{m,n}_o(0,\delta)v(t,\omega)\,v(t,\omega)\widetilde g(t,\tau)\,d\omega
dt\nonumber\\[4pt]
&&\qquad\quad + \int_{\RR}\int_{S^{n-1}} \left({\mathcal
L}^{m,n}_o(\partial_t,\delta)-{\mathcal
L}^{m,n}_o(0,\delta)\right)v(t,\omega)\,v(t,\omega)\widetilde g(t,\tau)\,d\omega
dt.
\end{eqnarray}

\noindent Furthermore, as in \eqref{eq4.21}--\eqref{eq4.22},
\begin{eqnarray}\label{eq6.22}\nonumber
&& \int_{\RR}\int_{S^{n-1}} ({\mathcal
L}^{m,n}_o(\partial_t,\delta)-{\mathcal
L}^{m,n}_o(0,\delta))v(t,\omega)\,v(t,\omega)\widetilde g(t,\tau)\,d\omega
dt\\[4pt]
&&\qquad = \frac 12\sum_{p=0}^\infty \sum_{l=-p}^{p} \int_{\RR} 
v_{pl}^2(t) \left({\mathcal L}^{m,n}(-\partial_t,-p(p+n-2))-{\mathcal
L}^{m,n}(0,-p(p+n-2))\right) \widetilde g(t,\tau)\,  dt\nonumber \\[4pt]
&&\qquad + \sum_{\stackrel{k\geq 1,\,i\geq 0}{2i+k\leq
2m}}\sum_{\stackrel{l\geq
1,\,j\geq 0}{2l+j\leq
k}}c_{lj}^{ki}\int_{\RR}\int_{S^{n-1}}
(\partial_t^l\nabla_\omega^iv)^2\,\partial_t^j\,\widetilde g(t,\tau)\,d\omega dt=:I_1^o+I_2^o,
\end{eqnarray}

\noindent for some constants $c_{lj}^{ki}$ depending on $m,n$ only. 

Let us now denote 
\begin{equation}\label{eq6.23}
p_0:=\left\{\begin{array}{l}
0,\qquad\,\, \mbox{when } m-n/2 \mbox{ is even},\\[8pt]
1,\qquad\,\, \mbox{when } m-n/2 \mbox{ is odd}.\\[6pt]
\end{array}
\right.
\end{equation}
 
Then 
\begin{eqnarray}\label{eq6.24}\nonumber
 I_1^o &=&\frac 12\sum_{p=0}^\infty \sum_{l=-p}^{p} \int_{\RR} 
v_{pl}^2(t) \Bigg({\mathcal L}^{m,n}_o(-\partial_t,-p(p+n-2))-{\mathcal
L}^{m,n}_o(0,-p(p+n-2))\nonumber \\[4pt] && \qquad\qquad\qquad \qquad\qquad
-{\mathcal L}^{m,n}_o(-\partial_t,-p_0(p_0+n-2)) \Bigg) \,\,\widetilde g(t,\tau)\,  dt 
\nonumber \\[4pt]
&& +\frac 12 \int_{S^{n-1}} 
v^2(\tau, \omega) \,d\omega, 
\end{eqnarray}

\noindent where we used \eqref{eq6.1}--\eqref{eq6.2}, \eqref{eq6.8}--\eqref{eq6.9}, as well as the fact that the operator ${\mathcal L}_o^{m,n}(-\partial_t,-p_0(p_0+n-2))$ kills constants and thus, 
\begin{equation}\label{eq6.25} {\mathcal
L}^{m,n}_o(-\partial_t,-p_0(p_0+n-2)) \widetilde g(t,\tau)={\mathcal
L}^{m,n}_o(-\partial_t,-p_0(p_0+n-2))h(t-\tau)=\delta(t-\tau).
\end{equation}

\vskip 0.08 in \noindent {\bf Step III}. The term in the representation of $I_1^o$ in \eqref{eq6.24} associated to ${\mathcal
L}^{m,n}_o(0,-p(p+n-2))$ and the first term on the right-hand side of \eqref{eq6.21} are, modulo a multiplicative constant, the same, and can be estimated as follows:
\begin{eqnarray}\label{eq6.26} 
&&\left|\int_{\RR}\int_{S^{n-1}}  {\mathcal
L}^{m,n}_o(0,\delta)
v(t,\omega)v(t,\omega)\widetilde g(t,\tau)\,d\omega dt\right|\nonumber \\[4pt]
&&\quad=\left|\sum_{p=0}^\infty \sum_{l=-p}^{p} \int_{\RR} 
v_{pl}^2(t) {\mathcal
L}^{m,n}_o(0,-p(p+n-2))\widetilde g(t,\tau)\,  dt\right|\nonumber \\[4pt]
&&\quad\leq \sum_{p=0}^\infty \sum_{l=-p}^{p} \int_{\RR} 
v_{pl}^2(t) {\mathcal
L}^{m,n}_o(0,-p(p+n-2))\,|\widetilde g(t,\tau)|\,  dt\nonumber \\[4pt]
&&\quad\leq \sum_{p=0}^\infty \sum_{l=-p}^{p} \int_{\RR} 
v_{pl}^2(t) {\mathcal
L}^{m,n}_o(0,-p(p+n-2))\,\left(C_0(m,n)+|\mu^{(4)}|\,(C_R+t)\right)\,  dt\nonumber \\[4pt]
&&\quad=\int_{\RR}\int_{S^{n-1}}  {\mathcal
L}^{m,n}_o(0,\delta)
v(t,\omega)v(t,\omega)\,\left(C_0(m,n)+|\mu^{(4)}|\,(C_R+t)\right)\,\,d\omega dt\nonumber \\[4pt]
&&\quad\leq C\int_{\RR^n}(-\Delta)^m u(x)\,u(x)\,\left(C_0(m,n)+|\mu^{(4)}|\,(C_R+\log |x|^{-1})\right) \,dx,
\end{eqnarray}

\noindent where we used the fact that the polynomial ${\mathcal
L}^{m,n}(0,-p(p+n-2))\geq 0$ (see \eqref{eq5.13}--\eqref{eq5.14}) for the first inequality, \eqref{eq6.19} for the second one, and 
\eqref{eq5.1}--\eqref{eq5.2} with \eqref{eq5.13}--\eqref{eq5.14} for the last inequality above. 

It remains to estimate $I_2^o$ (see \eqref{eq6.22} for definition) and 
\begin{equation}\label{eq6.27}\sum_{p=0}^\infty \sum_{l=-p}^{p} \int_{\RR} 
v_{pl}^2(t) \Bigg({\mathcal L}^{m,n}_o(-\partial_t,-p(p+n-2))
-{\mathcal L}^{m,n}_o(-\partial_t,-p_0(p_0+n-2)) \Bigg) \,\,\widetilde g(t,\tau)\,  dt. \end{equation}

\noindent For the terms in the expression above (and for $I_2^o$), it is sufficient to either prove positivity or, for those which are not necessarily positive,  a bound on (the sum of) their absolute values in terms of the right-hand side of \eqref{eq6.14}.

\vskip 0.08 in \noindent {\bf Step IV, $0\leq p\leq m-n/2$, $m-n/2$ even}. At this stage, the discussion splits according to whether $m-n/2$ is even or odd. In fact, the two cases are, in some sense, symmetric, but tracking both at the same time would make the discussion too cumbersome. We start with the situation when $m-n/2$ is even, that is, we are in the regimen of Lemma~\ref{l6.1}, and then we will list the necessary modifications for $m-n/2$ odd.

To outline in a few words the strategy, we are going to show that the terms in the expression above corresponding to $0\leq p\leq m-n/2$ even are positive by virtue of the theory developed in Section~\ref{s3}, while the terms corresponding to $0\leq p\leq m-n/2$ odd, as well as the sum of those corresponding to $p>m-n/2$, are bounded by the right-hand side of \eqref{eq6.14}. Essentially, for large $p$ we use the same considerations as in Section~\ref{s4}, while for $0\leq p\leq m-n/2$ odd we observe that the last integrals in \eqref{eq5.1}, \eqref{eq5.2} do not vanish and provide the desired bounds. Finally, $I_2^o$ can be analyzed much as in Section~\ref{s4}. Now we turn to the details.

\vskip 0.08 in \noindent {\bf Step IV(a), $m-n/2$ even, $0\leq p\leq m-n/2$ even}. In this case, we shall employ the strategy suggested by the results of Section~\ref{s3}. One would like to show that 
\begin{multline}\label{eq6.28}{\mathcal L}^{m,n}_o(-\partial_t,0)h=\delta \quad \Longrightarrow \\[4pt]\quad {\mathcal L}^{m,n}_o(-\partial_t,-p(p+n-2))h(t)\geq 0, \quad t\neq 0,  \quad \mbox{for all}\,\,0\leq p\leq m-n/2 \,\,\mbox{even}.\end{multline}
\noindent As we shall demonstrate below, ${\mathcal L}^{m,n}_o(-\partial_t,-p(p+n-2))$, viewed as a polynomial in $\partial_t$, has a double root at zero for all $0\leq p\leq m-n/2$ even and hence, the result of its action (as an operator) on $\widetilde g$ is the same as the result of its action on $h$. Thus, it is indeed enough to consider $h$ in place of $\widetilde g$ in the terms corresponding to such values of $p$.

Next, due to \eqref{eq5.6}--\eqref{eq5.7}, the roots of ${\mathcal L}^{m,n}_o(-\partial_t,-p(p+n-2))$ as a polynomial in $\partial_t$ are given by 
$$ \pm \Bigl(p+n/2+m-2j-2\Bigr), \quad j=0, ..., m-1.$$

\noindent Clearly, the sequence of the roots is symmetric with respect to zero and moreover, the roots can be ordered in non-decreasing order as 
single roots given by 
\begin{equation} \mbox{even numbers\quad  from }\quad -p-m-n/2+2\quad \mbox{ to }\quad p-m+n/2-2, \label{eq6.29}
\end{equation}

\noindent then double roots given by 
\begin{equation}\label{eq6.30}
\mbox{even numbers\quad  from }\quad p-m+n/2\quad \mbox{ to }\quad -p+m-n/2, 
\end{equation}

\noindent (including 0), and then single roots given by 
\begin{equation} \mbox{even numbers\quad  from }\quad -p+m-n/2+2\quad \mbox{ to }\quad p+m+n/2-2. \label{eq6.31}
\end{equation}

Now, in the spirit of the Proof of Proposition~\ref{p3.1}, let us denotes the non-decreasing sequence of the roots of ${\mathcal L}^{m,n}_o(-\partial_t,-p(p+n-2))$ for $p=0$ by
\begin{equation}\label{eq6.32}
b_{m-1}\leq b_{m-2}\leq ...\leq b_1<b_0=0=a_0<a_1\leq ...\leq a_m,
\end{equation}

\noindent and, respectively, the roots of ${\mathcal L}^{m,n}_o(-\partial_t,-p(p+n-2))$ for a given even $0<p\leq m-n/2$ by
\begin{equation}\label{eq6.33}
d_{m-1}\leq d_{m-2}\leq ...\leq d_1<d_0=0=c_0<c_1\leq...\leq c_m.
\end{equation}

It is not hard to see that for all even $0<p\leq m-n/2$ there holds  $[b_i,a_i]\subset [d_i,c_i]$ for all $i=0,...,m-1$, and moreover, $a_0=c_0=0=b_0=d_0$ are the only zero roots. We are now in a position to apply Corollary~\ref{c3.3}. The only problem lies in extra zero roots which appear in the current scenario. To avoid it, we do the following.

The goal is to establish \eqref{eq6.28}.  In the framework of Corollary~\ref{c3.3}, let us denote by $P$ the differential operator of the order $2(m-1)$ such that $\partial_t^2 P(\partial_t)=(-1)^m{\mathcal L}^{m,n}_o(-\partial_t,0)$, $f:=\delta$, and let $Q$ be the differential operator of the order $2(m-1)$ such that $\partial_t^2  Q(\partial_t)=(-1)^m{\mathcal L}^{m,n}_o(-\partial_t,-p(p+n-2))$ with $0<p\leq m-n/2$ even. Then the roots of $P$ and $Q$ are, respectively, $b_i$, $a_i$, $i=1,...,m-1$, and  $d_i$, $c_i$, $i=1,...,m-1$, as above. They satisfy $[b_i,a_i]\subset [d_i,c_i]$, $i=1,...,m-1$, and they are all non-zero. Hence, applying Corollary~\ref{c3.3} with $m-1$ in place of $m$ we have
\begin{equation}\label{eq6.34} P(\partial_t)\,H=\delta \quad \Longrightarrow\quad Q(\partial_t)\,H \,\,\mbox{is a positive distribution for all}\,\,0\leq p\leq m-n/2 \,\,\mbox{even},\end{equation}

\noindent assuming that $H$ is the (unique) solution to $P(\partial_t)\,H=\delta$ which is bounded and vanishes at $+\infty$. It remains to show that $H=(-1)^m \partial_t^2h$, where $h$ is the unique solution to ${\mathcal L}^{m,n}_o(-\partial_t,0)h=\delta$ which vanishes at $+\infty$ and has at most linear growth or decay at $-\infty$ (that is, $h$ is given by Lemma~\ref{l6.1}). However, evidently, $\partial_t^2 [(-1)^m h]$ is bounded and vanishes at $+\infty$. Since  a solution  with such decay properties  is unique, \eqref{eq6.34} indeed proves \eqref{eq6.28}.

\vskip 0.08 in \noindent {\bf Step IV(b), $m-n/2$ even, $0\leq p\leq m-n/2$ odd}. We first note that  
\begin{equation}\label{eq6.35}
{\mathcal L}^{m,n}_o(0,-p(p+n-2))\geq C, \quad\mbox{for $m-n/2$ even and $0\leq p\leq m-n/2$ odd,}
\end{equation}
 \noindent where $C>0$ is a {\it strictly} positive constant depending on $m$ and $n$ only. This can be seen immediately from \eqref{eq5.13}--\eqref{eq5.14}, or, alternatively, from \eqref{eq5.6}--\eqref{eq5.7}. Indeed, the polynomial ${\mathcal L}^{m,n}_o(0,-p(p+n-2))$ is a product of squares. Thus, it is non-negative. Moreover, if $m-n/2$ even, ${\mathcal L}^{m,n}_o(0,-p(p+n-2))$ can only be zero when $p$ is even too. Thus, one can take $C$ as a (strictly positive) minimum of ${\mathcal L}^{m,n}_o(0,-p(p+n-2))$ over $0\leq p\leq m-n/2$ odd, justifying \eqref{eq6.35}. 
 
 Thus, for any  $0\leq p\leq m-n/2$ odd
\begin{multline}\label{eq6.36}
\left|\Bigg({\mathcal L}^{m,n}_o(-\partial_t,-p(p+n-2))
-{\mathcal L}^{m,n}_o(-\partial_t,0) \Bigg) \,\,\widetilde g(t,\tau)\right|\leq C_1+C_2(C_R+t)\\[4pt] \leq {\mathcal L}^{m,n}_o(0,-p(p+n-2))\left(C_3+C_4(C_R+t)\right),\quad \quad t,\tau\geq \log (2R)^{-1},
\end{multline}

\noindent where we used \eqref{eq6.19}--\eqref{eq6.20} for the first inequality above and \eqref{eq6.35} for the second one. The constants $C_i$, $i=1,2,3,4$, depend on $m$ and $n$ only. Note that the operator ${\mathcal L}^{m,n}_o(-\partial_t,-p(p+n-2))
-{\mathcal L}^{m,n}_o(-\partial_t,0)$ has order strictly smaller than $2m$, thus, the corresponding derivatives of $\widetilde g$ are continuous at $t=\tau.$

Therefore, for any $0\leq p\leq m-n/2$ odd we have 
\begin{multline}\label{eq6.37}\left| \int_{\RR} 
v_{pl}^2(t) \Bigg({\mathcal L}^{m,n}_o(-\partial_t,-p(p+n-2))
-{\mathcal L}^{m,n}_o(-\partial_t,0) \Bigg) \,\,\widetilde g(t,\tau)\,  dt\right|\\[4pt]
\leq C  \int_{\RR} 
v_{pl}^2(t) {\mathcal L}^{m,n}_o(0,-p(p+n-2))\left(C_3+C_4(C_R+t)\right)\,  dt\\[4pt]
\leq C\int_{\RR^n}(-\Delta)^m u(x)\,u(x)\,\left(C_3+C_4\,(C_R+\log |x|^{-1})\right) \,dx,
 \end{multline}

\noindent as desired.

 \vskip 0.08 in \noindent {\bf Step V, $0\leq p\leq m-n/2$, $m-n/2$ odd}. The argument in this case is very similar to that in Step IV, except that the approach to even $p$ for $m-n/2$ odd resembles the approach to odd $p$ for $m-n/2$ even and vice versa. 
 
Indeed, when $m-n/2$ is odd, the case of $0\leq p\leq m-n/2$ even can be handled following verbatim the argument in Step IV(b), observing that 
\begin{equation}\label{eq6.38}
{\mathcal L}^{m,n}_o(0,-p(p+n-2))\geq C, \quad\mbox{for $m-n/2$ odd and $0\leq p\leq m-n/2$ even,}
\end{equation}
 \noindent where $C>0$ is a {\it strictly} positive constant depending on $m$ and $n$ only. Then, as in Step IV(b), when $m-n/2$ is odd,
for any $0\leq p\leq m-n/2$ even we have with $p_0=1$
\begin{multline}\label{eq6.39}\Bigg| \int_{\RR} 
v_{pl}^2(t) \Bigg({\mathcal L}^{m,n}_o(-\partial_t,-p(p+n-2))
-{\mathcal
L}^{m,n}_o(-\partial_t,-p_0(p_0+n-2))\Bigg) \,\,\widetilde g(t,\tau)\,  dt\Bigg|\\[4pt]
\leq C  \int_{\RR} 
v_{pl}^2(t) {\mathcal L}^{m,n}_o(0,-p(p+n-2))\left(C_3+C_4(C_R+t)\right)\,  dt\\[4pt]
\leq C\int_{\RR^n}(-\Delta)^m u(x)\,u(x)\,\left(C_3+C_4\,(C_R+\log |x|^{-1})\right) \,dx.
 \end{multline}

Thus, it remains to treat $0\leq p\leq m-n/2$ odd. The argument essentially follows Step IV(a). First of all, the roots of ${\mathcal L}^{m,n}_o(-\partial_t,-p(p+n-2))$ for any $0\leq p\leq m-n/2$ odd, $m-n/2$ odd, are exactly described by \eqref{eq6.29}--\eqref{eq6.31}, including the parity. Hence, for all such $p$ the operator ${\mathcal L}^{m,n}_o(-\partial_t,-p(p+n-2))$ kills constants and it is enough to show that
with $m-n/2$ odd
\begin{multline}\label{eq6.40}{\mathcal L}^{m,n}_o(-\partial_t,-1(1+n-2))h=\delta \quad \Longrightarrow \\[4pt] {\mathcal L}^{m,n}_o(-\partial_t,-p(p+n-2))h(t)\geq 0, \quad t\neq 0,  \quad \mbox{for all }\,\,0\leq p\leq m-n/2 \,\,\mbox{odd}.\end{multline}

As before, the crux of the matter is the structure of the roots of ${\mathcal L}^{m,n}_o(-\partial_t,-p(p+n-2))$ and ${\mathcal L}^{m,n}_o(-\partial_t,-1(1+n-2))$ in the context of the results of Section~\ref{s3}. However, if we denote the set of the numbers described by \eqref{eq6.29}--\eqref{eq6.31} as $S^e_{m,n,p}$ when $m-n/2$ is even and $0\leq p\leq m-n/2$ is even and as $S^o_{m,n,p}$ when $m-n/2$ is odd and $0\leq p\leq m-n/2$ is odd, then they can be identified via $S^o_{m,n,p}=S^e_{m,n+2,p-1}$. Hence, the results of Step IV(a)  translate into the present setting and give exactly \eqref{eq6.40}, as desired.

\vskip 0.08 in \noindent {\bf Step VI, $p> m-n/2$}. Our treatment of $p> m-n/2$ does not depend on whether $m-n/2$ is even or odd, and we record it here in full generality. 

First of all, by \eqref{eq5.8} we have ${\mathcal B}_j(q)^2\geq C\,q(q+n-2)$ for all $q>m-n/2$ and $j=0,...,m-1,$ with $C>0$ depending on $m$ and $n$ only. Thus, 
\begin{equation}\label{eq6.41}
{\mathcal L}^{m,n}_o(0,-p(p+n-2))\geq Cp^m(p+n-2)^m, \quad\mbox{for $p>m-n/2$,}
\end{equation}

\noindent with $C>0$ depending on $m$ and $n$ only. Secondly, recalling the definition of $p_0$ in \eqref{eq6.23}, we have
\begin{multline}\label{eq6.42} \left|\Bigg({\mathcal L}^{m,n}_o(-\partial_t,-p(p+n-2))
-{\mathcal L}^{m,n}_o(-\partial_t,-p_0(p_0+n-2)) \Bigg) \,\,\widetilde g(t,\tau)\right|\\[4pt]
\leq \left(C_1+C_2(C_R+t)\right) p^m(p+n-2)^m,  \end{multline}
\noindent by \eqref{eq6.19}--\eqref{eq6.20}. Combining \eqref{eq6.41} with \eqref{eq6.42} yields then  
for any $p> m-n/2$ 
\begin{multline}\label{eq6.43}\sum_{p=m-n/2+1}^\infty \sum_{l=-p}^{p}\left| \int_{\RR} 
v_{pl}^2(t) \Bigg({\mathcal L}^{m,n}_o(-\partial_t,-p(p+n-2))
-{\mathcal L}^{m,n}_o(-\partial_t,-p_0(p_0+n-2)) \Bigg) \,\,\widetilde g(t,\tau)\,  dt\right|\\[4pt]
\leq C \sum_{p=m-n/2+1}^\infty \sum_{l=-p}^{p} \int_{\RR} 
v_{pl}^2(t) p^m(p+n-2)^m\left(C_1+C_2(C_R+t)\right)\,  dt\\[4pt]
\leq\sum_{p=m-n/2+1}^\infty \sum_{l=-p}^{p} \int_{\RR} 
v_{pl}^2(t) {\mathcal L}^{m,n}_o(0,-p(p+n-2))\left(C_3+C_4(C_R+t)\right)\,  dt\\[4pt]
\leq C\int_{\RR^n}(-\Delta)^m u(x)\,u(x)\,\left(C_3+C_4\,(C_R+\log |x|^{-1})\right) \,dx.
 \end{multline}
 
This finishes the argument for $I_1^o$.

\vskip 0.08 in \noindent {\bf Step VII, the bound on $I_2^o$}. The estimate on $I_2^o$ is essentially already incorporated in \eqref{eq5.1}, \eqref{eq5.2}. Indeed, 
\begin{eqnarray}\label{eq6.44}\nonumber
I_2^o&\leq &  \sum_{\stackrel{k\geq 1,\,i\geq 0}{2i+k\leq
2m}}\sum_{\stackrel{l\geq
1,\,j\geq 0}{2l+j\leq
k}}|c_{lj}^{ki}|\int_{\RR}\int_{S^{n-1}}
(\partial_t^l\nabla_\omega^iv)^2\,|\partial_t^j\,\widetilde g (t,\tau)|\,dt\\[4pt]
&&\qquad\leq C \sum_{\stackrel{k\geq 1,\,i\geq 0}{ i+k\leq
m}}\int_{\RR}\int_{S^{n-1}} (\partial_t^k\nabla_\omega^i
v(t,\omega))^2 (C_1+C_2(C_R+t))\,d\omega dt\nonumber\\[4pt]
&&\qquad  \leq \int_{\RR^n}(-\Delta)^m u(x)\,u(x)\,\left(C_3+C_4\,(C_R+\log |x|^{-1})\right) \,dx,
\end{eqnarray}

\noindent employing \eqref{eq6.19}--\eqref{eq6.20} for the first inequality and  \eqref{eq5.1}, \eqref{eq5.2} for the second one. 

This finishes the argument. \ep

\section{Pointwise and local $L^2$ estimates for solutions to the polyharmonic equation}\label{s7}
\setcounter{equation}{0}

This section is devoted to the proof of Theorem~\ref{t1.1}. In addition, we
will establish sharp local estimates for solutions in a neighborhood of a boundary point and ``at infinity", that is, when moving away from a given $Q\in\po$. 

To start, we record for future reference a 
well-known result that follows from the energy estimates for solutions
of elliptic equations. 

\begin{lemma}\label{l7.1} Let $\Omega$ be an arbitrary
domain in $\RR^n$, $n\geq 2$, $Q\in\RR^n\setminus\Omega$ and $R>0$. Suppose
\begin{equation}\label{eq7.1}
(-\Delta)^m u=f \,\,{\mbox{in}}\,\,\Omega, \quad f\in
C_0^{\infty}(\Omega\setminus B_{4R}(Q)),\quad u\in \ring
W^{m,2}(\Omega).
\end{equation}

\noindent Then
\begin{equation}\label{eq7.2}
\sum_{i=1}^m\frac{1}{\rho^{2m-2i}}\int_{B_{\rho}(Q)\cap\Omega}|\nabla^i u|^2\,dx
\leq
\frac{C}{\rho^{2m}}\int_{C_{\rho,2\rho}(Q)\cap\Omega}|u|^2\,dx
\end{equation}

\noindent for every $\rho<2R$.
\end{lemma}

The following Proposition   reflects the rate of growth of solutions near a boundary point encoded in Theorems~\ref{t4.2}, \ref{t6.3}, and ultimately provides a passage to Theorem~\ref{t1.1}.

\begin{proposition}\label{p7.2}
Let $\Omega$ be a bounded domain in $\RR^n$, $2\leq n \leq 2m+1$,
$Q\in\RR^n\setminus\Omega$, and $R>0$. Suppose
\begin{equation}\label{eq7.3}
(-\Delta)^m u=f \,\,{\mbox{in}}\,\,\Omega, \quad f\in
C_0^{\infty}(\Omega\setminus B_{4R}(Q)),\quad u\in \ring
W^{m,2}(\Omega).
\end{equation}

\noindent Then
\begin{equation}\label{eq7.4}
\frac{1}{\rho^{2\lambda +n-1}}\int_{S_{\rho}(Q)\cap\Omega}|u(x)|^2\,d\sigma_x
\leq \frac{C}{R^{2\lambda +n}} \int_{C_{R,4R}(Q)\cap\Omega} |u(x)|^2\,dx\quad
{\mbox{ for every}}\quad \rho<R,
\end{equation}

\noindent where $C$ is a  constant depending on $m$ and $n$ only, and  $\lambda=\left[m-\frac n2+\frac 12\right]$, that is, 
\begin{equation}\label{eq7.5}
\lambda=\left\{\begin{array}{l} m-n/2+1/2 \quad \,\,\mbox{when $n$ is odd},\\[8pt] m-n/2\qquad\qquad\mbox{when $n$ is even}.\end{array}
\right.
\end{equation}

Moreover, for every $x\in B_{R/4}(Q)\cap\Omega$
\begin{equation}\label{eq7.15}
|\nabla^{i} u(x)|^2\leq C\, \frac{|x-Q|^{2\lambda-2i}}{R^{n+2\lambda}}\int_{C_{R/4,4R}(Q)\cap\Omega}
|u(y)|^2\,dy,\qquad 0\leq i\leq \lambda.
\end{equation}

In particular, for every bounded domain $\Omega\subset\RR^n$ the
solution to the boundary value problem {\rm (\ref{eq7.3})}
satisfies
\begin{equation}\label{eq7.16}
|\nabla^{m-n/2+1/2} u|\in L^\infty(\Omega)\,\,\mbox{when $n$ is odd\quad and \quad} |\nabla^{m-n/2} u|\in L^\infty(\Omega)\,\,\mbox{when $n$ is even}.
\end{equation}

\end{proposition}

This Proposition will conclude the proof of Theorem~\ref{t1.1}.

\vskip 0.08 in

\noindent {\it Proof of Proposition~\ref{p7.2}}. Without loss of generality  we can assume that $Q=O$. Let us
approximate $\Omega$ by a sequence of domains with smooth boundaries
$\left\{\Omega_n\right\}_{n=1}^\infty$ satisfying
\begin{equation}\label{eq7.6}
\bigcup_{n=1}^\infty \Omega_n=\Omega\quad\mbox{and}\quad
{\overline{\Omega}}_n\subset\Omega_{n+1} \quad \mbox{for every}
\quad n\in\NN.
\end{equation}

\noindent Choose $n_0\in\NN$ such that ${\rm supp}\,f\subset
\Omega_n$ for every $n\geq n_0$ and denote by $u_n$ a unique
solution of the Dirichlet problem
\begin{equation}\label{eq7.7}
(-\Delta)^m u_n = f\quad {\rm in} \quad \Omega_n, \quad  u_n\in \ring
W^{m,2}(\Omega_n), \quad n\geq n_0.
\end{equation}

\noindent The sequence  $\{u_n\}_{n=n_0}^\infty$ converges to $u$ in
$\ring W^{m,2}(\Omega)$ (see, e.g., \cite{Necas}, \S{6.6}).

Next,  take some $\eta\in C_0^\infty(B_{2R})$ such that
\begin{equation}\label{eq7.8}
0\leq\eta\leq 1\,\, {\rm in}\,\,B_{2R}, \quad \eta=1\,\, {\rm
in}\,\,B_R\quad {\rm and} \quad |\nabla^k \eta|\leq CR^{-k}, \quad
k\leq 2m.
\end{equation}

\noindent  Also, fix $\tau=\log\rho^{-1}$ and let $g$ be the
function defined in (\ref{eq4.12}) when $n$ is odd and by \eqref{eq6.15} when $n$ is even. 

We observe that when $n$ is odd, the formula \eqref{eq4.12} yields 
\begin{equation}\label{eq7.8.1}
\left|\nabla_x^k\, g(\log|x|^{-1},\log \rho^{-1})\right| \leq C |x|^{-k-1}, \quad 0\leq k\leq 2m, \quad x\in\RR^n\setminus\{0\}, \quad\rho\in (0,\infty),
\end{equation}

\noindent while for $n$ even by \eqref{eq6.15} and \eqref{eq6.19}--\eqref{eq6.20}
\begin{equation}\label{eq7.8.2}
\left|g(\log|x|^{-1},\log \rho^{-1})\right| \leq C_1+C_2 (C_R+\log |x|^{-1}), \quad 0<|x|, \rho<2R,
\end{equation}

\noindent and 
\begin{equation}\label{eq7.8.3}
\left|\nabla_x^k\, g(\log|x|^{-1},\log \rho^{-1})\right| \leq C |x|^{-k}, \quad 1\leq k\leq 2m, \quad 0<|x|, \rho<2R.
\end{equation}

\noindent Here, as usually, we assume $|x|\neq \rho$ when $k=2m$, and lower derivatives of $g$ as well as $g$ itself are defined at $x$ such that $|x|\neq \rho$ by continuity. Hence, in particular, 
when $n$ is odd, 
\begin{equation}\label{eq7.8.4}
\left|\nabla_x^k\, g(\log|x|^{-1},\log \rho^{-1})\right| \leq C R^{-k-1}, \quad 0\leq k\leq 2m, \quad x\in C_{R,2R}, \quad\rho<R,
\end{equation}

\noindent and when $n$ is even, 
\begin{equation}\label{eq7.8.5}
\left|\nabla_x^k\, g(\log|x|^{-1},\log \rho^{-1})\right| \leq C R^{-k}, \quad 0\leq k\leq 2m, \quad x\in C_{R,2R}, \quad\rho<R, \end{equation}

\noindent since for $x\in C_{R,2R}$ we have $$C_R+\log |x|^{-1}=\log(4R)+\log |x|^{-1}=\log \frac{4R}{|x|}\approx C.$$

Consider now 
\begin{equation}\label{eq7.9}
\int_{\RR^n}\bigg([(-\Delta)^m,\eta] u_n(x)\bigg) \bigg(\eta(x)u_n(x)g(\log|x|^{-1},\log \rho^{-1})\bigg)\,dx,
\end{equation}

\noindent where the brackets denote the commutator $$[(-\Delta)^m,\eta] u_n(x)=(-\Delta)^m\left(\eta(x)u_n(x)\right)-\eta(x)(-\Delta)^m u_n(x),$$ the integral in (\ref{eq7.9}) is understood in the sense of pairing between $\ring W^{m,2}(\Omega_n)$ and its dual. 
Evidently, the
support of the integrand is a subset of ${\rm
supp}\,\nabla\eta\subset C_{R,2R}$, and therefore, the expression
in (\ref{eq7.9}) is equal to 
\begin{multline}\label{eq7.10}
\int_{C_{R,2R}}\bigg([(-\Delta)^m,\eta] u_n(x)\bigg) \bigg(\eta(x)u_n(x)g(\log|x|^{-1},\log \rho^{-1})\bigg)\,dx\\[4pt]\leq 
C\sum_{i=0}^m \frac{1}{R^{2\lambda+n-2i}}\int_{C_{R,2R}} |\nabla^i
u_n(x)|^2\,dx\leq \frac{C}{R^{2\lambda+n}}\int_{C_{R,4R}} |u_n(x)|^2\,dx,
\end{multline}

\noindent using Cauchy-Schwarz inequality, \eqref{eq7.8.4}--\eqref{eq7.8.5}, and \eqref{eq7.8} for the first inequality and 
Lemma~\ref{l7.1} for the second one.

On the other hand, since $u_n$ is biharmonic in $\Omega_n\cap B_{4R}$ and $\eta$ is
supported in $B_{2R}$, one can see that $\eta\,(-\Delta)^m u_n =0$ and hence the integral in \eqref{eq7.9} is equal to 
\begin{equation}\label{eq7.11}
\int_{\RR^n}(-\Delta)^m\left(\eta(x)u_n(x)\right)\, \eta(x)u_n(x)g(\log|x|^{-1},\log \rho^{-1})\,dx.
\end{equation}
\noindent To estimate it we employ Theorems~\ref{t4.2} and \ref{t6.3}
with $u=\eta\,u_n$. The results \eqref{eq4.11} and \eqref{eq6.14} hold for such a choice
of $u$. This can be seen directly by inspection of the argument or
one can approximate each $u_n$ by a sequence of
$C_0^\infty(\Omega_n)$ functions in $\ring W^{m,2}(\Omega_n)$ and then
take a limit using that $O\not\in\overline\Omega_n$. Then
(\ref{eq7.11}) is bounded from below by
\begin{equation}\label{eq7.12}
\frac{C}{\rho^{2\lambda+n-1}}\int_{S_{\rho}}|\eta(x)u_n(x)|^2\,d\sigma_x.
\end{equation}

\noindent Hence, for every $\rho<R$
\begin{equation}\label{eq7.13}
\frac{1}{\rho^{2\lambda+n-1}}\int_{S_{\rho}}|u_n(x)|^2\,d\sigma_x \leq \frac{C}{R^{2\lambda+n}}\int_{C_{R,4R}} |u_n(x)|^2\,dx.
\end{equation}

\noindent Now the proof of \eqref{eq7.4} can be finished by taking the limit as $n\to \infty$.  

Going further, by virtue of the local estimates for solutions of elliptic equations, the bound in (\ref{eq7.4}) transforms into the uniform  pointwise estimates for $\nabla^\lambda u$ (observe that the exponent $\lambda$, as defined in (\ref{eq7.5}), is an integer number both when $n$ is odd and when $n$ is even).

Indeed, by interior estimates for solutions of the elliptic equations
(see \cite{ADN})
\begin{equation}\label{eq7.17}
|\nabla^i u(x)|^2\leq \frac{C}{d(x)^n} \int_{B_{d(x)/2}(x)} |\nabla^i
u(y)|^2\,dy,\qquad 0\leq i\leq m,
\end{equation}

\noindent where $d(x)$ denotes the distance from $x$ to $\po$. Let
$x_0$ be a point on the boundary of $\Omega$ such that
$d(x)=|x-x_0|$. Since $x\in B_{R/4}(Q)\cap\Omega$ and
$Q\in\RR^n\setminus\Omega$, we have $x\in B_{R/4}(x_0)$, and
therefore
\begin{equation}\label{eq7.18}
\frac{1}{d(x)^n}\int_{B_{d(x)/2}(x)} |\nabla^i u(y)|^2\,dy\leq
\frac{C}{d(x)^{n+2i}} \int_{B_{ 2d(x)}(x_0)} |u(y)|^2\,dy \leq C\,
\frac{d(x)^{-2i+2\lambda}}{R^{n+2\lambda}}\int_{C_{3R/4,3R}(x_0)} |u(y)|^2\,dy,
\end{equation}

\noindent as soon as $i\leq \lambda$. Indeed, the first inequality in \eqref{eq7.18} follows from 
Lemma~\ref{l7.1}. Turning to the second one, we observe that   $d(x)\leq R/4$ and
therefore, $2d(x)< 3R/4$. On the other hand, $u$ is polyharmonic  in
$B_{4R}(Q)\cap\Omega$ and $|Q-x_0|\leq |Q-x|+|x-x_0|\leq R/2.$
Hence, $u$ is polyharmonic in $B_{3R}(x_0)\cap\Omega$ and
Proposition~\ref{p7.2} holds with $x_0$ in place of $Q$, $3R/4$ in
place of $R$ and $\rho=2d(x)$. Furthermore, $C_{3R/4,3R}(x_0)\subset C_{R/4,4R}(Q),$
and that finishes the argument for the second inequality in (\ref{eq7.18}). Clearly, $d(x)\leq |x-Q|$, so that \eqref{eq7.17}--\eqref{eq7.18} entail \eqref{eq7.15}.\ep

At this point, we are ready to address  the behavior of solutions ``at infinity". 
\begin{proposition}\label{p7.4} Let $\Omega$ be a bounded domain in $\RR^n$, $2\leq n \leq 2m+1$,
$Q\in\RR^n\setminus\Omega$, $r>0$ and assume that
\begin{equation}\label{eq7.21}
(-\Delta)^m u=f \,\,{\mbox{in}}\,\,\Omega, \quad f\in
C_0^{\infty}(B_{r/4}(Q)\cap\Omega),\quad u\in \ring W^{m,2}(\Omega).
\end{equation}

\noindent Then
\begin{equation}\label{eq7.22}
\rho^{2\lambda+n+1-4m}\int_{S_{\rho}(Q)\cap\Omega}|u(x)|^2\,d\sigma_x
\leq C\, r^{2\lambda+n-4m}\int_{ C_{r/4,r}(Q)\cap\Omega} |u(x)|^2\,dx,
\end{equation}

\noindent for any $\rho>r$ and $\lambda$ given by \eqref{eq7.5}.

Furthermore, for any $x\in\Omega\setminus B_{4r}(Q)$
\begin{equation}\label{eq7.23}
|\nabla^i u(x)|^2\leq
C\,\frac{r^{2\lambda+n-4m}}{|x-Q|^{2\lambda+2n-4m+2i}}\int_{C_{r/4,4r}(Q)\cap\Omega} |u(y)|^2\,dy,\qquad 0\leq i\leq \lambda.
\end{equation}
\end{proposition}

\bp Without loss of generality, one can consider $Q=O$. Retain the
approximation of $\Omega$ with the sequence of smooth domains
$\Omega_n$ satisfying (\ref{eq7.6}) and define $u_n$ according to
(\ref{eq7.7}). We denote by ${\mathcal I}$ the inversion $x\mapsto
y=x/|x|^2$ and by $U_n$ the Kelvin transform of $u_n$,
\begin{equation}\label{eq7.24}
U_n(y):=|y|^{2m-n}\, u_n(y/|y|^2), \quad y\in {\mathcal I}(\Omega_n).
\end{equation}

\noindent Then
\begin{equation}\label{eq7.25}
(-\Delta)^m U_n(y)=|y|^{-n-2m} ((-\Delta)^m u_n)(y/|y|^2),
\end{equation}

\noindent and therefore, $U_n$ is polyharmonic in ${\mathcal
I}(\Omega_n)\cap B_{4/r}$. Moreover, 
$U_n\in \ring W^{m,2}({\mathcal I}(\Omega_n))$ if and only if $u_n\in \ring W^{m,2}(\Omega_n).$
Note also that $\Omega_n$ is a bounded domain with
$O\not\in \overline{\Omega}_n$, hence, so is ${\mathcal I}(\Omega_n)$
and $O\not\in \overline{{\mathcal I}(\Omega_n)}$.

At this point we can invoke  Proposition~\ref{p7.2} to show that
\begin{equation}\label{eq7.27}
\rho^{2\lambda+n-1}\int_{S_{1/\rho}}|U_n(y)|^2\,d\sigma_y \leq C\,
r^{2\lambda+n}\int_{C_{1/r,4/r}} |U_n(y)|^2\,dy,
\end{equation}

\noindent so that
\begin{equation}\label{eq7.28}
\rho^{2\lambda+n+1-4m}\int_{S_{\rho}}|u_n(x)|^2\,d\sigma_x \leq
C\,r^{2\lambda+n-4m} \int_{ C_{r/4,r}} |u_n(x)|^2\,dx.
\end{equation}
\noindent using the substitution (\ref{eq7.24}) and the
change of coordinates.

This furnishes the desired $L^2$ estimates. It remains to establish pointwise control. To this end, let us fix some $x\in\Omega\setminus B_{4r}(Q)$ and recall that 
$$u_n(x)=|x|^{2m-n}\, U_n(x/|x|^2),$$
\noindent hence,
\begin{equation}\label{eq7.29}
|\nabla^i u_n(x)|\leq C \sum_{k=0}^i |x|^{2m-n-i-k}\,(\nabla^k
U_n)(x/|x|^2).
\end{equation}

\noindent Therefore, combining
(\ref{eq7.29}) and Proposition~\ref{p7.2} applied to the function
$U_n$, we deduce that
\begin{equation}\label{eq7.30}
|\nabla^i u_n(x)|^2\leq C
\,\frac{r^{n+2\lambda}}{|x|^{2\lambda+2n-4m+2i}}\int_{C_{1/(4r),4/r}}|U_n(z)|^2\,dz=
C  \,\frac{r^{n+2\lambda-4m}}{|x|^{2\lambda+2n-4m+2i}}\int_{C_{r/4,4r}}|u_n(z)|^2\,dz,
\end{equation}

\noindent for all $0\leq i\leq\lambda$.
Now one  can use the limiting procedure to complete the
argument. Indeed, since $u_n$ converges to $u$ in $\ring
W^{m,2}(\Omega)$ and $i\leq\lambda\leq m$, the integrals on the right-hand side of (\ref{eq7.28}) and  (\ref{eq7.30})
converge to the corresponding integrals with
$u_n$ replaced by $u$. Turning to $|\nabla^i u_n(x)|$, we observe
that both $u_n$ and $u$ are $m$-harmonic in a neighborhood of $x$,
in particular, for sufficiently small $d$
\begin{equation}\label{eq7.31}
|\nabla^i (u_n(x)-u(x))|^2\leq
\,\frac{C}{d^{n+2i}}\int_{B_{d/2}(x)}|u_n(z)-u(z)|^2\,dz.
\end{equation}

\noindent As $n\to\infty$, the integral on the right-hand side of
(\ref{eq7.31}) vanishes and therefore, $|\nabla^i u_n(x)|\to |\nabla^i
u(x)|$. Similar considerations apply to $u_n(x)$ and the integral on the left-hand side of \eqref{eq7.28}. \ep

\section{Green's function estimates}\label{s8}
\setcounter{equation}{0}

The present section addresses sharp pointwise estimates on Green's function of the polyharmonic equation and the regular part of Green's function, that is, the difference between Green's function and the fundamental solution. We have discussed some of these bounds, as well as their applications for the solutions of the Dirichlet problem, in \cite{MayMazGr}. The estimates listed below are more refined compared to the statements in \cite{MayMazGr} and we present them here with full proofs. 

To start, let us recall the definition of the fundamental solution for the polyharmonic equation (see. e.g., \cite{PolyharmBook}). A fundamental solution for the $m$-Laplacian is a linear combination of the characteristic singular solution (defined below) and any $m$-harmonic function in $\RR^n$. The characteristic singular solution is
\begin{eqnarray}\label{eq8.2}
&C_{m,n}|x|^{2m-n},\quad &\mbox{if $n$ is odd, or if $n$ is even with $n\geq 2m+2$},\\[4pt]
\label{eq8.3}
&C_{m,n}|x|^{2m-n} \log |x|,\quad &\mbox{if $n$ is even with $n\leq 2m$}.
\end{eqnarray}

\noindent The exact expressions for constants  $C_{m,n}$ can be found in \cite{PolyharmBook}, p.8.  Hereafter we will use the fundamental solution given by 
\begin{equation}\label{eq8.4}
\Gamma(x)=C_{m,n}  \left\{\begin{array}{l}
|x|^{2m-n},\quad \mbox{if $n$ is odd},\\[4pt]
|x|^{2m-n}\log \frac{{\rm diam}\,\Omega}{|x|},\quad \mbox{if $n$ is even and $n\leq 2m$},\\[4pt]
|x|^{2m-n},\quad \mbox{if $n$ is even and $n\geq 2m+2$}.
\end{array}
\right.
\end{equation}

As in the Introduction, by $G$ we denote the Green function of the operator $(-\Delta)^m$ (see \eqref{eq8.1}) and $S$ stands for its regular part, that is $S(x,y)=G(x,y)-\Gamma(x-y)$, $x,y\in\Omega$.

\begin{theorem}\label{t8.1} Let $\Omega\subset \RR^n$ be an arbitrary bounded
domain, $m\in\NN$, $n\in [2,2m+1]\cap\NN$, and let $\lambda$ retain the significance of \eqref{eq7.5}.  Fix any number $N\geq 25$.
Then there exists a  constant $C$ depending on $m, n, N$ only such that for
every $x,y\in\Omega$ the following estimates hold.

If $n\in [3,2m+1] \cap\NN$ is odd then
\begin{equation}\label{eq8.25}
|\nabla_x^i \nabla_y^{j}G(x,y)|\leq
C\,\frac{d(y)^{\lambda-j}}{|x-y|^{\lambda+n-2m+i}}, \quad \mbox{when}\quad |x-y|\geq N\,d(y),\quad 0\leq i,j\leq \lambda,
\end{equation}

\noindent and 
\begin{equation}\label{eq8.26}
|\nabla_x^i \nabla_y^{j}G(x,y)|\leq
C\,\frac{d(x)^{\lambda-i}}{|x-y|^{\lambda+n-2m+j}}, \quad \mbox{when}\quad |x-y|\geq N\,d(x),\quad 0\leq i,j\leq \lambda.
\end{equation}

\noindent Next, 

\begin{multline}\label{eq8.42}
|\nabla_x^i\nabla_y^jG(x,y)|\leq
\frac{C}{|x-y|^{n-2m+i+j}}, \quad \\[4pt]
\mbox{when}\,\,|x-y|\leq N^{-1}\max\{d(x),d(y)\},\,\,\mbox{and}\,\,i+j\geq 2m-n, \,\,0\leq i,j\leq m-n/2+1/2,
\end{multline}

\noindent and 
\begin{multline}\label{eq8.43}
|\nabla_x^i\nabla_y^jG(x,y)|\leq 
{C}\,{\min\{d(x),d(y)\}^{2m-n-i-j}}, \quad\\[4pt] \mbox{when}\,\,|x-y|\leq N^{-1}\max\{d(x),d(y)\},\,\,\mbox{and}\,\,i+j\leq 2m-n,\,\,0\leq i,j\leq m-n/2+1/2.
\end{multline}

\noindent Finally, 
\begin{multline}\label{eq8.53}
|\nabla_x^i\nabla_y^jG(x,y)|\leq  \frac{C}{\min \{d(x),d(y),|x-y|\}^{n-2m+i+j}}\approx \frac{C}{\max \{d(x),d(y),|x-y|\}^{n-2m+i+j}},\\[4pt]
\mbox{when}\,\,N^{-1}\,d(x)\leq |x-y|\leq Nd(x)\,\, \mbox{and}\,\,
N^{-1}\,d(y)\leq |x-y|\leq Nd(y),\quad 0\leq i,j\leq \lambda.
\end{multline}

Furthermore, if $n\in [3,2m+1] \cap\NN$ is odd, the estimates on the regular part of the Green function $S$ are as follows:
\begin{equation}\label{eq8.30}
|\nabla_x^i\nabla_y^j S(x-y)|\leq \frac{C}{|x-y|^{n-2m+i+j}}\quad \mbox{when}\quad |x-y|\geq N \min\{d(x),d(y)\},\quad  0\leq i,j\leq \lambda.\end{equation}
\noindent  

\noindent Next, 
\begin{multline}\label{eq8.40}
|\nabla_x^i\nabla_y^jS(x,y)|\leq
\frac{C}{\max\{d(x),d(y)\}^{n-2m+i+j}}, \quad\\[4pt]
\mbox{when}\,\,|x-y|\leq N^{-1}\max\{d(x),d(y)\},\,\,\mbox{and}\,\,i+j\geq 2m-n, \,\,0\leq i,j\leq m-n/2+1/2,
\end{multline}

\noindent and  
\begin{multline}\label{eq8.41}
|\nabla_x^i\nabla_y^jS(x,y)|\leq
{C}\,{\min\{d(x),d(y)\}^{2m-n-i-j}},  \quad\\[4pt] \mbox{when}\,\,|x-y|\leq N^{-1}\max\{d(x),d(y)\},\,\,\mbox{and}\,\,i+j\leq 2m-n,\,\,0\leq i,j\leq m-n/2+1/2.
\end{multline}

Finally, 
\begin{multline}\label{eq8.54}
|\nabla_x^i\nabla_y^jS(x,y)|\leq  \frac{C}{\min \{d(x),d(y),|x-y|\}^{n-2m+i+j}}\approx \frac{C}{\max \{d(x),d(y),|x-y|\}^{n-2m+i+j}},\\[4pt]
\mbox{when}\,\,N^{-1}\,d(x)\leq |x-y|\leq Nd(x)\,\, \mbox{and}\,\,
N^{-1}\,d(y)\leq |x-y|\leq Nd(y),\quad 0\leq i,j\leq m-n/2+1/2.
\end{multline}

If $n\in [2,2m]\cap\NN$ is even, then \eqref{eq8.25}--\eqref{eq8.26} and \eqref{eq8.53} are valid with $\lambda=m-\frac n2$, and 
\begin{multline}\label{eq8.49}
|\nabla_x^i\nabla_y^jG(x,y)|\leq  
{C}\,{\min\{d(x),d(y)\}^{2m-n-i-j}}\,\left( C'+\log \frac{\min\{d(x),d(y)\}}{|x-y|}\right), \, \\[4pt]\mbox{when}\,|x-y|\leq N^{-1}\max\{d(x),d(y)\}\, \,\,\mbox{and}\,\, 0\leq i,j\leq m-n/2.\end{multline}

Furthermore, if $n\in [2,2m]\cap\NN$ is even, the estimates on the regular part of the Green function $S$ are as follows:
\begin{multline}\label{eq8.31}
|\nabla_x^i\nabla_y^j S(x-y)|\leq C\,|x-y|^{-n+2m-i-j}\left(C'+\log \frac{{\rm diam}\,{(\Omega)}}{|x-y|}\right)\\[4pt]\quad \mbox{when}\quad |x-y|\geq N \min\{d(x),d(y)\},\quad 0\leq i,j\leq m-n/2.
\end{multline}

\noindent Next,
\begin{multline}\label{eq8.46}
|\nabla_x^i\nabla_y^jS(x,y)|\leq
{C}\,{\min\{d(x),d(y)\}^{2m-n-i-j}}\left( C'+\log \frac{{\rm diam}\,\Omega}{{\max\{d(x),d(y)\}}}\right),  \\[4pt] \mbox{when\,\, $|x-y|\leq N^{-1}\max\{d(x),d(y)\}$, \quad $0\leq i,j\leq m-n/2$.}
\end{multline}

\noindent Finally, 
\begin{multline}\label{eq8.55}|\nabla_x^i\nabla_y^jS(x,y)|\leq C\,{\min \{d(x),d(y),|x-y|\}^{2m-n-i-j}}\left(C'+\log \frac{{\rm diam}\,\Omega}{\max \{d(x),d(y),|x-y|\}^{n-2m+i+j}}\right)
\\[4pt]
\mbox{when}\,\,N^{-1}\,d(x)\leq |x-y|\leq Nd(x)\,\, \mbox{and}\,\,
N^{-1}\,d(y)\leq |x-y|\leq Nd(y),\quad 0\leq i,j\leq m-n/2.
\end{multline}

\end{theorem}

Before passing to the proof of the Theorem, we would like to point out that  the bounds on the highest order derivatives highlighted in Theorem~\ref{t1.2} is a particular case of Theorem~\ref{t8.1}. 
%Note that both the Green function and the fundamental solution are symmetric functions of $x$ and $y$ and thus, it is sufficient to consider the cases $d(x)\leq d(y) \leq |x-y|$,  $d(x) \leq |x-y|\leq d(y)$, and $|x-y|\leq d(x) \leq d(y)$, $x,y\in\Omega$, as the estimates for all remaining situations can be recovered from symmetry.  Modulo these symmetry considerations, the bounds in the Corollary below are equivalent to  \eqref{eq8.5}, \eqref{eq8.7} and \eqref{eq8.9}, \eqref{eq8.10} for $i=j=\lambda$ with $\lambda,$ as usually, defined by \eqref{eq7.5}.

\vskip 0.08 in
\noindent {\it Proof of Theorem~\ref{t8.1}}. Recall the definition of the fundamental solution $\Gamma$ in \eqref{eq8.4}. For any  $\alpha$, a multi-index of length less than or equal to $\lambda$, the function $\partial^{\alpha}_y\,\Gamma(x-y)$ can be written as 
\begin{equation}\label{eq8.11}
\partial^{\alpha}_y\,\Gamma(x-y)=P^{\alpha}(x-y)\log \frac{{\rm diam}\,\Omega}{|x-y|}+Q^{\alpha}(x-y).
\end{equation}

\noindent Here, when the dimension is odd, $P^{\alpha}\equiv 0$ and $Q^{\alpha}$ is a homogeneous function of order $2m-n-|\alpha|$. If the dimension is even (and less than or equal to $2m$ by the assumptions of the theorem) then $P^{\alpha}$ and $Q^{\alpha}$ are homogeneous polynomials of order $2m-n-|\alpha|$ as long as $|\alpha|\leq 2m-n$. In both cases, $P^{\alpha}$ and $Q^{\alpha}$ do not depend in any way on the domain $\Omega$.

Consider a
function $\eta$ such that
\begin{equation}\label{eq8.12}
\eta \in C_0^\infty(B_{1/2})\quad\mbox{and}\quad \eta=1
\quad\mbox{in}\quad B_{1/4},
\end{equation}

\noindent and define 
\begin{equation}\label{eq8.13}
{\mathcal R}_{\alpha} (x,y):=\partial^{\alpha}_y\,
G(x,y)-\eta\left(\frac{x-y}{d(y)}\right)\left(-P^{\alpha}(x-y)\log \frac{|x-y|}{d(y)}+Q^{\alpha}(x-y)\right),\qquad x,y\in\Omega.
\end{equation}

For every fixed $y\in\Omega$ the function $x\mapsto
{\mathcal R}_{\alpha}(x,y)$ is a solution of the boundary value problem
\begin{equation}\label{eq8.16}
(-\Delta_x)^m {\mathcal R}_\alpha (x,y)= f_\alpha(x,y) \,\,{\mbox{in}}\,\,\Omega,
\quad f_\alpha(\cdot,y)\in C_0^{\infty}(\Omega),\quad {\mathcal
R}_\alpha(\cdot,y)\in \ring W^{m,2}(\Omega),
\end{equation}

\noindent  where
\begin{equation}\label{eq8.14}
f_{\alpha}(x,y):= (-\Delta_x)^m {\mathcal R}_{\alpha}
(x,y)=-\left[(-\Delta_x)^m,\eta\left(\frac{x-y}{d(y)}\right)
\right]\left(P^{\alpha}(x-y)\log \frac{|x-y|}{d(y)}+Q^{\alpha}(x-y)\right).
\end{equation}

\noindent Indeed, it is not hard to see that for every $\alpha$
\begin{equation}\label{eq8.15}
f_{\alpha}(\cdot,y)\in C_0^\infty(C_{d(y)/4,
d(y)/2}(y))\qquad\mbox{and}\qquad |f_\alpha(x,y)|\leq Cd(y)^{-n-|\alpha|}, \quad
x,y\in\Omega,
\end{equation}

\noindent   so that, in particular, $f_\alpha(\cdot,y)\in C_0^{\infty}(\Omega)$ as stated in \eqref{eq8.16}. Furthermore, by \eqref{eq8.14},
\begin{equation}\label{eq8.17}
\left\|\nabla_x^m{\mathcal
R}_\alpha(\cdot,y)\right\|_{L^2(\Omega)}=\left\|{\mathcal
R}_\alpha(\cdot,y)\right\|_{W^{m,2}(\Omega)}\leq
C\|f_\alpha(\cdot,y)\|_{W^{-m,2}(\Omega)},
\end{equation}

\noindent where $W^{-m,2}(\Omega)$ stands for the Banach space dual
of $\ring W^{m,2}(\Omega)$, i.e.
\begin{equation}\label{eq8.18}
\|f_\alpha(\cdot,y)\|_{W^{-m,2}(\Omega)}=\sup_{v\in \ring
W^{m,2}(\Omega):\,\|v\|_{\ring
W^{m,2}(\Omega)}=1}\int_{\Omega}f_\alpha(x,y)v(x)\,dx.
\end{equation}

 Recall that by Hardy's inequality
\begin{equation}\label{eq8.19}
\left\|\frac{v}{|\cdot-\,Q|^{m}}\right\|_{L^2(\Omega)}\leq C
\left\|\nabla^m v\right\|_{L^2(\Omega)}\quad\mbox{for every}\quad
v\in\ring W^{m,2}(\Omega),\quad Q\in\po.
\end{equation}

\noindent Then by \eqref{eq8.18}
\begin{multline}\label{eq8.20}
\|f_\alpha(\cdot,y)\|_{W^{-m,2}(\Omega)}=\sup_{v\in \ring
W^{m,2}(\Omega):\,\|v\|_{\ring
W^{m,2}(\Omega)}=1} \int_{\Omega}f_\alpha(x,y)v(x)\,dx\\[4pt]
\leq C \sup_{v\in \ring
W^{m,2}(\Omega):\,\|v\|_{\ring
W^{m,2}(\Omega)}=1}
\left\|\frac{v}{|\cdot-y_0|^{m}}\right\|_{L^2(\Omega)}\left\|f_\alpha(\cdot,y)\,\,
|\cdot-y_0|^{m}\right\|_{L^2(\Omega)} \\[4pt]
\qquad \leq C \sup_{v\in \ring
W^{m,2}(\Omega):\,\|v\|_{\ring
W^{m,2}(\Omega)}=1} d(y)^{m} \left\|\nabla^m v\right\|_{L^2(\Omega)}
\|f_\alpha(\cdot,y) \|_{L^2(C_{d(y)/4, d(y)/2}(y))}\\[4pt]
\qquad \leq C  d(y)^{m} 
\|f_\alpha(\cdot,y) \|_{L^2(C_{d(y)/4, d(y)/2}(y))},
\end{multline}

\noindent where $y_0\in\po$ is such that $|y-y_0|=d(y)$. Therefore, by (\ref{eq8.15})
\begin{equation}\label{eq8.21}
\left\|\nabla_x^m{\mathcal R}_\alpha(\cdot,y)\right\|_{L^2(\Omega)}\leq C
d(y)^{m-|\alpha|-n/2}.
\end{equation}

Now we split the discussion into a few cases. 

\vskip 0.08in \noindent {\bf Case I: either $|x-y|\geq N\, d(y)$ or $|x-y|\geq N\,d(x)$ for some $N\geq 25$.}

Let us  first assume that
$|x-y|\geq N d(y)$, $N\geq 25$. In this case we have $\eta\left(\frac{x-y}{d(y)}\right)=0$ and hence,
\begin{equation}\label{eq8.24.1}
\nabla_x^i {\mathcal R}_{\alpha}(x,y)=\nabla_x^i \partial_y^\alpha G(x,y), \quad\mbox{for $x, y\in\Omega$ such that }|x-y|\geq N\, d(y), \quad 0\leq i,|\alpha|\leq \lambda.
\end{equation}

As
before, we denote by $y_0$ some point on the boundary such that
$|y-y_0|=d(y)$. Then by \eqref{eq8.16}, (\ref{eq8.15}) the function $x\mapsto {\mathcal
R}_\alpha(x,y)$ is $m$-harmonic in $\Omega\setminus B_{3d(y)/2}(y_0)$.
Hence, by Proposition~\ref{p7.4} with $r=6d(y)$
\begin{equation}\label{eq8.22}
|\nabla_x^i {\mathcal R}_{\alpha}(x,y)|^2\leq
C\,\frac{d(y)^{2\lambda+n-4m}}{|x-y_0|^{2\lambda+2n-4m+2i}}\,\int_{C_{3d(y)/2,24d(y)}(y_0)}|{\mathcal
R}_\alpha(z,y)|^2\,dz,\quad\mbox{for}\quad 0\leq i\leq \lambda,
\end{equation}

\noindent provided $|x-y_0|\geq 4r=24\,d(y).$ The latter condition is automatically satisfied for $x$ such that $|x-y|\geq N\, d(y)$ with $N\geq 25$, which is the current assumption on $x,y$.
The
right-hand side of (\ref{eq8.22}) is bounded by
\begin{eqnarray}\label{eq8.23}
&& C\,\frac{d(y)^{2\lambda+n-2m}}{|x-y_0|^{2\lambda+2n-4m+2i}}\,\int_{C_{3d(y)/2,24d(y)}(y_0)}\frac{|{\mathcal
R}_\alpha(z,y)|^2}{|z-y_0|^{2m}}\,dz\nonumber\\[4pt]
&&\qquad \leq C\,\frac{d(y)^{2\lambda+n-2m}}{|x-y_0|^{2\lambda+2n-4m+2i}}\,\int_{\Omega}|\nabla_z^m{\mathcal R}_\alpha(z,y)|^2\,dz
\leq C\,\frac{d(y)^{2\lambda-2|\alpha|}}{|x-y|^{2\lambda+2n-4m+2i}},
\end{eqnarray}

\noindent by Hardy's inequality and (\ref{eq8.21}). Therefore, 
\begin{equation}\label{eq8.24}
|\nabla_x^i {\mathcal R}_{\alpha}(x,y)|^2\leq
C\,\frac{d(y)^{2\lambda-2|\alpha|}}{|x-y|^{2\lambda+2n-4m+2i}}, \quad \mbox{when}\quad |x-y|\geq N\, d(y),\quad 0\leq i,|\alpha|\leq \lambda.
\end{equation}

\noindent By \eqref{eq8.24.1}, the estimate \eqref{eq8.24} with $j:=|\alpha|$ implies \eqref{eq8.25}. 
Also, by the symmetry of Green's function we automatically deduce \eqref{eq8.26}.
This furnishes the desired estimates on Green's function when either $|x-y|\geq N\, d(y)$ or $|x-y|\geq N \,d(x)$, both in the case when $n$ is odd and when $n$ is even. 

Turning to the estimates on the regular part of Green's function, we observe that, in  particular, \eqref{eq8.25} and \eqref{eq8.26} combined give the estimate
\begin{equation}\label{eq8.27}
|\nabla_x^i \nabla_y^{j}G(x,y)|^2\leq
\frac{C}{|x-y|^{2n-4m+2i+2j}}, \quad \mbox{when}\quad |x-y|\geq N \min\{d(x),d(y)\},\quad 0\leq i,j\leq \lambda.
\end{equation}

Furthermore, if $n$ is odd, then 
\begin{equation}\label{eq8.28}
|\nabla_x^i\nabla_y^j \Gamma(x-y)|\leq \frac{C}{|x-y|^{n-2m+i+j}}\qquad\mbox{for
all}\qquad x,y\in\Omega, \quad i,j\geq 0,
\end{equation}

\noindent while if $n$ is even, then 
\begin{equation}\label{eq8.29}
|\nabla_x^i\nabla_y^j \Gamma(x-y)|\leq C_1\,|x-y|^{-n+2m-i-j}\log \frac{{\rm diam}\,{(\Omega)}}{|x-y|}+ C_2\, |x-y|^{-n+2m-i-j}\quad\mbox{if} \,\, 0\leq i+j\leq 2m-n,
\end{equation}

\noindent for all $x,y\in\Omega$. 

Combining this  with  (\ref{eq8.27}) we
deduce \eqref{eq8.30} and \eqref{eq8.31} in the cases when the dimension is odd and even, respectively.

\comment{
\todo{Note: The formulas \eqref{eq8.30}, \eqref{eq8.31} are only strong for the highest derivatives, i.e., when $i=j=\lambda$, otherwise the power $-n+2m-i-j$ is positive, and we actually estimate by a bigger quantity. Yet, as far as $S$ goes, nothing stronger can be obtained from the Green function bounds, because of the estimates on the fundamental solution.}}

\vskip 0.08in \noindent {\bf Case II: either $|x-y|\leq N^{-1} d(y)$ or $|x-y|\leq N^{-1} d(x)$ for some $N\geq 25$}. 

Assume that $|x-y|\leq N^{-1}d(y)$. For such $x$ we have
$\eta\bigl(\frac{x-y}{d(y)}\bigr)= 1$ and therefore
\begin{equation}\label{eq8.32}
{\mathcal R}_{\alpha} (x,y)=\partial^{\alpha}_y\,
G(x,y)+P^{\alpha}(x-y)\log \frac{|x-y|}{d(y)}-Q^{\alpha}(x-y).
\end{equation}

\noindent Hence, if $n$ is odd, 
\begin{equation}\label{eq8.33}
{\mathcal R}_{\alpha} (x,y)=\partial^{\alpha}_y\left(
G(x,y)-
\Gamma(x-y)\right),\quad \mbox{when}\quad |x-y|\leq N^{-1}d(y),
\end{equation}
\noindent and if $n$ is even, 
\begin{equation}\label{eq8.34}
{\mathcal R}_{\alpha} (x,y)=\partial^{\alpha}_y\left(
G(x,y)-
\Gamma(x-y)\right)+P^\alpha(x-y)\log \frac{{\rm diam}\,\Omega}{d(y)},\quad \mbox{when}\quad |x-y|\leq N^{-1}d(y).
\end{equation}

\noindent By the interior estimates for solutions of elliptic
equations for any $i\leq m$
\begin{equation}\label{eq8.35}
|\nabla_x^i{\mathcal R}_\alpha(x,y)|^2\leq
\frac{C}{d(y)^{n+2i}}\int_{B_{d(y)/8}(x)}|{\mathcal R}(z,y)|^2\,dz,
\end{equation}

\noindent since by \eqref{eq8.15} for any fixed $y$ the function ${\mathcal R}_\alpha(\cdot, y)$ is $m$-harmonic in
$B_{d(y)/4}(y)\supset B_{d(y)/8}(x)$. Next, since for every $z\in B_{d(y)/8}(x)\subset B_{d(y)/4}(y)$ and any $y_0\in\po$ such that $|y-y_0|=d(y)$ we have $|z-y_0|\approx d(y),$
one can  bound the expression
above by
\begin{equation}\label{eq8.36}
\frac{C}{d(y)^{n+2i-2m}}\int_{B_{d(y)/4}(y)}\frac{|{\mathcal
R}(z,y)|^2}{|z-y_0|^{2m}}\,dz \leq \frac{C}{d(y)^{n+2i-2m}}\left\|\nabla_x^m
{\mathcal R}(\cdot,y)\right\|_{L^2(\Omega)}^2\leq \frac{C}{d(y)^{2n-4m+2i+2|\alpha|}},
\end{equation}

\noindent with $0\leq |\alpha|\leq \lambda$. Then, overall, 
\begin{equation}\label{eq8.36.1}
|\nabla_x^i{\mathcal R}_\alpha(x,y)|^2\leq \frac{C}{d(y)^{2n-4m+2i+2|\alpha|}},\quad |x-y|\leq N^{-1}d(y), \quad 0\leq |\alpha|\leq \lambda.
\end{equation}

Having this at hand, we turn to the estimates on the regular part of Green's function, starting with the case when $n$ is odd. It follows from \eqref{eq8.33} and  \eqref{eq8.36.1} that   
\begin{equation}\label{eq8.37}
|\nabla_x^i\nabla_y^jS(x,y)|\leq
\frac{C}{d(y)^{n-2m+i+j}}, \quad 0\leq i\leq m, \quad 0\leq j\leq \lambda, \quad |x-y|\leq N^{-1}d(y).
\end{equation}

\noindent and hence, by symmetry,
\begin{equation}\label{eq8.38}
|\nabla_x^i\nabla_y^jS(x,y)|\leq
\frac{C}{d(x)^{n-2m+i+j}}, \quad 0\leq i\leq\lambda, \quad 0\leq j\leq m, \quad |x-y|\leq N^{-1}d(x). 
\end{equation}

\noindent However,  we have
\begin{equation}\label{eq8.39}
|x-y|\leq N^{-1}d(y)\quad\Longrightarrow \quad (N-1)\, d(y)\leq N d(x)\leq (N+1)\, d(y),
\end{equation}

\noindent i.e. $d(y)\approx d(x)\approx \min\{d(y), d(x)\}\approx \max\{d(y), d(x)\}$ whenever $|x-y|$ is less than or equal to $N^{-1}d(y)$ or $N^{-1}d(x)$. 
Therefore, when the dimension is odd, we obtain \eqref{eq8.40} and \eqref{eq8.41} depending on the range of $i,j$. 

\comment{
\todo{Once again, the formula \eqref{eq8.40} is the only strong one, it applies  when $i=j=\lambda$. The formula \eqref{eq8.41} is rather weak - estimating by the bigger quantity.}}

Moreover, by \eqref{eq8.28},  for the Green function itself, we then arrive at \eqref{eq8.42}, \eqref{eq8.43}, once again depending on the range of $i,j$.

Similar considerations apply to the case when the dimension is even, leading by \eqref{eq8.34} to the following results: 
\begin{multline}\label{eq8.44}
 |\nabla_x^i\nabla_y^jS(x,y)|\leq  
C\,{d(y)^{2m-n-i-j}}\,\left( C'+\log \frac{{\rm diam}\,\Omega}{d(y)}\right), \\[4pt]
\mbox{for}\, 0\leq i\leq m, \, 0\leq j\leq m-\frac n2, \, |x-y|\leq N^{-1}d(y), 
\end{multline}

\noindent and 
\begin{multline} \label{eq8.45}
|\nabla_x^i\nabla_y^jS(x,y)|\leq 
C\,{d(x)^{2m-n-i-j}}\,\left( C'+\log \frac{{\rm diam}\,\Omega}{d(x)}\right),\\[4pt]
\mbox{for} \, 0\leq i\leq m-\frac n2, \, 0\leq j\leq m, \, |x-y|\leq N^{-1}d(x).
\end{multline}

\noindent Overall, in view of \eqref{eq8.39}, and the fact that $2m-n-i-j\geq 0$ whenever $0\leq i,j \leq m-\frac n2$, \eqref{eq8.44}--\eqref{eq8.45} yield \eqref{eq8.46}.

\comment{
\todo{Once again, above is only strong in the case when $i=j=m-n/2$, when in fact there is no power in front of the $C+\log$ term.}}

As for the Green's function estimates, when $n$ is even, \eqref{eq8.32} and \eqref{eq8.36.1} lead to the bounds 
\begin{eqnarray}\label{eq8.47}
&& |\nabla_x^i\nabla_y^jG(x,y)|\leq  
{C}\,{d(y)^{2m-n-i-j}}\,\left( C'+\log \frac{d(y)}{|x-y|}\right), \, 0\leq i,j\leq m-\frac n2, \, |x-y|\leq N^{-1}d(y), 
\\[4pt] \label{eq8.48}
&& |\nabla_x^i\nabla_y^jG(x,y)|\leq  
{C}\,{d(x)^{2m-n-i-j}}\,\left( C'+\log \frac{d(x)}{|x-y|}\right), \, 0\leq i,j\leq m-\frac n2, \, |x-y|\leq N^{-1}d(x),
\end{eqnarray}

\noindent and therefore, one obtains \eqref{eq8.49}.

Finally, it remains to consider  

\vskip 0.08in \noindent {\bf Case III: $|x-y|\approx d(y)\approx d(x)$},
 or more precisely, the situation when
\begin{equation}\label{eq8.50}
N^{-1}\,d(x)\leq |x-y|\leq Nd(x)\quad \mbox{and}\quad
N^{-1}\,d(y)\leq |x-y|\leq Nd(y)\quad\mbox{for some} \quad N\geq 25.
\end{equation}

Fix any $x,y$ as in \eqref{eq8.50}. Then the mapping $z\mapsto G(z,y)$ is $m$-harmonic in
$B_{d(x)/(2N)}(x)$. Assume first that $n$ is odd. Then by the interior estimates, with $x_0\in\po$
such that $|x-x_0|=d(x)$, we have
\begin{eqnarray}\label{eq8.51}\nonumber
|\nabla_x^i\nabla_y^{|\alpha|} G(x,y)|^2 &\leq &
\frac{C}{d(x)^{n+2i}}\,\int_{B_{d(x)/(2N)}(x)}|\nabla_y^{|\alpha|}
G(z,y)|^2\,dz\\[4pt]
&\leq&  \frac{C}{d(x)^{n+2i}}\,\int_{B_{d(x)/(2N)}(x)}|\nabla_y^{|\alpha|}
\Gamma(z-y)|^2\,dz+
\frac{C}{d(x)^{n+2i-2m}}\,\int_{B_{2d(x)}(x_0)}\frac{|{\mathcal
R}_{\alpha}(z,y)|^2}{|z-x_0|^{2m}}\,dz
\end{eqnarray}

\noindent by definition \eqref{eq8.13}. Then the expression above is bounded by
\begin{eqnarray}\label{eq8.51.1}
&&  \frac{C}{d(x)^{n+2i}}\,\int_{B_{d(x)/(2N)}(x)}|\nabla_y^{|\alpha|}
\Gamma(z-y)|^2\,dz+ \frac{C}{d(x)^{n+2i-2m}}\,\int_{\Omega}|\nabla_z^m {\mathcal
R}_{\alpha}(z,y)|^2\,dz\nonumber \\[4pt]
&& \quad \leq \frac{C}{d(x)^{2n-4m+2i+2|\alpha|}}+\frac{C}{d(x)^{n-2m+2i}d(y)^{n-2m+2|\alpha|}},
\end{eqnarray}

\noindent provided that $0\leq i\leq m$, $0\leq |\alpha|\leq m-\frac n2+\frac 12$. 

The same estimate on derivatives of Green's function holds when $n$ is even, since

\begin{eqnarray}\label{eq8.52}\nonumber
&&|\nabla_x^i\nabla_y^{|\alpha|} G(x,y)|^2\leq 
\frac{C}{d(x)^{n+2i}}\,\int_{B_{d(x)/(2N)}(x)}|\nabla_y^{|\alpha|}
G(z,y)|^2\,dz \leq \frac{C}{d(x)^{n+2i}}\,\int_{B_{d(x)/(2N)}(x)}|P^{\alpha}(z-y)|^2\,dz \\[4pt]
&&\qquad +\frac{C}{d(x)^{n+2i}}\,\int_{B_{d(x)/(2N)}(x)}|Q^{\alpha}(z-y)|^2\,dz+
\frac{C}{d(x)^{n+2i-2m}}\,\int_{B_{2d(x)}(x_0)}\frac{|{\mathcal
R}_{\alpha}(z,y)|^2}{|z-x_0|^{2m}}\,dz,
\end{eqnarray}

\noindent since the absolute value of $\log \frac{|z-y|}{d(y)}$ is bounded by a constant for $z,x,y$ as in \eqref{eq8.52}, \eqref{eq8.50}. This yields \eqref{eq8.53}, both in the case when $n$ is even and odd.

Furthermore, if $n$ is odd, the same argument (or simply an estimate on the difference between the Green function and $S$, that is, the fundamental solution) can be used to deduce
\begin{multline}\label{eq8.54.1}
|\nabla_x^i\nabla_y^jS(x,y)|\leq  \frac{C}{\min \{d(x),d(y),|x-y|\}^{n-2m+i+j}}\approx \frac{C}{\max \{d(x),d(y),|x-y|\}^{n-2m+i+j}},\\[4pt]
\mbox{when $x,y$ satisfy \eqref{eq8.50}},\quad 0\leq i,j\leq m-n/2+1/2,
\end{multline}

\noindent which gives \eqref{eq8.54}.
However, if $n$ is even, using \eqref{eq8.50},  we are led to the bound
\begin{eqnarray}\label{eq8.55.1}\nonumber
|\nabla_x^i\nabla_y^jS(x,y)|&\leq & \frac{C}{\min \{d(x),d(y),|x-y|\}^{n-2m+i+j}}\left(C'+\log \frac{{\rm diam}\,\Omega}{\max \{d(x),d(y),|x-y|\}^{n-2m+i+j}}\right)\\[4pt]
&\approx & \frac{C}{\max \{d(x),d(y),|x-y|\}^{n-2m+i+j}}\left(C'+\log \frac{{\rm diam}\,\Omega}{\max \{d(x),d(y),|x-y|\}^{n-2m+i+j}}\right),
\end{eqnarray}
\noindent for $0\leq i,j\leq m-\frac n2$. This establishes \eqref{eq8.55} and finishes the argument. \ep

\comment{
\todo{This is useful information from a left-over:

 It helps to observe that the cases \eqref{eq8.25} --  \eqref{eq8.26} are disjoint from \eqref{eq8.42} -- \eqref{eq8.43} -- \eqref{eq8.49}, that is, cannot happen simultaneously. The condition $|x-y|\leq N^{-1}\max\{d(x),d(y)\}$ excludes the possibility of both \eqref{eq8.25} and \eqref{eq8.26}. Also, the bound \eqref{eq8.53} is the same as \eqref{eq8.25}, \eqref{eq8.26}, \eqref{eq8.42}, \eqref{eq8.43} for the case when $d(x)$, $d(y)$ and $|x-y|$ are all comparable. Hence, it can be suitably absorbed. Finally, it is straightforward to check that 
\begin{equation}\label{eq8.56}
C'+\log \frac{\min\{d(x),d(y)\}}{|x-y|}\approx  \log \left(1+\frac{\min\{d(x),d(y)\}}{|x-y|}\right)
\end{equation}
\noindent  for all $0\leq i,j\leq \lambda,$ and  $ |x-y|\leq N^{-1}\max\{d(x),d(y)\}.$
}
}

%DOUBLE CHECK THE HARDY INEQUALITY... KUFNER?

\section{Estimates for solutions of the Dirichlet problem} \label{s9} \setcounter{equation}{0}

The results of Section~\ref{s7} provide certain local monotonicity estimates for solutions of the Dirichlet problem. They can, in principle, be translated into the bounds on solution $u$ in terms of data $f$ when $f$ is a suitably supported $C_0^\infty$ function. However, the pointwise control of the Green function that we obtained in Section~\ref{s8}
allows us to estimate the solutions of the Dirichlet problem
for a wide, in some sense, optimal,  class of data. 

\begin{proposition}\label{p9.1} Let $\Omega\subset \RR^n$ be an arbitrary bounded
domain, $m\in\NN$, $n\in [2,2m+1]\cap\NN$, and let $\lambda$ retain the significance of \eqref{eq7.5}. Consider the
boundary value problem
\begin{equation}\label{eq9.1}
(-\Delta)^m u=\sum_{|\alpha|\leq\lambda}c_\alpha \partial^{\alpha}f_{\alpha}, \quad u\in \ring W^{m,2}(\Omega).
\end{equation}

\noindent Then the solution satisfies the following estimates. 

If  $n\in [3,2m+1] \cap\NN$ is odd, then
\begin{eqnarray}\label{eq9.2}|\nabla^{m-\frac n2+\frac 12} u(x)|&\leq & C_{m,n}  \sum_{|\alpha|\leq  m-\frac n2+\frac 12}
\int_\Omega \frac{d(y)^{m-\frac n2+\frac 12-|\alpha|}}{|x-y|}\,\left|f_\alpha(y)\right|\,dy,\qquad  x\in\Omega,
\end{eqnarray}

\noindent whenever the integrals on the right-hand side of \eqref{eq9.2} are finite. In particular, 
\begin{eqnarray}\label{eq9.3}\|\nabla^{m-\frac n2+\frac 12} u\|_{L^\infty(\Omega)}&\leq & C_{m,n,\Omega}    \sum_{|\alpha|\leq  m-\frac n2+\frac 12}
\| d(\cdot)^{m-\frac n2-\frac 12-|\alpha|}f_{\alpha}\|_{L^p(\Omega)},\,\, p>\frac{n}{n-1},
\end{eqnarray}

\noindent provided that the norms on the right-hand side of \eqref{eq9.3} are finite.

If  $n\in [2,2m] \cap\NN$ is even, then 
\begin{eqnarray}\label{eq9.12}|\nabla^{m-\frac n2} u(x)|&\leq & C_{m,n}  \sum_{|\alpha|\leq  m-\frac n2}
\int_\Omega d(y)^{m-\frac n2-|\alpha|}\,\log\left(1+\frac{d(y)}{|x-y|}\right)\,\left|f_\alpha(y)\right|\,dy,\qquad  x\in\Omega,
\end{eqnarray}

\noindent whenever the integrals on the right-hand side of \eqref{eq9.12} are finite. In particular, 
\begin{eqnarray}\label{eq9.13}\|\nabla^{m-\frac n2} u\|_{L^\infty(\Omega)}&\leq & C_{m,n,\Omega}    \sum_{|\alpha|\leq  m-\frac n2}
\| d(\cdot)^{m-\frac n2-|\alpha|}f_{\alpha}\|_{L^p(\Omega)},\,\, p>1,
\end{eqnarray}

\noindent provided that the norms on the right-hand side of \eqref{eq9.13} are finite.

The constants $C_{m,n}$ above depend on $m$ and $n$ only, while the constants denoted by $C_{m,n,\Omega}$ depend on $m$, $n$, and the diameter of the domain $\Omega$.
\end{proposition}

\comment{
\todo{Understand where is $f$ and in which sense we understand the solution. It might be not $u\in\ring W^{m,2}(\Omega)$!!!!}}

\bp
First of all, the integral representation formula
\begin{equation}\label{eq9.4}
u(x)=\int_{\Omega} G(x,y) \sum_{|\alpha|\leq \lambda}c_\alpha \partial^{\alpha}f_{\alpha}(y)
\,dy,\qquad x\in\Omega,
\end{equation}

\noindent follows directly from the definition of Green's function. It implies that
\begin{eqnarray}\label{eq9.5}
\nabla^{\lambda} u(x)&=& \sum_{|\alpha|\leq \lambda}c_\alpha (-1)^{|\alpha|} \int_{\Omega}\nabla_x^{\lambda}\partial^\alpha_y  G(x,y)\,f_\alpha(y)\,dy\qquad x\in\Omega.\end{eqnarray}

Let us now assume that $n\in [3,2m+1] \cap\NN$ is odd. We claim that in this case
\begin{eqnarray}\label{eq9.5.1}
\left|\nabla_x^{m-\frac n2+\frac 12}\partial^\alpha_y  G(x,y)\right|\leq  C\,\frac{d(y)^{m-\frac n2+\frac 12-|\alpha|}}{|x-y|},\qquad  \mbox{ for all\,\,}x,y\in\Omega, \quad 0\leq |\alpha|\leq m-\frac n2+\frac 12.\end{eqnarray}

Indeed, when $|\alpha|=m-\frac n2+\frac 12$, \eqref{eq9.5.1} is equivalent to \eqref{eq8.5.1}. Thus, it remains consider $|\alpha|\leq m-\frac n2-\frac 12$. Now we split the cases, essentially according to Theorem~\ref{t8.1}.  When $|x-y|\geq N d(y),$ for some $N\geq 25$,  \eqref{eq8.25} directly yields \eqref{eq9.5.1}. If $x,y\in \Omega$ are such that $|x-y|\leq N^{-1} d(y)$ then \eqref{eq9.5.1}  follows from \eqref{eq8.43}. Finally, if $N^{-1}d(y) \leq |x-y|\leq Nd(y)$, we observe that 
$d(x)\leq |x-y|+d(y)\leq (1+N)|x-y|$. The latter estimate shows that  \eqref{eq8.26}--\eqref{eq8.53} imply the bound on the left-hand side of \eqref{eq9.5.1} by $C\,|x-y|^{m-\frac n2-\frac 12-|\alpha|}$, and hence, \eqref{eq9.5.1}, since here $d(y)\approx |x-y|$.

Now \eqref{eq9.5.1}
allows us to deduce \eqref{eq9.2}, and \eqref{eq9.3} follows from it via the mapping properties of the Riesz potential. 

Let us now turn to the case when $n\in [2,2m] \cap\NN$ is even. Then we have to show that 
\begin{eqnarray}\label{eq9.5.2}
\left|\nabla_x^{m-\frac n2}\partial^\alpha_y  G(x,y)\right|\leq  C\,d(y)^{m-\frac n2-|\alpha|}\,\log\left(1+\frac{d(y)}{|x-y|}\right),\,\,  \mbox{for all\,\,}x,y\in\Omega, \,\, 0\leq |\alpha|\leq m-\frac n2.\end{eqnarray}

Once again, we split to cases according to Theorem~\ref{t8.1}. When $|x-y|\geq N d(y),$ for some $N\geq 25$,  then $1+\frac{d(y)}{|x-y|}$ is bounded from below and above by a positive constant, so that \eqref{eq8.25} is the same as \eqref{eq9.5.2}. If $x,y\in \Omega$ are such that $|x-y|\leq N^{-1} d(y)$ then \eqref{eq9.5.2}  follows directly from \eqref{eq8.49}. Finally, if $N^{-1}d(y) \leq |x-y|\leq Nd(y)$, then, as before,  $d(x)\leq  (1+N)|x-y|$ and $1+\frac{d(y)}{|x-y|}$ is bounded from below and above by a positive constant. Then  \eqref{eq8.26}, \eqref{eq8.49}, and \eqref{eq8.53} imply the bound on the left-hand side of \eqref{eq9.5.2} by $C\,|x-y|^{m-\frac n2-|\alpha|}$, and hence, \eqref{eq9.5.2}, using that $d(y)\approx |x-y|$. This finishes the proof of \eqref{eq9.12}.

Finally, we obtain from \eqref{eq9.12} the bound 
\begin{eqnarray}\label{eq9.17}|\nabla^{m-\frac n2} u(x)|&\leq & C  \sum_{|\alpha|\leq  m-\frac n2}
\int_\Omega d(y)^{m-\frac n2-|\alpha|}\,\left(\frac{d(y)}{|x-y|}\right)^\eps |f_\alpha(y)|\,dy,\qquad  x\in\Omega.\quad \eps>0,
\end{eqnarray}

\noindent Then, by the mapping properties of the Riesz potential we recover an estimate 
\begin{eqnarray}\label{eq9.18}\|\nabla^{m-\frac n2} u\|_{L^\infty(\Omega)}&\leq & C_\Omega  \sum_{|\alpha|\leq  m-\frac n2}
\left\| d(y)^{m-\frac n2-|\alpha|+\eps} f_\alpha\right\|_{L^p(\Omega)},\quad  \mbox{provided that}\,\, \eps>0,\,\, p>\frac{n}{n-\eps},
\end{eqnarray}

\noindent which, in turn, leads to \eqref{eq9.13}.\ep

\vskip 0.08in \noindent --------------------------------------
\vskip 0.10in

\noindent {\it Svitlana Mayboroda}

\noindent School of Mathematics, University of Minnesota,\\
127 Vincent Hall, 206 Church St SE, Minneapolis, MN, 55408, USA

\noindent {\tt svitlana\@@math.umn.edu}

\vskip 0.10in

\noindent {\it Vladimir Maz'ya}

\noindent Department of Mathematical Sciences, M\&O Building,\\
University of Liverpool, Liverpool L69 3BX, UK

\vskip 0.05in

\noindent Department of Mathematics, Link\"oping University,\\
SE-581 83 Link\"oping, Sweden

\noindent {\tt vlmaz\@@mai.liu.se}

\end{document}